\documentclass{amsart}
\usepackage{graphicx}
\usepackage[all]{xy}
\usepackage{latexsym}
\usepackage{amssymb}
\usepackage{amsmath}
\usepackage{color}

\xyoption{arc}

 \newfam\cyrfam

  \font\tencyr=wncyr10

  \font\sevencyr=wncyr7

  \font\fivecyr=wncyr5

  \textfont\cyrfam=\tencyr \scriptfont\cyrfam=\sevencyr

    \scriptscriptfont\cyrfam=\fivecyr

  \newfam\cyifam

  \font\tencyi=wncyi10

  \font\sevencyi=wncyi7

  \font\fivecyi=wncyi5

  \textfont\cyifam=\tencyi \scriptfont\cyifam=\sevencyi

    \scriptscriptfont\cyifam=\fivecyi

\def\id{{\mbox{1 \hskip -7pt 1}}}
\newcommand{\sgn}{{\mathit s  \mathit g\mathit  n}}
 \newcommand{\lon}{\longrightarrow}
 \newcommand{\bu}{\bullet}
 
 \newcommand{\rar}{\rightarrow}
 \newcommand{\hook}{\hookrightarrow}

\newcommand{\p}{{\partial}}
\newcommand{\Id}{{\mathrm{Id}}}

\newcommand{\Der}{\mathrm{Der}}

\newcommand{\tw}{\mathsf{tw}}
\newcommand{\RGra}{\mathcal{R} \mathcal{G} ra}
\newcommand{\GRav}{\mathcal{G} \mathcal{R}av}
\newcommand{\Grav}{\mathcal{G} rav}

 \newcommand{\Z}{{\mathbb Z}}
 \newcommand{\bS}{{\mathbb S}}
 \renewcommand{\P}{{\mathbb P}}
 \newcommand{\C}{{\mathbb C}}
 \newcommand{\R}{{\mathbb R}}
 
 \newcommand{\K}{{\mathbb K}}

 \newcommand{\ot}{\otimes}

\newcommand{\sd}{{\mathsf d}}

\newcommand{\Def}{\mathsf{Def}}
\newcommand{\fGC}{\mathsf{fGC}}
\newcommand{\GC}{\mathsf{GC}}
\newcommand{\OGC}{\mathsf{OGC}}
\newcommand{\dGC}{\mathsf{dGC}}
\newcommand{\dfGC}{\mathsf{dfGC}}
\newcommand{\grt}{\fg\fr\ft}
 \newcommand{\Beq}{\begin{equation}}
 \newcommand{\Eeq}{\end{equation}}
 \newcommand{\Beqr}{\begin{eqnarray}}
 \newcommand{\Eeqr}{\end{eqnarray}}
 \newcommand{\Beqrn}{\begin{eqnarray*}}
 \newcommand{\Eeqrn}{\end{eqnarray*}}
 \newcommand{\Ba}{\begin{array}}
 \newcommand{\Ea}{\end{array}}
 \newcommand{\Bi}{\begin{itemize}}
 \newcommand{\Ei}{\end{itemize}}
 \newcommand{\Bc}{\begin{center}}
 \newcommand{\Ec}{\end{center}}

 \newcommand{\fg}{{\mathfrak g}}

\newcommand{\fr}{{\mathfrak r}}

\newcommand{\ft}{{\mathfrak t}}

 \newcommand{\cC}{{\mathcal C}}
 \newcommand{\caD}{{\mathcal D}}
 \newcommand{\cE}{{\mathcal E}}
 \newcommand{\cF}{{\mathcal F}}
 \newcommand{\cG}{{\mathcal G}}

 \newcommand{\caL}{{\mathcal L}}
 \newcommand{\cM}{{\mathcal M}}
 \newcommand{\cN}{{\mathcal N}}
 \newcommand{\cP}{{\mathcal P}}
 
 \newcommand{\cR}{{\mathcal R}}
 
 \newcommand{\cT}{{\mathcal T}}

 \newcommand{\ga}{\gamma}
 
 \newcommand{\Ga}{\Gamma}

 \newcommand{\var}{\varepsilon}
 
 \newcommand{\om}{\omega}

 \newcommand{\Ker}{{\mathsf K \mathsf e \mathsf r}\, }
 \newcommand{\Img}{{\mathsf I\mathsf m}\, }

 \newcommand{\sip}{\smallskip}
 \newcommand{\bip}{\bigskip}
 \newcommand{\mip}{\vspace{2.5mm}}

\newcommand{\LB}{\mathcal{L}\mathit{ieb}}

\newcommand{\LBcd}{\mathcal{L}\mathit{ieb}_{c,d}}

\newcommand{\HoLBcd}{\mathcal{H}\mathit{olieb}_{c,d}}

\newcommand{\wHoLBcd}{\widehat{\mathcal{H}\mathit{olieb}}_{c,d}}

\newcommand{\HoLB}{\mathcal{H}\mathit{olieb}}

\newcommand{\HoqLBcd}{\mathcal{H}\mathit{oqlieb}_{c,d}}

\newcommand{\wHoqLBcd}{\widehat{\mathcal{H}\mathit{oqlieb}}_{c,d}}

\newcommand{\RGC}{\mathsf{RGC}}

\newcommand{\Lie}{\mathcal{L} \mathit{ie}}

\theoremstyle{plain}
\swapnumbers

\newtheorem{prop-def}[theorem]{Proposition-definition}

\newtheorem{f-theorem}{Formality Theorem}[section]
\newtheorem{main-theorem}{Main~Theorem}[section]
\newtheorem{section-theorem}{Theorem}[section]

\theoremstyle{definition}

\usepackage{tikz}
\usetikzlibrary{matrix}
\usetikzlibrary{shapes}
\usetikzlibrary{arrows}
\usetikzlibrary{calc,3d}
\usetikzlibrary{decorations,decorations.pathmorphing,decorations.pathreplacing,decorations.markings}
\usetikzlibrary{through}

\tikzset{ext/.style={circle, draw,inner sep=1pt},int/.style={circle,draw,fill,inner sep=1.4pt},nil/.style={inner sep=1pt}}
\tikzset{cy/.style={circle,draw,fill,inner sep=2pt},scy/.style={circle,draw,inner sep=2pt},scyx/.style={draw,cross out,inner sep=2pt},scyt/.style={draw,regular polygon,regular polygon sides=3,inner sep=0.95pt}}
\tikzset{exte/.style={circle, draw,inner sep=3pt},inte/.style={circle,draw,fill,inner sep=3pt}}
\tikzset{diagram/.style={matrix of math nodes, row sep=3em, column sep=2.5em, text height=1.5ex, text depth=0.25ex}}
\tikzset{diagram2/.style={matrix of math nodes, row sep=0.5em, column sep=0.5em, text height=1.5ex, text depth=0.25ex}}

 \tikzset{
  rightblack/.style={
    decoration={markings,mark=at position .8 with {\arrow[scale=1.2,black]{latex}}},
    postaction={decorate},
    shorten >=0.4pt}}
\tikzset{
  leftblack/.style={
    decoration={markings,mark=at position .55 with {\arrowreversed[scale=1.2,black]{latex}}},
    postaction={decorate},
    shorten >=0.4pt}}

\begin{document}

 \sloppy

 \newenvironment{proo}{\begin{trivlist} \item{\sc {Proof.}}}
  {\hfill $\square$ \end{trivlist}}

\long\def\symbolfootnote[#1]#2{\begingroup%
\def\thefootnote{\fnsymbol{footnote}}\footnote[#1]{#2}\endgroup}

\title{Gravity properad and moduli spaces  $\cM_{g,n}$}

\author{Sergei~A.\ Merkulov}
\address{Sergei~A.\ Merkulov:  Department of Mathematics, Luxembourg University,  Grand Duchy of Luxembourg }
\email{sergei.merkulov@uni.lu}

\begin{abstract}
Let $\cM_{g,m+n}$ be the moduli space of algebraic curves of genus $g$ with $m\geq 1$ boundaries and   $n\geq 0$ marked points, and $H_c^{\bu}(\cM_{m+n})$  its compactly supported cohomology group.  We prove that the collection of  $\bS_m^{op}\times \bS_n$-modules
$$
\left\{\prod_{2g+(m+n)\geq 3}H_c^{\bu}\left(\cM_{g,m+n}\times \R_+^m\right)
 \right\}=:\GRav,
$$
has the structure of a {\it properad} (called the {\it gravity properad})
such that  it  contains  the E.\ Getzler's {\it gravity operad} as  the sub-collection  $\{H_c^{\bu-1}(\cM_{0,1+n})\}_{n\geq 2}$.
The properadic structure in $\GRav$ is highly non-trivial and generates higher genus
cohomology classes from lower ones (which is demonstrated on infinitely many non-trivial examples producing higher genus cohomology classes from just zero genus ones).
Moreover, we prove that the generators of the 1-dimensional cohomology groups
$H_c^{\bu-1}(\cM_{0,1+2})$, $H_c^{\bu-2}(\cM_{0,2+1})\ot \sgn_2$ and  $H_c^{\bu-3}(\cM_{0,3+0})\ot \sgn_3$
satisfy with respect to this properadic structure the relations of the (degree shifted) quasi-Lie bialgebra, a fact making the totality of cohomology groups
$$
\prod_{g,m,n} H_c^{\bu}(\cM_{g,m+n})\ot_{\bS_m^{op}\times \bS_n} (\sgn_m\ot \id_n)
$$
into a complex with the differential fully determined by the just mentioned three cohomology classes. It is proven that this complex contains infinitely many non-trivial cohomology classes, all coming from  M.\ Kontsevich's odd graph complex.

\sip

The prop structure in $\GRav$ is established with the help of T.\ Willwacher's twisting endofunctor $\tw$ (in the category of properads under the operad of Lie algebras) and K. Costello's theory of moduli spaces of nodal disks with marked boundaries and internal marked points.

\sip

\sip
\noindent {\sc Mathematics Subject Classifications} (2000). 14D22,  18M85, 18G85.

\noindent {\sc Key words}. Moduli spaces of algebraic curves, graph complexes, properads, Lie bialgebras.
\end{abstract}

 \maketitle
 \markboth{}{}

{\small
{\small
\tableofcontents
}
}

{\Large
\section{\bf Introduction}
}

Let $\cM_{g,m+n}$ be the moduli space of algebraic curves of genus $g$ with $m\geq 1$ boundaries and  $n\geq 0$ marked points\footnote{$\cM_{g,m+n}$ is best understood in the context of this paper as the moduli spaces of complete hyperbolic surfaces with $m\geq 1$ geodesic boundaries (of any fixed length) and $n$ cusps.}, and
$$
H_c^{\bu}(\cM_{g,m+n})\simeq
H_{\bu+6-6g-2(m+n)}(\cM_{g,m+n})
$$
 its compactly supported cohomology groups. Let $\R_+\subset \R$ be the subspace of positive real numbers.  We prove that the collection of  $\bS_m^{op}\times \bS_n$-modules (with $\bS_m$ acting diagonally on the (tensor) products),
\Beq\label{1 GRav as family of H(M_g,n)}
\left\{\prod_{g\geq 0\atop 2g+(m+n)\geq 3}H_c^{\bu}(\cM_{g,m+n}\times \R_+^m)\simeq H_c^{\bu-m}(\cM_{g,m+n})\ot \sgn_m 
 \right\}=:\GRav,
\Eeq
 has the structure of a {\it properad} --- called the {\it gravity properad} --- such that
\Bi
\item[(i)] $\GRav$  contains as a sub-properad the (degree shifted) Ezra Getzler's {\it gravity operad}\, \cite{Ge} as the sub-collection of cohomology groups of moduli spaces of genus zero algebraic curves,
    $$
    \{H_c^{\bu-1}(\cM_{0,1+n})\}_{n\geq 2} \subset \GRav,
    $$

\item[(ii)] $\GRav$  comes equipped with a
morphism from the properad of (degree shifted) Lie quasi-bialgebras\footnote{The properad controlling degree shifted quasi-Lie bialgebras with  Lie bracket of degree $1-d$ and the quasi-Lie cobracket of degree $1-c$  is denoted by $q\LBcd$.}
\Beq\label{1: map qLB to Grav}
\Ba{rccc}
j: & q\LB_{-1,0} & \lon & \GRav \vspace{1mm}\\
&
\Ba{c}\begin{xy}
 <0mm,0.66mm>*{};<0mm,3mm>*{}**@{-},
 <0.39mm,-0.39mm>*{};<2.2mm,-2.2mm>*{}**@{-},
 <-0.35mm,-0.35mm>*{};<-2.2mm,-2.2mm>*{}**@{-},
 <0mm,0mm>*{\bu};<0mm,0mm>*{}**@{},
   <0.39mm,-0.39mm>*{};<2.9mm,-4mm>*{^{_2}}**@{},
   <-0.35mm,-0.35mm>*{};<-2.8mm,-4mm>*{^{_1}}**@{},
\end{xy}\Ea
 &\lon &
1\in H_c^{\bu}(\cM_{0,1+2}\times \R_+)\simeq \ \id_2[-1] \subset \GRav(1,2)\vspace{1mm} \\
&
\Ba{c}\begin{xy}
 <0mm,-0.55mm>*{};<0mm,-2.5mm>*{}**@{-},
 <0.5mm,0.5mm>*{};<2.2mm,2.2mm>*{}**@{-},
 <-0.48mm,0.48mm>*{};<-2.2mm,2.2mm>*{}**@{-},
 <0mm,0mm>*{\bu};<0mm,0mm>*{}**@{},
 <0.5mm,0.5mm>*{};<2.7mm,2.8mm>*{^{_2}}**@{},
 <-0.48mm,0.48mm>*{};<-2.7mm,2.8mm>*{^{_1}}**@{},
 \end{xy}\Ea
 & \lon &
1\in H_c^{\bu-2}(\cM_{0,2+1}\times \R_+^2)\simeq \sgn_2[-2]\subset \GRav(2,1)\vspace{2mm} \\
&
\Ba{c}\begin{xy}
 <0mm,-1mm>*{\bu};<-4mm,3mm>*{^{_1}}**@{-},
 <0mm,-1mm>*{\bu};<0mm,3mm>*{^{_2}}**@{-},
 <0mm,-1mm>*{\bu};<4mm,3mm>*{^{_3}}**@{-},
 \end{xy}\Ea
 & \lon &
1\in H_c^{\bu-3}(\cM_{0,3+0}\times \R_+^3)\simeq \sgn_3[-3]\subset \GRav(3,0)\\
\Ea
\Eeq
The images on the r.h.s.\  should be understood as the unique cohomology classes of the moduli spaces of  hyperbolic spheres $S^2$ with (i) one geodesic boundary and two cusps, (ii) two geodesic boundaries and one cusp, (iii) and three geodesic boundaries, respectively.
\Ei

\sip

The existence of a properad structure on the family of cohomology groups (\ref{1 GRav as family of H(M_g,n)}) is  established  with the help of
(a) Thomas Willwacher's twisting endofunctor $\tw$   in the category of properads under the operad of Lie algebras \cite{W} which is applied in this paper to
 the properad of ribbon graphs $\cR\cG ra_d$, $d\in \Z$, introduced earlier in \cite{MW} by Thomas Willwacher and the author,
 $$
 \tw: \RGra_d \lon \tw\RGra_d=\{\tw\RGra_d(m,n)\}_{m\geq 1, n\geq 0},
 $$
 and
(b) Kevin Costello's theory of moduli spaces $\{\caD_{g,m,0,n}\}_{m\geq 1,n\geq 0}$ of nodal disks  with $m$  marked boundaries and $n$ internal marked points (such that each disk contains at most one internal marked point) \cite{Co1,Co2}.

\sip

Every $\bS_m^{op}\times \bS_n$-module $\tw\RGra_d(m,n)$
of the dg properad $\tw  \cR\cG ra_d=\{\tw\RGra_d(m,n)\}$ is generated by oriented ribbon graphs with $m$ labelled boundaries, and with vertices of two types --- $n$ labelled white vertices and any number of black (unlabelled) vertices which are assigned the cohomological degree $d$. The differential and the properadic structure in $\tw  \cR\cG ra_d$ are of purely combinatorial nature and can be given by  explicit formulae (see \S 3). Though {\it orientations}\, of ribbon graphs depend much on the parity of $d$ we show that, up to a degree shift, all the twisted dg prop(erad)s $\tw  \cR\cG ra_d$ for various $d$ are isomorphic to each other so that it is enough to study the case $d=0$.\footnote{It is worth however keeping in mind that the integer parameter $d$ can be arbitrary as the dg properad $\tw  \cR\cG ra_d$ admits a canonical representation \cite{Me} in the cyclic  Hochschild complex of an arbitrary Poincare duality algebra of degree $d$; in particular, the (degree shifted) gravity properad $\GRav$ acts on the reduced equivariant homology $\bar{H}_\bu^{S^1}(LM)$ of the free loop space $LM$ of any closed semisimple $d$-dimensional manifold $M$ (cf.\ \cite{CFL}).}  The dg properad $\tw  \cR\cG ra_0$
contains a dg sub-properad $\tw  \cR\cG ra_0^{\geq 3}$ spanned by ribbon graphs with all {\it black}\, vertices at least trivalent. The inclusion
$$
\tw  \cR\cG ra_0^{\geq 3} \lon \tw  \cR\cG ra_0
$$
is a quasi-isomorphism modulo the family of polytope-like cohomology classes (for a precise formulation see Proposition {\ref{3: prop on tw^3Rgra to twRgra}}). The main problem at this point is to compute the associated  cohomology properad $H^\bu(\tw  \cR\cG ra_0^{\geq 3})$ which is solved in \S 3 with the help of the remarkable Kevin Costello's
homotopy equivalence \cite{Co1,Co2} of moduli spaces.

\mip

{\sc{Theorem A}}. $H^\bu(\tw  \cR\cG ra_0^{\geq 3})= \GRav$.

\mip

This result implies the above claim about the canonical properadic structure on the above family of cohomology groups $H_c^{\bu}(\cM_{g,m+n}\times \R_+^m)$. Its proof is relatively short and is based on the observation (see \S 3 for full details)  that, as a family of complexes, the twisted properad $\tw\cR\cG ra_0^{\geq 3}=\{\tw\cR\cG ra_0^{\geq 3}(m,n)\}_{m\geq 1, n\geq 0}$  can be identified  with the cell complexes of the family of moduli spaces $\{\caD_{g,m,0,n}\}_{m\geq 1,n\geq 0}$ described explicitly in \cite{Co2}.
This observation gives us an explicit properadic structure on the $\bS$-bimodule
$$
\{\cC hains(\caD_{g,m,0,n})\}_{m\geq 1,n\geq 0, 2g+m+n\geq 3}
$$
which is called the {\em chain gravity properad}. It contains as a suboperad the cell complexes
of the topological $\bS$-bimodule $\{\caD_{0,1,0,n}\}$ whose cohomology can in turn be identified, using the results of \cite{Wa}, with the gravity operad  $\{H_{\bu-1}(\cM_{0,1+n})\}_{n\geq 2}$ first introduced in \cite{Ge}.

\sip

 According to the results \cite{Co1,Co2} by Kevin Costello  and the uniformization theorem, the moduli spaces  $\caD_{g,m,0,n}$ are homotopy equivalent to $\cM_{g,m+n}$ so that the above Theorem A follows.
The claim about the map (\ref{1: map qLB to Grav}) is proven in \S 3 by direct calculations using the complex of twisted ribbon graphs. A very useful fact is that the map (\ref{1: map qLB to Grav}) can be described explicitly in terms of the cocycles representing the cohomology classes
on the r.h.s.\ as follows,
$$
\Ba{rccc}
j: & q\LB_{-1,0} & \lon & \text{co-cycles in}\ \tw\RGra_0^{\geq 3} \vspace{1mm}\\
&
\Ba{c}\begin{xy}
 <0mm,0.66mm>*{};<0mm,4mm>*{^{_{\bar{1}}}}**@{-},
 <0.39mm,-0.39mm>*{};<2.2mm,-2.2mm>*{}**@{-},
 <-0.35mm,-0.35mm>*{};<-2.2mm,-2.2mm>*{}**@{-},
 <0mm,0mm>*{\bu};<0mm,0mm>*{}**@{},
   <0.39mm,-0.39mm>*{};<2.9mm,-4mm>*{^{_2}}**@{},
   <-0.35mm,-0.35mm>*{};<-2.8mm,-4mm>*{^{_1}}**@{},
\end{xy}\Ea
 &\lon &
  \Ba{c}\resizebox{10mm}{!}{  \xy
 (3.5,4)*{^{\bar{1}}};
 (7,0)*+{_2}*\frm{o}="A";
 (0,0)*+{_1}*\frm{o}="B";
 \ar @{-} "A";"B" <0pt>
\endxy} \Ea  \in \tw\RGra_0^{\geq 3}(1,2)\vspace{1mm} \\
&
\Ba{c}\begin{xy}
 <0mm,-0.55mm>*{};<0mm,-2.5mm>*{}**@{-},
 <0.5mm,0.5mm>*{};<2.2mm,2.2mm>*{}**@{-},
 <-0.48mm,0.48mm>*{};<-2.2mm,2.2mm>*{}**@{-},
 <0mm,0mm>*{\bu};<0mm,0mm>*{}**@{},
 <0mm,-0.55mm>*{};<0mm,-3.8mm>*{_1}**@{},
 <0.5mm,0.5mm>*{};<2.7mm,2.8mm>*{^{_{\bar{2}}}}**@{},
 <-0.48mm,0.48mm>*{};<-2.7mm,2.8mm>*{^{_{\bar{1}}}}**@{},
 \end{xy}\Ea
 & \lon &
 \frac{1}{2}
 \left(
\Ba{c}\resizebox{6mm}{!}{
\mbox{$\xy
 (0.5,0.9)*{^{{^{\bar{1}}}}},
(0.5,5)*{^{{^{\bar{2}}}}},
 (0,-8)*+{_{_1}}*\frm{o}="C";
(0,-2)*{\bu}="A";
(0,-2)*{\bu}="B";
"A"; "B" **\crv{(6,6) & (-6,6)};
 \ar @{-} "A";"C" <0pt>
\endxy$}}
\Ea
-
\Ba{c}\resizebox{6mm}{!}{
\mbox{$\xy
 (0.5,0.9)*{^{{^{\bar{2}}}}},
(0.5,5)*{^{{^{\bar{1}}}}},
 (0,-8)*+{_{_1}}*\frm{o}="C";
(0,-2)*{\bu}="A";
(0,-2)*{\bu}="B";
"A"; "B" **\crv{(6,6) & (-6,6)};
 \ar @{-} "A";"C" <0pt>
\endxy$}}
\Ea
\right) \in \tw\RGra_0^{\geq 3}(2,1)\vspace{2mm} \\
&
\Ba{c}\begin{xy}
 <0mm,-1mm>*{\bu};<-4mm,3mm>*{^{_{\bar{1}}}}**@{-},
 <0mm,-1mm>*{\bu};<0mm,3mm>*{^{_{\bar{2}}}}**@{-},
 <0mm,-1mm>*{\bu};<4mm,3mm>*{^{_{\bar{3}}}}**@{-},
 \end{xy}\Ea
 & \lon &  \frac{1}{2}\left(-
 \Ba{c}\resizebox{11mm}{!}{
\xy
(2.0,-3.5)*{^{{^{\bar{2}}}}},
 (-2,-3.5)*{^{{^{\bar{1}}}}},
(0.5,4)*{^{{^{\bar{3}}}}},
 (0,-8)*{\bu}="C";
(0,3)*{\bu}="A1";
(0,3)*{\bu}="A2";
"C"; "A1" **\crv{(-5,-9) & (-5,4)};
"C"; "A2" **\crv{(5,-9) & (5,4)};
 \ar @{-} "A1";"C" <0pt>
\endxy}
\Ea
+
\Ba{c}\resizebox{11mm}{!}{
\xy
(2.0,-3.5)*{^{{^{\bar{1}}}}},
 (-2,-3.5)*{^{{^{\bar{2}}}}},
(0.5,4)*{^{{^{\bar{3}}}}},
 (0,-8)*{\bu}="C";
(0,3)*{\bu}="A1";
(0,3)*{\bu}="A2";
"C"; "A1" **\crv{(-5,-9) & (-5,4)};
"C"; "A2" **\crv{(5,-9) & (5,4)};
 \ar @{-} "A1";"C" <0pt>
\endxy}
\Ea\right)
\in \tw\RGra_0^{\geq 3}(3,0)\\
\Ea
$$

Note that the morphism of properads
$$
j: \LB_{-1,0} \rar \GRav
$$
can {\it not}\, be lifted to a morphism from $\LB_{-1,0}$ into the
chain gravity properad --- it holds true only at the cohomology level.
It is not clear at the moment if this morphism is an injection.
 However we can be sure that it is highly non-trivial

\sip

{\sc{Theorem}} B. {\it The morphism $j: \LB_{-1,0} \rar \GRav$  is non-trivial on infinitely many elements of $\LB_{-1,0}$ with loop number $1$ producing thereby via gravity properadic compositions from the unique class in $H^0(\cM_{0,3})$  infinitely many  cohomology  classes in $\prod_{m\geq 1, n\geq 0} H_c^\bu(\cM_{1,m+n})$}.

\sip

The deformation complex of the above morphism $j$ has, according to the general theory \cite{MV} (see also \cite{MW2} for its small but important ``plus" adjustment), the following structure
\Beq\label{1: totality of Mg,n+m as Def}
\Def(q\LB_{-1,0} \stackrel{j}{\rar} \GRav)\simeq \prod_{g,n\geq 0,m\geq 1\atop 2g+n+m\geq 3} H^{\bu-1}(\cM_{g,m+n})\ot_{\bS_m^{op}\times \bS_n} (\sgn_m\ot \id_n)
\Eeq
The induced differential has three parts corresponding to the above mentioned three incarnations of the pair of pants.

\sip

{\sc{Theorem}} C. {\it The complex (\ref{1: totality of Mg,n+m as Def}) has infinitely many cohomology classes.}

\mip

Moreover, these particular cohomology classes can described explicitly in terms of twisted ribbon graphs (see \S 4 below). They all have genus $g=1$ and there are indications that they all come from $H_c^\bu(\cM_1)$;  we conjecture that {\it the cohomology of the complex (\ref{1: totality of Mg,n+m as Def}) is equal to $\prod_{g\geq 1} H_c^\bu(\cM_g)$},
where $\cM_g$ is the moduli space of (unpunctured) algebraic curves of genus $g$ (cf.\ \cite{AWZ}).



\subsection{Notation}
 The set $\{1,2, \ldots, n\}$ is abbreviated to $[n]$;  its group of automorphisms is
denoted by $\bS_n$;
the trivial one-dimensional representation of
 $\bS_n$ is denoted by $\id_n$, while its one dimensional sign representation is
 denoted by $\sgn_n$. The cardinality of a finite set $I$ is denoted by $\# I$.

\sip

We work throughout in the category of differential $\Z$-graded vector spaces over a field $\K$
of characteristic zero; all our differentials have degree +1 (in particular,  the chain complexes of topological spaces and their homology groups are {\it non-positively}\, graded).
 If $V=\oplus_{i\in \Z} V^i$ is a graded vector space, then
$V[k]$ stands for the graded vector space with $V[k]^i:=V^{i+k}$; for $v\in V^i$ we set $|v|:=i$.
For a properad $\cP$ we denote by $\cP\{k\}$ the unique properad which has the following property:
for any graded vector space $V$ there is a one-to-one correspondence between representations of
$\cP\{k\}$ in $V$ and representations of
$\cP$ in $V[-k]$; in particular, $\cE nd_V\{k\}=\cE nd_{V[k]}$.

\mip

{\bf Acknowledgement}.
It is a great pleasure to thank Alexey Kalugin,  Anton Khoroshkin, Sergey Shadrin and especially Thomas Willwacher for valuable communications and discussions.

\bip

\bip

{\Large
\section{\bf On twisting of properads under $\Lie_d$ (after Thomas Willwacher \cite{W})}
}

\mip

\subsection{Operad of degree shifted Lie algebras} For any $d\in \Z$ the operad of degree  shifted Lie algebras is defined as a quotient,
$$
\Lie_{d}:=\cF ree\langle E\rangle/\langle\cR\rangle,
$$
of the free operad generated by an  $\bS$-module $E=\{E(n)\}_{n\geq 2}$ with
 all $E(n)=0$ except\footnote{When representing elements of operads and props
   as (decorated) graphs we tacitly assume that all edges and legs are {\em directed}\, along the flow going from the bottom of the graph to the top.}
$$
E(2):= \sgn_2^{d}\ot \id_1[d-1]=\mbox{span}\left\langle
\Ba{c}\begin{xy}
 <0mm,0.66mm>*{};<0mm,3mm>*{}**@{-},
 <0.39mm,-0.39mm>*{};<2.2mm,-2.2mm>*{}**@{-},
 <-0.35mm,-0.35mm>*{};<-2.2mm,-2.2mm>*{}**@{-},
 <0mm,0mm>*{\bu};<0mm,0mm>*{}**@{},
   <0.39mm,-0.39mm>*{};<2.9mm,-4mm>*{^{_2}}**@{},
   <-0.35mm,-0.35mm>*{};<-2.8mm,-4mm>*{^{_1}}**@{},
\end{xy}\Ea
=(-1)^{d}
\Ba{c}\begin{xy}
 <0mm,0.66mm>*{};<0mm,3mm>*{}**@{-},
 <0.39mm,-0.39mm>*{};<2.2mm,-2.2mm>*{}**@{-},
 <-0.35mm,-0.35mm>*{};<-2.2mm,-2.2mm>*{}**@{-},
 <0mm,0mm>*{\bu};<0mm,0mm>*{}**@{},
   <0.39mm,-0.39mm>*{};<2.9mm,-4mm>*{^{_1}}**@{},
   <-0.35mm,-0.35mm>*{};<-2.8mm,-4mm>*{^{_2}}**@{},
\end{xy}\Ea
\right\rangle
$$
by the ideal generated by the following relation\footnote{Given any module $V$ of the $\K$-algebra $\K[\bS_3]$, the action of the element $\sum_{k=0}^{2} (123)^k\in \K[\bS_3]$ on any $a\in V$ is denoted by $\oint_{123} a$.}
$$
\oint_{123} \Ba{c}\resizebox{10mm}{!}{ \begin{xy}
 <0mm,0mm>*{\bu};<0mm,0mm>*{}**@{},
 <0mm,0.69mm>*{};<0mm,3.0mm>*{}**@{-},
 <0.39mm,-0.39mm>*{};<2.4mm,-2.4mm>*{}**@{-},
 <-0.35mm,-0.35mm>*{};<-1.9mm,-1.9mm>*{}**@{-},
 <-2.4mm,-2.4mm>*{\bu};<-2.4mm,-2.4mm>*{}**@{},
 <-2.0mm,-2.8mm>*{};<0mm,-4.9mm>*{}**@{-},
 <-2.8mm,-2.9mm>*{};<-4.7mm,-4.9mm>*{}**@{-},
    <0.39mm,-0.39mm>*{};<3.3mm,-4.0mm>*{^3}**@{},
    <-2.0mm,-2.8mm>*{};<0.5mm,-6.7mm>*{^2}**@{},
    <-2.8mm,-2.9mm>*{};<-5.2mm,-6.7mm>*{^1}**@{},
 \end{xy}}\Ea=0
$$
Its representations in a dg vector space $V$ are in one-to-one correspondence with Lie algebra structures on $V$ which have Lie brackets,
$$
[\ ,\ ]: \odot^2(V[d])\rar V[1+d]
$$
of degree $1-d$; the case $d=1$ corresponds to the ordinary Lie algebra structure in $V$.

\subsection{Twisting of (prop)operads under $\Lie_d$}\label{2: subsec on twisting of (prop)erads under Lie} Let $\cP=\{\cP(m,n)\}_{m,n\geq 0}$ be a dg properad with the differential denoted by $\sd$. We represent its  generic elements pictorially as $(m,n)$-corollas
\Beq\label{2: generic elements of cP as (m,n)-corollas}
\Ba{c}\resizebox{16mm}{!}{
 \begin{xy}
 <0mm,0mm>*{\circ};<-8mm,6mm>*{^1}**@{-},
 <0mm,0mm>*{\circ};<-4.5mm,6mm>*{^2}**@{-},
 <0mm,0mm>*{\circ};<0mm,5.5mm>*{\ldots}**@{},
 <0mm,0mm>*{\circ};<3.5mm,5mm>*{}**@{-},
 <0mm,0mm>*{\circ};<8mm,6mm>*{^m}**@{-},
 <0mm,0mm>*{\circ};<-8mm,-6mm>*{_1}**@{-},
 <0mm,0mm>*{\circ};<-4.5mm,-6mm>*{_2}**@{-},
 <0mm,0mm>*{\circ};<0mm,-5.5mm>*{\ldots}**@{},
 <0mm,0mm>*{\circ};<4.5mm,-6mm>*+{}**@{-},
 <0mm,0mm>*{\circ};<8mm,-6mm>*{_n}**@{-},
   \end{xy}}\Ea
\Eeq
whose white vertex is decorated by an element of $\cP(m,n)$. Properadic compositions in $\cP$ are represented pictorially by gluing out-legs of such decorated corollas to in-legs  of another decorated corollas.

\sip

Assume that $\cP$ comes equipped with a morphism of properads
$$
\Ba{rccc}
\imath: & (\Lie_d, 0) & \lon & (\cP,\sd)\\
& \Ba{c}\begin{xy}
 <0mm,0.66mm>*{};<0mm,4mm>*{}**@{-},
 <0.39mm,-0.39mm>*{};<2.2mm,-2.2mm>*{}**@{-},
 <-0.35mm,-0.35mm>*{};<-2.2mm,-2.2mm>*{}**@{-},
 <0mm,0mm>*{\bu};<0mm,0mm>*{}**@{},
   <0.39mm,-0.39mm>*{};<2.9mm,-4mm>*{^{_2}}**@{},
   <-0.35mm,-0.35mm>*{};<-2.8mm,-4mm>*{^{_1}}**@{},
\end{xy}\Ea & \lon & \Ba{c}\resizebox{8mm}{!}{  \xy
(-5,6)*{}="1";
    (-5,+1)*{\circledcirc}="L";
  (-8,-5)*+{_1}="C";
   (-2,-5)*+{_2}="D";
\ar @{-} "D";"L" <0pt>
\ar @{-} "C";"L" <0pt>
\ar @{-} "1";"L" <0pt>
 \endxy}
 \Ea
\Ea
$$
where $\Lie_d$ is understood as a differential operad with trivial differential. The image of the generator of $\Lie_d$ is a special element in $\cP$ represented as a $(1,2)$-corolla with the vertex denoted by $\circledcirc$.

 \sip

 Following \cite{W}, one defines a twisted dg properad $(\tw \cP, \sd_\centerdot)$
 as the unique properad freely generated by $\cP$ and an extra $(1,0)$ generator
 $
  \Ba{c}\resizebox{2mm}{!}{\begin{xy}
 <0mm,0.5mm>*{};<0mm,6mm>*{}**@{-},
 <0mm,0mm>*{\blacksquare};<0mm,0mm>*{}**@{},
 \end{xy}}\Ea
 $
of cohomological degree $d$. The differential $\sd_\centerdot$ in $\tw \cP$ is given explicitly on the elements of $\cP$ by
 \Beq\label{2: d_centerdot on twP under Lie}
\sd_\centerdot \Ba{c}\resizebox{16mm}{!}{
 \begin{xy}
 <0mm,0mm>*{\circ};<-8mm,6mm>*{^1}**@{-},
 <0mm,0mm>*{\circ};<-4.5mm,6mm>*{^2}**@{-},
 <0mm,0mm>*{\circ};<0mm,5.5mm>*{\ldots}**@{},
 <0mm,0mm>*{\circ};<3.5mm,5mm>*{}**@{-},
 <0mm,0mm>*{\circ};<8mm,6mm>*{^m}**@{-},
 <0mm,0mm>*{\circ};<-8mm,-6mm>*{_1}**@{-},
 <0mm,0mm>*{\circ};<-4.5mm,-6mm>*{_2}**@{-},
 <0mm,0mm>*{\circ};<0mm,-5.5mm>*{\ldots}**@{},
 <0mm,0mm>*{\circ};<4.5mm,-6mm>*+{}**@{-},
 <0mm,0mm>*{\circ};<8mm,-6mm>*{_n}**@{-},
   \end{xy}}\Ea
=
\sd \Ba{c}\resizebox{16mm}{!}{
 \begin{xy}
 <0mm,0mm>*{\circ};<-8mm,6mm>*{^1}**@{-},
 <0mm,0mm>*{\circ};<-4.5mm,6mm>*{^2}**@{-},
 <0mm,0mm>*{\circ};<0mm,5.5mm>*{\ldots}**@{},
 <0mm,0mm>*{\circ};<3.5mm,5mm>*{}**@{-},
 <0mm,0mm>*{\circ};<8mm,6mm>*{^m}**@{-},
 <0mm,0mm>*{\circ};<-8mm,-6mm>*{_1}**@{-},
 <0mm,0mm>*{\circ};<-4.5mm,-6mm>*{_2}**@{-},
 <0mm,0mm>*{\circ};<0mm,-5.5mm>*{\ldots}**@{},
 <0mm,0mm>*{\circ};<4.5mm,-6mm>*+{}**@{-},
 <0mm,0mm>*{\circ};<8mm,-6mm>*{_n}**@{-},
   \end{xy}}\Ea
+
\overset{m-1}{\underset{i=0}{\sum}}
\Ba{c}\resizebox{17mm}{!}{
\begin{xy}
 <0mm,0mm>*{\circ};<-8mm,5mm>*{}**@{-},
 <0mm,0mm>*{\circ};<-3.5mm,5mm>*{}**@{-},
 <0mm,0mm>*{\circ};<-6mm,5mm>*{..}**@{},
 <0mm,0mm>*{\circ};<0mm,5mm>*{}**@{-},
  <0mm,13mm>*{\circledcirc};
  <0mm,13mm>*{};<5mm,10mm>*{_\blacksquare}**@{-},
  <0mm,5mm>*{};<0mm,12mm>*{}**@{-},
  <0mm,14mm>*{};<0mm,17mm>*{}**@{-},
  <0mm,5mm>*{};<0mm,19mm>*{^{i\hspace{-0.2mm}+\hspace{-0.5mm}1}}**@{},
<0mm,0mm>*{\circ};<8mm,5mm>*{}**@{-},
<0mm,0mm>*{\circ};<3.5mm,5mm>*{}**@{-},
<6mm,5mm>*{..}**@{},
<-8.5mm,5.5mm>*{^1}**@{},
<-4mm,5.5mm>*{^i}**@{},
<9.0mm,5.5mm>*{^m}**@{},
 <0mm,0mm>*{\circ};<-8mm,-5mm>*{}**@{-},
 <0mm,0mm>*{\circ};<-4.5mm,-5mm>*{}**@{-},
 <-1mm,-5mm>*{\ldots}**@{},
 <0mm,0mm>*{\circ};<4.5mm,-5mm>*{}**@{-},
 <0mm,0mm>*{\circ};<8mm,-5mm>*{}**@{-},
<-8.5mm,-6.9mm>*{^1}**@{},
<-5mm,-6.9mm>*{^2}**@{},
<4.5mm,-6.9mm>*{^{n\hspace{-0.5mm}-\hspace{-0.5mm}1}}**@{},
<9.0mm,-6.9mm>*{^n}**@{},
 \end{xy}}\Ea
 - (-1)^{|a|}
\overset{n-1}{\underset{i=0}{\sum}}
 \Ba{c}\resizebox{17mm}{!}{\begin{xy}
 <0mm,0mm>*{\circ};<-8mm,-5mm>*{}**@{-},
 <0mm,0mm>*{\circ};<-3.5mm,-5mm>*{}**@{-},
 <0mm,0mm>*{\circ};<-6mm,-5mm>*{..}**@{},
 <0mm,0mm>*{\circ};<0mm,-5mm>*{}**@{-},
   <0mm,-11mm>*{\circledcirc};
  <0mm,-12mm>*{};<5mm,-16mm>*{_\blacksquare}**@{-},
  <0mm,-5mm>*{};<0mm,-10mm>*{}**@{-},
  <0mm,-12mm>*{};<0mm,-17mm>*{}**@{-},
  <0mm,-5mm>*{};<0mm,-19mm>*{^{i\hspace{-0.2mm}+\hspace{-0.5mm}1}}**@{},
<0mm,0mm>*{\circ};<8mm,-5mm>*{}**@{-},
<0mm,0mm>*{\circ};<3.5mm,-5mm>*{}**@{-},
 <6mm,-5mm>*{..}**@{},
<-8.5mm,-6.9mm>*{^1}**@{},
<-4mm,-6.9mm>*{^i}**@{},
<9.0mm,-6.9mm>*{^n}**@{},
 <0mm,0mm>*{\circ};<-8mm,5mm>*{}**@{-},
 <0mm,0mm>*{\circ};<-4.5mm,5mm>*{}**@{-},
<-1mm,5mm>*{\ldots}**@{},
 <0mm,0mm>*{\circ};<4.5mm,5mm>*{}**@{-},
 <0mm,0mm>*{\circ};<8mm,5mm>*{}**@{-},
<-8.5mm,5.5mm>*{^1}**@{},
<-5mm,5.5mm>*{^2}**@{},
<4.5mm,5.5mm>*{^{m\hspace{-0.5mm}-\hspace{-0.5mm}1}}**@{},
<9.0mm,5.5mm>*{^m}**@{},
 \end{xy}}\Ea,
\Eeq
 and on the extra generator (called the {\it MC generator}) by
 \Beq\label{2: d_c on boxdot}
 \sd_\centerdot
 \Ba{c}\resizebox{2mm}{!}{\begin{xy}
 <0mm,0.5mm>*{};<0mm,6mm>*{}**@{-},
 <0mm,0mm>*{\blacksquare};<0mm,0mm>*{}**@{},
 \end{xy}}\Ea
 =
\frac{1}{2}
\Ba{c}\resizebox{8mm}{!}{  \xy
(-5,6)*{}="1";
    (-5,+1)*{\circledcirc}="L";
  (-8,-5)*{_\blacksquare}="C";
   (-2,-5)*{_\blacksquare}="D";
\ar @{-} "D";"L" <0pt>
\ar @{-} "C";"L" <0pt>
\ar @{-} "1";"L" <0pt>
%
 %
 \endxy}
 \Ea.
 \Eeq
 The twisted properad comes equipped with a natural epimorphism of dg properads
 $$
 (\tw\cP, \sd_\centerdot) \lon (\cP, \sd)
 $$
 which sends the MC generator to zero.

 \sip

 It is easy to check using Jacobi identity in $\Lie_d$ that the element  $\Ba{c}\resizebox{8mm}{!}{  \xy
(-5,6)*{}="1";
    (-5,+1)*{\circledcirc}="L";
  (-8,-5)*+{_1}="C";
   (-2,-5)*+{_2}="D";
\ar @{-} "D";"L" <0pt>
\ar @{-} "C";"L" <0pt>
\ar @{-} "1";"L" <0pt>
 \endxy}
 \Ea$ remains a cocycle even after the twisting of the original differential, $\sd\rar \sd_\centerdot$, in $\cP$,
 $$
 \delta_\centerdot\hspace{-3mm}
\Ba{c}\resizebox{8mm}{!}{  \xy
(-5,6)*{}="1";
    (-5,+1)*{\circledcirc}="L";
  (-8,-5)*+{_1}="C";
   (-2,-5)*+{_2}="D";
\ar @{-} "D";"L" <0pt>
\ar @{-} "C";"L" <0pt>
\ar @{-} "1";"L" <0pt>
 \endxy}
 \Ea
 \equiv
 \Ba{c}\resizebox{11.5mm}{!}{  \xy
(-9,8)*{}="1";
    (-9,+3)*{\circledcirc}="L";
 (-14,-3.5)*{\circledcirc}="B";
 (-18,-10)*+{_1}="b1";
 (-10,-10)*+{_2}="b2";
  (-3,-4)*{\blacksquare}="C";
\ar @{-} "C";"L" <0pt>
\ar @{-} "B";"L" <0pt>
\ar @{-} "B";"b1" <0pt>
\ar @{-} "B";"b2" <0pt>
\ar @{-} "1";"L" <0pt>
 \endxy}
 \Ea
 +(-1)^d
  \Ba{c}\resizebox{10.5mm}{!}{  \xy
(-9,8)*{}="1";
    (-9,+3)*{\circledcirc}="L";
 (-14,-3.5)*{\circledcirc}="B";
 (-18,-10)*{_1}="b1";
 (-10,-10)*{\blacksquare}="b2";
  (-3,-5)*{_2}="C";
\ar @{-} "C";"L" <0pt>
\ar @{-} "B";"L" <0pt>
\ar @{-} "B";"b1" <0pt>
\ar @{-} "B";"b2" <0pt>
\ar @{-} "1";"L" <0pt>
 \endxy}
 \Ea
 +
   \Ba{c}\resizebox{10.5mm}{!}{  \xy
(-9,8)*{}="1";
    (-9,+3)*{\circledcirc}="L";
 (-14,-3.5)*{\circledcirc}="B";
 (-18,-10)*{_2}="b1";
 (-10,-10)*{\blacksquare}="b2";
  (-3,-5)*{_1}="C";
\ar @{-} "C";"L" <0pt>
\ar @{-} "B";"L" <0pt>
\ar @{-} "B";"b1" <0pt>
\ar @{-} "B";"b2" <0pt>
\ar @{-} "1";"L" <0pt>
%
 %
 \endxy}
 \Ea=0.
$$
Hence the original morphism $i$ extends to the twisted version by the same formula,
\Beq\label{2: map i from Lie to twP}
\Ba{rccc}
\imath: & (\Lie_d, 0) & \lon & (\tw\cP,\sd_\centerdot)\\
& \Ba{c}\begin{xy}
 <0mm,0.66mm>*{};<0mm,4mm>*{}**@{-},
 <0.39mm,-0.39mm>*{};<2.2mm,-2.2mm>*{}**@{-},
 <-0.35mm,-0.35mm>*{};<-2.2mm,-2.2mm>*{}**@{-},
 <0mm,0mm>*{\bu};<0mm,0mm>*{}**@{},
   <0.39mm,-0.39mm>*{};<2.9mm,-4mm>*{^{_2}}**@{},
   <-0.35mm,-0.35mm>*{};<-2.8mm,-4mm>*{^{_1}}**@{},
\end{xy}\Ea & \lon & \Ba{c}\resizebox{8mm}{!}{  \xy
(-5,6)*{}="1";
    (-5,+1)*{\circledcirc}="L";
  (-8,-5)*+{_1}="C";
   (-2,-5)*+{_2}="D";
\ar @{-} "D";"L" <0pt>
\ar @{-} "C";"L" <0pt>
\ar @{-} "1";"L" <0pt>
%
 %
 \endxy}
 \Ea
\Ea.
\Eeq

Given a representation,
$$
\rho: \cP \rar \cE nd_V,
$$
of  a prop(erad) $\cP$ in a dg vector space $V$, then composing $\rho$ with the map $i$ we get a dg Lie algebra structure $[\ ,\ ]$ in $V$. Assume $\ga$ is  a Maurer-Cartan element of that dg Lie algebra, that is a degree $d$ element $\ga\in V$ satisfying the equation
$$
d\ga + \frac{1}{2}[\ga,\ga]=0.
$$
Then the pair $(\rho, \ga)$ gives rise to a representation of $\tw\cP$ in $V$ which is given on the elements from $\cP$ by the map $\rho$, and on the extra MC generator by the association,
$$
\Ba{c}\resizebox{2mm}{!}{\begin{xy}
 <0mm,0.5mm>*{};<0mm,6mm>*{}**@{-},
 <0mm,0mm>*{\blacksquare};<0mm,0mm>*{}**@{},
 \end{xy}}\Ea \lon \ga,
$$
which explains the main idea of the twisting endofunctor and the terminology.

\subsection{Twisted operad of ribbon trees and the gravity operad}\label{2: subsec on grav operad and  Rtrees} Let $\cM_{0,1+n}$, $n \geq 2$, be the moduli space of rational curves with $1+n$ distinct points marked by the set $[n^+]:=\{0,1,\ldots, n\}$, that is, the space of injections of $[n^+]$ into the complex projective line,
$$
\cM_{0,1+n}:=\frac{\{[n^+] \hook \C\P^1\}}{PGL(2,\C)},
$$
modulo the action of the group $PGL(2,\C)=SL(2,\C)/\{\Id, -\Id\}$. E.\ Getzler has shown in \cite{Ge} that the $\bS$-module
$$
\cG ravity:=\{H^{\bu-1}_c(\cM_{0,1+n})\simeq H_{\bu-2n+3}(\cM_{0,1+n})\}_{n\geq 2}
$$
has the structure of an operad. Its representation in a graded vector space $V$ is given by a collection of linear maps,
$$
\{\ ,\ \ldots ,\ \}_n: \odot^n V \lon V[2n-3], \ \ \ n\geq 2,
$$
satisfying the relations
 $$
\{\{v_1,\ldots,v_k\},v_{k+1},\ldots,v_{k+l}\}=\hfill
$$
$$\sum_{1\leq i< j\leq k} (-1)^{\var(i,j)} \{\{v_i,v_j\},v_1,\ldots,\widehat{v_i},\ldots, \widehat{v_j},\ldots, v_k, v_{k+1},\ldots,v_{k+l}\}
 $$
for any $v_1,\ldots, v_{k+l}\in V$; here $(-1)^{\var(i,j)}$ is the standard Koszul sign.
 We refer to \cite{KMS} for a further study of the gravity operad, its Koszul dual operad, and its relation to the operad of Gerstenhaber algebras.

\sip

There is a nice chain representation of the gravity operad given in terms of twisted ribbon trees which was constructed in \cite{Wa}. Let $\RGra_d$, $d\in \Z$, be the properad of {\em connected}\, ribbon graphs  introduced in \cite{MW} (see \S {\ref{4: subsec on RGra}} below for a short reminder of its definition and main properties), and let $\cR\cT ree_d$ be its sub-operad generated by ribbon trees, that is, by ribbon graphs of genus zero with precisely one boundary, e.g.
$$
\Ba{c}\resizebox{15mm}{!}{
\mbox{$\xy
 (0,0)*+{_{_4}}*\frm{o}="C";
  (9,0)*+{_{_1}}*\frm{o}="1";
(-7,8)*+{_{_2}}*\frm{o}="2";
(-7,-8)*+{_{_3}}*\frm{o}="3";
 \ar @{-} "1";"C" <0pt>
 \ar @{-} "2";"C" <0pt>
  \ar @{-} "3";"C" <0pt>
\endxy$}}
\Ea
\in \cR\cT ree_d(1,4).
$$
This operad comes equipped with a morphism
$$
\Ba{ccc}
\Lie_{d} & \lon & \cR\cT ree_{d}\\
\Ba{c}\resizebox{6mm}{!}{\begin{xy}
 <0mm,0.66mm>*{};<0mm,3mm>*{}**@{-},
 <0.39mm,-0.39mm>*{};<2.2mm,-2.2mm>*{}**@{-},
 <-0.35mm,-0.35mm>*{};<-2.2mm,-2.2mm>*{}**@{-},
 <0mm,0mm>*{\bu};<0mm,0mm>*{}**@{},
   <0.39mm,-0.39mm>*{};<2.9mm,-4mm>*{^2}**@{},
   <-0.35mm,-0.35mm>*{};<-2.8mm,-4mm>*{^1}**@{},
\end{xy}}\Ea
&\lon&
\Ba{c} \resizebox{9mm}{!}{\xy
 (7,0)*+{_2}*\frm{o}="B";
 (0,0)*+{_1}*\frm{o}="A";
 \ar @{-} "A";"B" <0pt>
\endxy} \Ea
\Ea
$$
and hence can be twisted using T.\ Willwacher's endofunctor $\tw$. The associated dg operad
$\tw\cR\cT ree_d$ is generated by $\cR\cT ree_d$ and an extra element of cohomological degree $d$ with no inputs and one boundary which we identify with the ribbon graph $\bu$ consisting of one black vertex and one boundary; composing that boundary with, for example, the input vertex 4 of the above graph in $\cR\cT ree_d(1,4)$ one gets a ribbon graph
$$
\Ba{c}\resizebox{15mm}{!}{
\mbox{$\xy
 (0,0)*{\bu}="C";
  (7.9,0)*{_{_1}}*+\frm{o}="1";
(-6,7)*{_{_2}}*+\frm{o}="2";
(-6,-7)*{_{_3}}*+\frm{o}="3";
 \ar @{-} "1";"C" <0pt>
 \ar @{-} "2";"C" <0pt>
  \ar @{-} "3";"C" <0pt>
\endxy$}}
\Ea
\in \tw\cR\cT ree_d(1,3).
$$
Thus $\tw\cR\cT ree_d(1,3)$ is generated by ribbon trees with vertices of two types,
white ones which are labelled and black ones which are unlabelled and have the cohomological degree $d$ assigned. According to the general formula (\ref{2: d_centerdot on twP under Lie}), the induced differential in
$\tw\cR\cT ree_d(1,3)$ is given by
\Beq\label{5: d in Tw(RTree)}
d_\centerdot\Ga:=
\Ba{c}\resizebox{3mm}{!}{  \xy
 (0,6)*{\bu}="A";
 (0,0)*+{_1}*\frm{o}="B";
 \ar @{-} "A";"B" <0pt>
\endxy} \Ea
  \circ_1 \Ga\ \
- \ \ (-1)^{|\Ga|} \sum_{v\in V_\circ(\Ga)}  \Ga \circ_v
\Ba{c}\resizebox{3mm}{!}{  \xy
 (0,6)*{\bu}="A";
 (0,0)*+{_1}*\frm{o}="B";
 \ar @{-} "A";"B" <0pt>
\endxy} \Ea
\  -(-1)^{|\Ga|}\ \frac{1}{2} \sum_{v\in V_\bu(\Ga)} \Ga\circ_v  \left(\xy
 (0,0)*{\bullet}="a",
(5,0)*{\bu}="b",
\ar @{-} "a";"b" <0pt>
\endxy\right)
\Eeq
where the symbol $\ga_1 \circ_v \ga_2$ means substitution of the ribbon tree $\ga_2$ into the vertex $v$ of the ribbon tree $\ga_1$ and re-distributing edges attached to $v$ among all (i.e. both white and black (if any)) vertices of $\ga_2$ in all possible ways while respecting their cyclic order while moving along its unique boundary.

\sip

It has been proven in \cite{Wa} that the cohomology of this twisted operad is precisely the degree shifted gravity operad,
$$
H^\bu(\tw \cR\cT ree_d) =\cG ravity\{d\}.
$$
 The elements of $\cG ravity\{d\}$ in this chain complex  get represented as equivalence classes of twisted ribbon trees with all white vertices being univalent and all black vertices (if any) being trivalent. It is worth noting that $\tw \cR\cT ree_d$ is quasi-isomorphic as a dg operad to its suboperad generated by ribbon trees with all black vertices being at least trivalent (cf.\ Proposition {\ref{3: prop on tw^3Rgra to twRgra}} below).

\sip

\bip

\bip


{\Large
\section{\bf Twisted properad of ribbon graphs and Kevin Costello's moduli spaces}
}

\mip

\subsection{Degree shifted quasi-Lie bialgebras and their deformation complex}
Recall \cite{MW2} that the prop(erad) of degree shifted Lie bialgebras is defined, for any pair of integer  $c,d\in \Z$, as the quotient
$$
\LB_{c,d}:=\cF ree\langle E_0\rangle/\langle\cR\rangle,
$$
of the free prop(erad) generated by an  $\bS$-bimodule $E=\{E(m,n)\}_{m,n\geq 0}$ with
 all $E_0(m,n)=0$ except\footnote{All graphs considered in \S 3.1 are assumed tacitly to have directed edges and legs with the direction flow running from bottom to the top.}
  $$
E_0(2,1):=\id_1\ot (\sgn_2)^{\ot |c|}[c-1]=\mbox{span}\left\langle
\Ba{c}\begin{xy}
 <0mm,-0.55mm>*{};<0mm,-2.5mm>*{}**@{-},
 <0.5mm,0.5mm>*{};<2.2mm,2.2mm>*{}**@{-},
 <-0.48mm,0.48mm>*{};<-2.2mm,2.2mm>*{}**@{-},
 <0mm,0mm>*{\bu};<0mm,0mm>*{}**@{},
 <0.5mm,0.5mm>*{};<2.7mm,2.8mm>*{^{_2}}**@{},
 <-0.48mm,0.48mm>*{};<-2.7mm,2.8mm>*{^{_1}}**@{},
 \end{xy}\Ea
=(-1)^{c}
\Ba{c}\begin{xy}
 <0mm,-0.55mm>*{};<0mm,-2.5mm>*{}**@{-},
 <0.5mm,0.5mm>*{};<2.2mm,2.2mm>*{}**@{-},
 <-0.48mm,0.48mm>*{};<-2.2mm,2.2mm>*{}**@{-},
 <0mm,0mm>*{\bu};<0mm,0mm>*{}**@{},
 <0.5mm,0.5mm>*{};<2.7mm,2.8mm>*{^{_1}}**@{},
 <-0.48mm,0.48mm>*{};<-2.7mm,2.8mm>*{^{_2}}**@{},
 \end{xy}\Ea
   \right\rangle
$$
$$
E_0(1,2):= (\sgn_2)^{\ot |d|}\ot \id_1[d-1]=\mbox{span}\left\langle
\Ba{c}\begin{xy}
 <0mm,0.66mm>*{};<0mm,3mm>*{}**@{-},
 <0.39mm,-0.39mm>*{};<2.2mm,-2.2mm>*{}**@{-},
 <-0.35mm,-0.35mm>*{};<-2.2mm,-2.2mm>*{}**@{-},
 <0mm,0mm>*{\bu};<0mm,0mm>*{}**@{},
   <0.39mm,-0.39mm>*{};<2.9mm,-4mm>*{^{_2}}**@{},
   <-0.35mm,-0.35mm>*{};<-2.8mm,-4mm>*{^{_1}}**@{},
\end{xy}\Ea
=(-1)^{d}
\Ba{c}\begin{xy}
 <0mm,0.66mm>*{};<0mm,3mm>*{}**@{-},
 <0.39mm,-0.39mm>*{};<2.2mm,-2.2mm>*{}**@{-},
 <-0.35mm,-0.35mm>*{};<-2.2mm,-2.2mm>*{}**@{-},
 <0mm,0mm>*{\bu};<0mm,0mm>*{}**@{},
   <0.39mm,-0.39mm>*{};<2.9mm,-4mm>*{^{_1}}**@{},
   <-0.35mm,-0.35mm>*{};<-2.8mm,-4mm>*{^{_2}}**@{},
\end{xy}\Ea
\right\rangle
$$
by the ideal generated by the following relations
\Beq\label{3: R for LieB}
\cR:\left\{
\Ba{c}
\displaystyle
\oint_{123} \Ba{c}\resizebox{8.4mm}{!}{
\begin{xy}
 <0mm,0mm>*{\bu};<0mm,0mm>*{}**@{},
 <0mm,-0.49mm>*{};<0mm,-3.0mm>*{}**@{-},
 <0.49mm,0.49mm>*{};<1.9mm,1.9mm>*{}**@{-},
 <-0.5mm,0.5mm>*{};<-1.9mm,1.9mm>*{}**@{-},
 <-2.3mm,2.3mm>*{\bu};<-2.3mm,2.3mm>*{}**@{},
 <-1.8mm,2.8mm>*{};<0mm,4.9mm>*{}**@{-},
 <-2.8mm,2.9mm>*{};<-4.6mm,4.9mm>*{}**@{-},
   <0.49mm,0.49mm>*{};<2.7mm,2.3mm>*{^3}**@{},
   <-1.8mm,2.8mm>*{};<0.4mm,5.3mm>*{^2}**@{},
   <-2.8mm,2.9mm>*{};<-5.1mm,5.3mm>*{^1}**@{},
 \end{xy}}\Ea =0
 \ , \ \ \ \ \
\oint_{123} \Ba{c}\resizebox{8.4mm}{!}{ \begin{xy}
 <0mm,0mm>*{\bu};<0mm,0mm>*{}**@{},
 <0mm,0.69mm>*{};<0mm,3.0mm>*{}**@{-},
 <0.39mm,-0.39mm>*{};<2.4mm,-2.4mm>*{}**@{-},
 <-0.35mm,-0.35mm>*{};<-1.9mm,-1.9mm>*{}**@{-},
 <-2.4mm,-2.4mm>*{\bu};<-2.4mm,-2.4mm>*{}**@{},
 <-2.0mm,-2.8mm>*{};<0mm,-4.9mm>*{}**@{-},
 <-2.8mm,-2.9mm>*{};<-4.7mm,-4.9mm>*{}**@{-},
    <0.39mm,-0.39mm>*{};<3.3mm,-4.0mm>*{^3}**@{},
    <-2.0mm,-2.8mm>*{};<0.5mm,-6.7mm>*{^2}**@{},
    <-2.8mm,-2.9mm>*{};<-5.2mm,-6.7mm>*{^1}**@{},
 \end{xy}}\Ea =0
\vspace{2mm} \\
 \Ba{c}\resizebox{6mm}{!}{\begin{xy}
 <0mm,2.47mm>*{};<0mm,0.12mm>*{}**@{-},
 <0.5mm,3.5mm>*{};<2.2mm,5.2mm>*{}**@{-},
 <-0.48mm,3.48mm>*{};<-2.2mm,5.2mm>*{}**@{-},
 <0mm,3mm>*{\bu};<0mm,3mm>*{}**@{},
  <0mm,-0.8mm>*{\bu};<0mm,-0.8mm>*{}**@{},
<-0.39mm,-1.2mm>*{};<-2.2mm,-3.5mm>*{}**@{-},
 <0.39mm,-1.2mm>*{};<2.2mm,-3.5mm>*{}**@{-},
     <0.5mm,3.5mm>*{};<2.8mm,5.7mm>*{^2}**@{},
     <-0.48mm,3.48mm>*{};<-2.8mm,5.7mm>*{^1}**@{},
   <0mm,-0.8mm>*{};<-2.7mm,-5.2mm>*{^1}**@{},
   <0mm,-0.8mm>*{};<2.7mm,-5.2mm>*{^2}**@{},
\end{xy}}\Ea
+(-1)^{cd+c+d}\left(
\Ba{c}\resizebox{7mm}{!}{\begin{xy}
 <0mm,-1.3mm>*{};<0mm,-3.5mm>*{}**@{-},
 <0.38mm,-0.2mm>*{};<2.0mm,2.0mm>*{}**@{-},
 <-0.38mm,-0.2mm>*{};<-2.2mm,2.2mm>*{}**@{-},
<0mm,-0.8mm>*{\bu};<0mm,0.8mm>*{}**@{},
 <2.4mm,2.4mm>*{\bu};<2.4mm,2.4mm>*{}**@{},
 <2.77mm,2.0mm>*{};<4.4mm,-0.8mm>*{}**@{-},
 <2.4mm,3mm>*{};<2.4mm,5.2mm>*{}**@{-},
     <0mm,-1.3mm>*{};<0mm,-5.3mm>*{^1}**@{},
     <2.5mm,2.3mm>*{};<5.1mm,-2.6mm>*{^2}**@{},
    <2.4mm,2.5mm>*{};<2.4mm,5.7mm>*{^2}**@{},
    <-0.38mm,-0.2mm>*{};<-2.8mm,2.5mm>*{^1}**@{},
    \end{xy}}\Ea
  + (-1)^{d}
\Ba{c}\resizebox{7mm}{!}{\begin{xy}
 <0mm,-1.3mm>*{};<0mm,-3.5mm>*{}**@{-},
 <0.38mm,-0.2mm>*{};<2.0mm,2.0mm>*{}**@{-},
 <-0.38mm,-0.2mm>*{};<-2.2mm,2.2mm>*{}**@{-},
<0mm,-0.8mm>*{\bu};<0mm,0.8mm>*{}**@{},
 <2.4mm,2.4mm>*{\bu};<2.4mm,2.4mm>*{}**@{},
 <2.77mm,2.0mm>*{};<4.4mm,-0.8mm>*{}**@{-},
 <2.4mm,3mm>*{};<2.4mm,5.2mm>*{}**@{-},
     <0mm,-1.3mm>*{};<0mm,-5.3mm>*{^2}**@{},
     <2.5mm,2.3mm>*{};<5.1mm,-2.6mm>*{^1}**@{},
    <2.4mm,2.5mm>*{};<2.4mm,5.7mm>*{^2}**@{},
    <-0.38mm,-0.2mm>*{};<-2.8mm,2.5mm>*{^1}**@{},
    \end{xy}}\Ea
  + (-1)^{d+c}
\Ba{c}\resizebox{7mm}{!}{\begin{xy}
 <0mm,-1.3mm>*{};<0mm,-3.5mm>*{}**@{-},
 <0.38mm,-0.2mm>*{};<2.0mm,2.0mm>*{}**@{-},
 <-0.38mm,-0.2mm>*{};<-2.2mm,2.2mm>*{}**@{-},
<0mm,-0.8mm>*{\bu};<0mm,0.8mm>*{}**@{},
 <2.4mm,2.4mm>*{\bu};<2.4mm,2.4mm>*{}**@{},
 <2.77mm,2.0mm>*{};<4.4mm,-0.8mm>*{}**@{-},
 <2.4mm,3mm>*{};<2.4mm,5.2mm>*{}**@{-},
     <0mm,-1.3mm>*{};<0mm,-5.3mm>*{^2}**@{},
     <2.5mm,2.3mm>*{};<5.1mm,-2.6mm>*{^1}**@{},
    <2.4mm,2.5mm>*{};<2.4mm,5.7mm>*{^1}**@{},
    <-0.38mm,-0.2mm>*{};<-2.8mm,2.5mm>*{^2}**@{},
    \end{xy}}\Ea
 + (-1)^{c}
\Ba{c}\resizebox{7mm}{!}{\begin{xy}
 <0mm,-1.3mm>*{};<0mm,-3.5mm>*{}**@{-},
 <0.38mm,-0.2mm>*{};<2.0mm,2.0mm>*{}**@{-},
 <-0.38mm,-0.2mm>*{};<-2.2mm,2.2mm>*{}**@{-},
<0mm,-0.8mm>*{\bu};<0mm,0.8mm>*{}**@{},
 <2.4mm,2.4mm>*{\bu};<2.4mm,2.4mm>*{}**@{},
 <2.77mm,2.0mm>*{};<4.4mm,-0.8mm>*{}**@{-},
 <2.4mm,3mm>*{};<2.4mm,5.2mm>*{}**@{-},
     <0mm,-1.3mm>*{};<0mm,-5.3mm>*{^1}**@{},
     <2.5mm,2.3mm>*{};<5.1mm,-2.6mm>*{^2}**@{},
    <2.4mm,2.5mm>*{};<2.4mm,5.7mm>*{^1}**@{},
    <-0.38mm,-0.2mm>*{};<-2.8mm,2.5mm>*{^2}**@{},
    \end{xy}}\Ea\right)=0.
    \Ea
\right.
\Eeq
where the vertices are assumed to be ordered in such a way that the ones on the top come first.
Its arbitrary representation,
$$
\rho: \LBcd \lon \cE nd_V,
$$
in a dg vector space $V$ is uniquely determined
by the values of $\rho$ on the generators,
$$
\rho \left(
 \begin{xy}
 <0mm,-0.55mm>*{};<0mm,-2.5mm>*{}**@{-},
 <0.5mm,0.5mm>*{};<2.2mm,2.2mm>*{}**@{-},
 <-0.48mm,0.48mm>*{};<-2.2mm,2.2mm>*{}**@{-},
 <0mm,0mm>*{\bu};<0mm,0mm>*{}**@{},
 \end{xy}\right): V[-c]\rar \odot^2(V[-c])[1], \  \ \ \
\left(
 \begin{xy}
 <0mm,0.66mm>*{};<0mm,3mm>*{}**@{-},
 <0.39mm,-0.39mm>*{};<2.2mm,-2.2mm>*{}**@{-},
 <-0.35mm,-0.35mm>*{};<-2.2mm,-2.2mm>*{}**@{-},
 <0mm,0mm>*{\bu};<0mm,0mm>*{}**@{},
 \end{xy}
 \right): \odot^2(V[d]) \rar V[1+d],
$$
which make  $V$ into a (degree shifted) dg Lie algebra and dg Lie coalgebra, with Lie bracket and Lie cobracket  satisfying a compatibility condition. This notion for $c=d=1$ was introduced by V.\ Drinfeld in the context of the deformation theory of the universal enveloping algebras in the category of quantum groups \cite{Dr1}. The homotopy theory of Lie $(c,d)$-bialgebras is highly non-trivial --- it is controlled by M.\ Kontsevich's graph complex $\GC_{c+d}$
(see \cite{MW2} for the proof); in particular, for $c=d=1$  the famous and mysterious Grothendieck-Teichm\"uller group $GRT$ \cite{Dr3} acts on the genus completed version of the properad $\LB_{1,1}$ as homotopy non-trivial automorphisms. The minimal resolution of $\LBcd$ is denoted by $\HoLBcd$.

\sip

V.\ Drinfeld introduced also in the same context of quantum groups the notion of {\it quasi-Lie bialgebra}\ or {\it Lie quasi-bialgebra} \cite{Dr2}. We construct  in this paper a new example of this structure in the context of the theory of cohomology groups of moduli spaces $\cM_{g,n}$.
The prop(erad) of degree shifted quasi-Lie bialgebras can be defined, for any pair of integer  $c,d\in \Z$, as the quotient
$$
q\LB_{c,d}:=\cF ree\langle E_q\rangle/\langle\cR_q\rangle,
$$
of the free prop(erad) generated by an  $\bS$-bimodule $E_q=\{E_q(m,n)\}_{m,n\geq 0}$ with
 all $E_q(m,n)=0$ except
\Beqrn
E_q(2,1)&:=&\id_1\ot \sgn_2^{c}[c-1]=\mbox{span}\left\langle
\Ba{c}\begin{xy}
 <0mm,-0.55mm>*{};<0mm,-3mm>*{}**@{-},
 <0.5mm,0.5mm>*{};<2.2mm,2.2mm>*{}**@{-},
 <-0.48mm,0.48mm>*{};<-2.2mm,2.2mm>*{}**@{-},
 <0mm,0mm>*{\bu};<0mm,0mm>*{}**@{},
 <0.5mm,0.5mm>*{};<2.7mm,2.8mm>*{^{_2}}**@{},
 <-0.48mm,0.48mm>*{};<-2.7mm,2.8mm>*{^{_1}}**@{},
 \end{xy}\Ea
=(-1)^{c}
\Ba{c}\begin{xy}
 <0mm,-0.55mm>*{};<0mm,-3mm>*{}**@{-},
 <0.5mm,0.5mm>*{};<2.2mm,2.2mm>*{}**@{-},
 <-0.48mm,0.48mm>*{};<-2.2mm,2.2mm>*{}**@{-},
 <0mm,0mm>*{\bu};<0mm,0mm>*{}**@{},
 <0.5mm,0.5mm>*{};<2.7mm,2.8mm>*{^{_1}}**@{},
 <-0.48mm,0.48mm>*{};<-2.7mm,2.8mm>*{^{_2}}**@{},
 \end{xy}\Ea
   \right\rangle,\\
E_q(1,2)&:=& \sgn_2^{d}\ot \id_1[d-1]=\mbox{span}\left\langle
\Ba{c}\begin{xy}
 <0mm,0.66mm>*{};<0mm,3mm>*{}**@{-},
 <0.39mm,-0.39mm>*{};<2.2mm,-2.2mm>*{}**@{-},
 <-0.35mm,-0.35mm>*{};<-2.2mm,-2.2mm>*{}**@{-},
 <0mm,0mm>*{\bu};<0mm,0mm>*{}**@{},
   <0.39mm,-0.39mm>*{};<2.9mm,-4mm>*{^{_2}}**@{},
   <-0.35mm,-0.35mm>*{};<-2.8mm,-4mm>*{^{_1}}**@{},
\end{xy}\Ea
=(-1)^{d}
\Ba{c}\begin{xy}
 <0mm,0.66mm>*{};<0mm,3mm>*{}**@{-},
 <0.39mm,-0.39mm>*{};<2.2mm,-2.2mm>*{}**@{-},
 <-0.35mm,-0.35mm>*{};<-2.2mm,-2.2mm>*{}**@{-},
 <0mm,0mm>*{\bu};<0mm,0mm>*{}**@{},
   <0.39mm,-0.39mm>*{};<2.9mm,-4mm>*{^{_1}}**@{},
   <-0.35mm,-0.35mm>*{};<-2.8mm,-4mm>*{^{_2}}**@{},
\end{xy}\Ea
\right\rangle,\\
E_q(3,0)&:=& (\sgn_3)^{\ot|c|}[2c-d-1]=\mbox{span}\left\langle
\Ba{c}\begin{xy}
 <0mm,-1mm>*{\bu};<-4mm,3mm>*{^{_1}}**@{-},
 <0mm,-1mm>*{\bu};<0mm,3mm>*{^{_2}}**@{-},
 <0mm,-1mm>*{\bu};<4mm,3mm>*{^{_3}}**@{-},
 \end{xy}\Ea= (-1)^{c|\sigma|}
 \Ba{c}\begin{xy}
 <0mm,-1mm>*{\bu};<-6mm,3mm>*{^{_{\sigma(1)}}}**@{-},
 <0mm,-1mm>*{\bu};<0mm,3mm>*{^{_{\sigma(2)}}}**@{-},
 <0mm,-1mm>*{\bu};<6mm,3mm>*{^{_{\sigma(3)}}}**@{-},
 \end{xy}\Ea\forall\sigma\in \bS_3
\right\rangle,
\Eeqrn
modulo the ideal generated by the following relations
$$
\cR_q:\left\{
\Ba{c}
\displaystyle
\oint_{123}
\left(\hspace{-2mm} \Ba{c}\resizebox{9.4mm}{!}{
\begin{xy}
 <0mm,0mm>*{\bu};<0mm,0mm>*{}**@{},
 <0mm,-0.49mm>*{};<0mm,-3.0mm>*{}**@{-},
 <0.49mm,0.49mm>*{};<1.9mm,1.9mm>*{}**@{-},
 <-0.5mm,0.5mm>*{};<-1.9mm,1.9mm>*{}**@{-},
 <-2.3mm,2.3mm>*{\bu};<-2.3mm,2.3mm>*{}**@{},
 <-1.8mm,2.8mm>*{};<0mm,4.9mm>*{}**@{-},
 <-2.8mm,2.9mm>*{};<-4.6mm,4.9mm>*{}**@{-},
   <0.49mm,0.49mm>*{};<2.7mm,2.3mm>*{^3}**@{},
   <-1.8mm,2.8mm>*{};<0.4mm,5.3mm>*{^2}**@{},
   <-2.8mm,2.9mm>*{};<-5.1mm,5.3mm>*{^1}**@{},
 \end{xy}}\Ea
+
\Ba{c}\resizebox{12.5mm}{!}{
\begin{xy}
(0,0)*{\bu};(-4,5)*{^{1}}**@{-},
(0,0)*{\bu};(0,5)*{^{2}}**@{-},
(0,0)*{\bu};(4,5)*{\bu}**@{-},
(4,5)*{\bu};(4,10)*{^{3}}**@{-},
(4,5)*{\bu};(8,0)*{_{\, 1}}**@{-},
\end{xy}
}
\Ea
\hspace{-2mm}
\right)
=0
 \ , \ \ \
\oint_{123}\hspace{-1mm} \Ba{c}\resizebox{10.0mm}{!}{ \begin{xy}
 <0mm,0mm>*{\bu};<0mm,0mm>*{}**@{},
 <0mm,0.69mm>*{};<0mm,3.0mm>*{}**@{-},
 <0.39mm,-0.39mm>*{};<2.4mm,-2.4mm>*{}**@{-},
 <-0.35mm,-0.35mm>*{};<-1.9mm,-1.9mm>*{}**@{-},
 <-2.4mm,-2.4mm>*{\bu};<-2.4mm,-2.4mm>*{}**@{},
 <-2.0mm,-2.8mm>*{};<0mm,-4.9mm>*{}**@{-},
 <-2.8mm,-2.9mm>*{};<-4.7mm,-4.9mm>*{}**@{-},
    <0.39mm,-0.39mm>*{};<3.3mm,-4.0mm>*{^3}**@{},
    <-2.0mm,-2.8mm>*{};<0.5mm,-6.7mm>*{^2}**@{},
    <-2.8mm,-2.9mm>*{};<-5.2mm,-6.7mm>*{^1}**@{},
 \end{xy}}\Ea =0
\ , \ \ \
\oint_{123}\left(\hspace{-1.5mm}
\Ba{c}\resizebox{11.4mm}{!}{
\begin{xy}
(0,0)*{\bu};(-4,5)*{^{1}}**@{-},
(0,0)*{\bu};(0,5)*{^{2}}**@{-},
(0,0)*{\bu};(4,5)*{\bu}**@{-},
(4,5)*{\bu};(1.5,10)*{^{3}}**@{-},
(4,5)*{\bu};(6.5,10)*{^{4}}**@{-},
\end{xy}}
\Ea
+
(-1)^c
\Ba{c}\resizebox{11.4mm}{!}{
\begin{xy}
(0,0)*{\bu};(-4,5)*{^{4}}**@{-},
(0,0)*{\bu};(0,5)*{^{1}}**@{-},
(0,0)*{\bu};(4,5)*{\bu}**@{-},
(4,5)*{\bu};(1.5,10)*{^{2}}**@{-},
(4,5)*{\bu};(6.5,10)*{^{3}}**@{-},
\end{xy}}
\Ea\hspace{-1mm}\right),
=0 \vspace{2mm}\\
 \Ba{c}\resizebox{6mm}{!}{\begin{xy}
 <0mm,2.47mm>*{};<0mm,0.12mm>*{}**@{-},
 <0.5mm,3.5mm>*{};<2.2mm,5.2mm>*{}**@{-},
 <-0.48mm,3.48mm>*{};<-2.2mm,5.2mm>*{}**@{-},
 <0mm,3mm>*{\bu};<0mm,3mm>*{}**@{},
  <0mm,-0.8mm>*{\bu};<0mm,-0.8mm>*{}**@{},
<-0.39mm,-1.2mm>*{};<-2.2mm,-3.5mm>*{}**@{-},
 <0.39mm,-1.2mm>*{};<2.2mm,-3.5mm>*{}**@{-},
     <0.5mm,3.5mm>*{};<2.8mm,5.7mm>*{^2}**@{},
     <-0.48mm,3.48mm>*{};<-2.8mm,5.7mm>*{^1}**@{},
   <0mm,-0.8mm>*{};<-2.7mm,-5.2mm>*{^1}**@{},
   <0mm,-0.8mm>*{};<2.7mm,-5.2mm>*{^2}**@{},
\end{xy}}\Ea
+(-1)^{cd+c+d}\left(
\Ba{c}\resizebox{7mm}{!}{\begin{xy}
 <0mm,-1.3mm>*{};<0mm,-3.5mm>*{}**@{-},
 <0.38mm,-0.2mm>*{};<2.0mm,2.0mm>*{}**@{-},
 <-0.38mm,-0.2mm>*{};<-2.2mm,2.2mm>*{}**@{-},
<0mm,-0.8mm>*{\bu};<0mm,0.8mm>*{}**@{},
 <2.4mm,2.4mm>*{\bu};<2.4mm,2.4mm>*{}**@{},
 <2.77mm,2.0mm>*{};<4.4mm,-0.8mm>*{}**@{-},
 <2.4mm,3mm>*{};<2.4mm,5.2mm>*{}**@{-},
     <0mm,-1.3mm>*{};<0mm,-5.3mm>*{^1}**@{},
     <2.5mm,2.3mm>*{};<5.1mm,-2.6mm>*{^2}**@{},
    <2.4mm,2.5mm>*{};<2.4mm,5.7mm>*{^2}**@{},
    <-0.38mm,-0.2mm>*{};<-2.8mm,2.5mm>*{^1}**@{},
    \end{xy}}\Ea
  + (-1)^{d}
\Ba{c}\resizebox{7mm}{!}{\begin{xy}
 <0mm,-1.3mm>*{};<0mm,-3.5mm>*{}**@{-},
 <0.38mm,-0.2mm>*{};<2.0mm,2.0mm>*{}**@{-},
 <-0.38mm,-0.2mm>*{};<-2.2mm,2.2mm>*{}**@{-},
<0mm,-0.8mm>*{\bu};<0mm,0.8mm>*{}**@{},
 <2.4mm,2.4mm>*{\bu};<2.4mm,2.4mm>*{}**@{},
 <2.77mm,2.0mm>*{};<4.4mm,-0.8mm>*{}**@{-},
 <2.4mm,3mm>*{};<2.4mm,5.2mm>*{}**@{-},
     <0mm,-1.3mm>*{};<0mm,-5.3mm>*{^2}**@{},
     <2.5mm,2.3mm>*{};<5.1mm,-2.6mm>*{^1}**@{},
    <2.4mm,2.5mm>*{};<2.4mm,5.7mm>*{^2}**@{},
    <-0.38mm,-0.2mm>*{};<-2.8mm,2.5mm>*{^1}**@{},
    \end{xy}}\Ea
  + (-1)^{d+c}
\Ba{c}\resizebox{7mm}{!}{\begin{xy}
 <0mm,-1.3mm>*{};<0mm,-3.5mm>*{}**@{-},
 <0.38mm,-0.2mm>*{};<2.0mm,2.0mm>*{}**@{-},
 <-0.38mm,-0.2mm>*{};<-2.2mm,2.2mm>*{}**@{-},
<0mm,-0.8mm>*{\bu};<0mm,0.8mm>*{}**@{},
 <2.4mm,2.4mm>*{\bu};<2.4mm,2.4mm>*{}**@{},
 <2.77mm,2.0mm>*{};<4.4mm,-0.8mm>*{}**@{-},
 <2.4mm,3mm>*{};<2.4mm,5.2mm>*{}**@{-},
     <0mm,-1.3mm>*{};<0mm,-5.3mm>*{^2}**@{},
     <2.5mm,2.3mm>*{};<5.1mm,-2.6mm>*{^1}**@{},
    <2.4mm,2.5mm>*{};<2.4mm,5.7mm>*{^1}**@{},
    <-0.38mm,-0.2mm>*{};<-2.8mm,2.5mm>*{^2}**@{},
    \end{xy}}\Ea
 + (-1)^{c}
\Ba{c}\resizebox{7mm}{!}{\begin{xy}
 <0mm,-1.3mm>*{};<0mm,-3.5mm>*{}**@{-},
 <0.38mm,-0.2mm>*{};<2.0mm,2.0mm>*{}**@{-},
 <-0.38mm,-0.2mm>*{};<-2.2mm,2.2mm>*{}**@{-},
<0mm,-0.8mm>*{\bu};<0mm,0.8mm>*{}**@{},
 <2.4mm,2.4mm>*{\bu};<2.4mm,2.4mm>*{}**@{},
 <2.77mm,2.0mm>*{};<4.4mm,-0.8mm>*{}**@{-},
 <2.4mm,3mm>*{};<2.4mm,5.2mm>*{}**@{-},
     <0mm,-1.3mm>*{};<0mm,-5.3mm>*{^1}**@{},
     <2.5mm,2.3mm>*{};<5.1mm,-2.6mm>*{^2}**@{},
    <2.4mm,2.5mm>*{};<2.4mm,5.7mm>*{^1}**@{},
    <-0.38mm,-0.2mm>*{};<-2.8mm,2.5mm>*{^2}**@{},
    \end{xy}}\Ea\right)=0.
    \Ea
\right.
$$
To understand the deformation theory of properads under $q\LBcd$ one needs its minimal resolution, $\HoqLBcd$,
which is easy to construct using the same ideas  which work fine for ordinary Lie bialgebras. It is a dg quasi-free properad,
$$
\HoqLBcd:=\cF ree \left\langle E\right\rangle
$$
 generated by an $\bS$-bimodule
 $$
 E=\{E(m,n)\}_{m\geq 1, n\geq 0, m+n\geq 3}
 $$
  given by
 $$
{E}(m,n):=\sgn_m^{\ot |c|}\ot \sgn_n^{|d|}[cm+dn-1-c-d]\equiv\text{span}\left\langle\hspace{-1mm}
\Ba{c}\resizebox{17mm}{!}{\begin{xy}
 <0mm,0mm>*{\bu};<0mm,0mm>*{}**@{},
 <-0.6mm,0.44mm>*{};<-8mm,5mm>*{}**@{-},
 <-0.4mm,0.7mm>*{};<-4.5mm,5mm>*{}**@{-},
 <0mm,0mm>*{};<1mm,5mm>*{\ldots}**@{},
 <0.4mm,0.7mm>*{};<4.5mm,5mm>*{}**@{-},
 <0.6mm,0.44mm>*{};<8mm,5mm>*{}**@{-},
   <0mm,0mm>*{};<-10.5mm,5.9mm>*{^{\sigma(1)}}**@{},
   <0mm,0mm>*{};<-4mm,5.9mm>*{^{\sigma(2)}}**@{},
   <0mm,0mm>*{};<10.0mm,5.9mm>*{^{\sigma(m)}}**@{},
 <-0.6mm,-0.44mm>*{};<-8mm,-5mm>*{}**@{-},
 <-0.4mm,-0.7mm>*{};<-4.5mm,-5mm>*{}**@{-},
 <0mm,0mm>*{};<1mm,-5mm>*{\ldots}**@{},
 <0.4mm,-0.7mm>*{};<4.5mm,-5mm>*{}**@{-},
 <0.6mm,-0.44mm>*{};<8mm,-5mm>*{}**@{-},
   <0mm,0mm>*{};<-10.5mm,-6.9mm>*{^{\tau(1)}}**@{},
   <0mm,0mm>*{};<-4mm,-6.9mm>*{^{\tau(2)}}**@{},
   <0mm,0mm>*{};<10.0mm,-6.9mm>*{^{\tau(n)}}**@{},
 \end{xy}}\Ea\hspace{-2mm}
=(-1)^{c|\sigma|+d|\tau|}
\Ba{c}\resizebox{14mm}{!}{\begin{xy}
 <0mm,0mm>*{\bu};<0mm,0mm>*{}**@{},
 <-0.6mm,0.44mm>*{};<-8mm,5mm>*{}**@{-},
 <-0.4mm,0.7mm>*{};<-4.5mm,5mm>*{}**@{-},
 <0mm,0mm>*{};<-1mm,5mm>*{\ldots}**@{},
 <0.4mm,0.7mm>*{};<4.5mm,5mm>*{}**@{-},
 <0.6mm,0.44mm>*{};<8mm,5mm>*{}**@{-},
   <0mm,0mm>*{};<-8.5mm,5.5mm>*{^1}**@{},
   <0mm,0mm>*{};<-5mm,5.5mm>*{^2}**@{},
   <0mm,0mm>*{};<4.5mm,5.5mm>*{^{m\hspace{-0.5mm}-\hspace{-0.5mm}1}}**@{},
   <0mm,0mm>*{};<9.0mm,5.5mm>*{^m}**@{},
 <-0.6mm,-0.44mm>*{};<-8mm,-5mm>*{}**@{-},
 <-0.4mm,-0.7mm>*{};<-4.5mm,-5mm>*{}**@{-},
 <0mm,0mm>*{};<-1mm,-5mm>*{\ldots}**@{},
 <0.4mm,-0.7mm>*{};<4.5mm,-5mm>*{}**@{-},
 <0.6mm,-0.44mm>*{};<8mm,-5mm>*{}**@{-},
   <0mm,0mm>*{};<-8.5mm,-6.9mm>*{^1}**@{},
   <0mm,0mm>*{};<-5mm,-6.9mm>*{^2}**@{},
   <0mm,0mm>*{};<4.5mm,-6.9mm>*{^{n\hspace{-0.5mm}-\hspace{-0.5mm}1}}**@{},
   <0mm,0mm>*{};<9.0mm,-6.9mm>*{^n}**@{},
 \end{xy}}\Ea
 \right\rangle_{ \forall \sigma\in \bS_m \atop \forall\tau\in \bS_n}.
$$
 for $m+n\geq 3$, $m\geq 1$, $n\geq 0$. The differential in $\HoqLBcd$ is given on the generators by
 \Beq\label{3: d in qLBcd_infty}
\delta
\Ba{c}\resizebox{14mm}{!}{\begin{xy}
 <0mm,0mm>*{\bu};<0mm,0mm>*{}**@{},
 <-0.6mm,0.44mm>*{};<-8mm,5mm>*{}**@{-},
 <-0.4mm,0.7mm>*{};<-4.5mm,5mm>*{}**@{-},
 <0mm,0mm>*{};<-1mm,5mm>*{\ldots}**@{},
 <0.4mm,0.7mm>*{};<4.5mm,5mm>*{}**@{-},
 <0.6mm,0.44mm>*{};<8mm,5mm>*{}**@{-},
   <0mm,0mm>*{};<-8.5mm,5.5mm>*{^1}**@{},
   <0mm,0mm>*{};<-5mm,5.5mm>*{^2}**@{},
   <0mm,0mm>*{};<4.5mm,5.5mm>*{^{m\hspace{-0.5mm}-\hspace{-0.5mm}1}}**@{},
   <0mm,0mm>*{};<9.0mm,5.5mm>*{^m}**@{},
 <-0.6mm,-0.44
 mm>*{};<-8mm,-5mm>*{}**@{-},
 <-0.4mm,-0.7mm>*{};<-4.5mm,-5mm>*{}**@{-},
 <0mm,0mm>*{};<-1mm,-5mm>*{\ldots}**@{},
 <0.4mm,-0.7mm>*{};<4.5mm,-5mm>*{}**@{-},
 <0.6mm,-0.44mm>*{};<8mm,-5mm>*{}**@{-},
   <0mm,0mm>*{};<-8.5mm,-6.9mm>*{^1}**@{},
   <0mm,0mm>*{};<-5mm,-6.9mm>*{^2}**@{},
   <0mm,0mm>*{};<4.5mm,-6.9mm>*{^{n\hspace{-0.5mm}-\hspace{-0.5mm}1}}**@{},
   <0mm,0mm>*{};<9.0mm,-6.9mm>*{^n}**@{},
 \end{xy}}\Ea
\ \ = \ \
 \sum_{[m]=I_1\sqcup I_2\atop
 {|I_1|\geq 0, |I_2|\geq 1}}
 \sum_{[n]=J_1\sqcup J_2\atop
 {|J_1|, |J_2|\geq 0}
}\hspace{0mm}
\pm
\Ba{c}\resizebox{22mm}{!}{ \begin{xy}
 <0mm,0mm>*{\bu};<0mm,0mm>*{}**@{},
 <-0.6mm,0.44mm>*{};<-8mm,5mm>*{}**@{-},
 <-0.4mm,0.7mm>*{};<-4.5mm,5mm>*{}**@{-},
 <0mm,0mm>*{};<0mm,5mm>*{\ldots}**@{},
 <0.4mm,0.7mm>*{};<4.5mm,5mm>*{}**@{-},
 <0.6mm,0.44mm>*{};<12.4mm,4.8mm>*{}**@{-},
     <0mm,0mm>*{};<-2mm,7mm>*{\overbrace{\ \ \ \ \ \ \ \ \ \ \ \ }}**@{},
     <0mm,0mm>*{};<-2mm,9mm>*{^{I_1}}**@{},
 <-0.6mm,-0.44mm>*{};<-8mm,-5mm>*{}**@{-},
 <-0.4mm,-0.7mm>*{};<-4.5mm,-5mm>*{}**@{-},
 <0mm,0mm>*{};<-1mm,-5mm>*{\ldots}**@{},
 <0.4mm,-0.7mm>*{};<4.5mm,-5mm>*{}**@{-},
 <0.6mm,-0.44mm>*{};<8mm,-5mm>*{}**@{-},
      <0mm,0mm>*{};<0mm,-7mm>*{\underbrace{\ \ \ \ \ \ \ \ \ \ \ \ \ \ \
      }}**@{},
      <0mm,0mm>*{};<0mm,-10.6mm>*{_{J_1}}**@{},
 <13mm,5mm>*{};<13mm,5mm>*{\bu}**@{},
 <12.6mm,5.44mm>*{};<5mm,10mm>*{}**@{-},
 <12.6mm,5.7mm>*{};<8.5mm,10mm>*{}**@{-},
 <13mm,5mm>*{};<13mm,10mm>*{\ldots}**@{},
 <13.4mm,5.7mm>*{};<16.5mm,10mm>*{}**@{-},
 <13.6mm,5.44mm>*{};<20mm,10mm>*{}**@{-},
      <13mm,5mm>*{};<13mm,12mm>*{\overbrace{\ \ \ \ \ \ \ \ \ \ \ \ \ \ }}**@{},
      <13mm,5mm>*{};<13mm,14mm>*{^{I_2}}**@{},
 <12.4mm,4.3mm>*{};<8mm,0mm>*{}**@{-},
 <12.6mm,4.3mm>*{};<12mm,0mm>*{\ldots}**@{},
 <13.4mm,4.5mm>*{};<16.5mm,0mm>*{}**@{-},
 <13.6mm,4.8mm>*{};<20mm,0mm>*{}**@{-},
     <13mm,5mm>*{};<14.3mm,-2mm>*{\underbrace{\ \ \ \ \ \ \ \ \ \ \ }}**@{},
     <13mm,5mm>*{};<14.3mm,-4.5mm>*{_{J_2}}**@{},
 \end{xy}}\Ea
\Eeq
If $c,d\in 2\Z$, all the signs in the r.h.s.\ are equal to $-1$. Taking the quotient of $\HoqLBcd$ by the ideal generated by all $(m,0)$-corollas, $m\geq 3$, gives the minimal model $\HoLBcd$ of the Lie bialgebras properad $\LBcd$.

\subsection{Reminder on the properad of ribbon graphs (after \cite{MW})}\label{4: subsec on RGra} A ribbon graph is a graph in which the set of half-edges attached to each vertex is equipped with a cyclic ordering. More precisely,
a ribbon graph is a triple
$\Ga=(H(\Ga), \tau, \sigma)$ consisting of (i) a finite set  $H(\Ga)$ called the
set of {half-edges}, (ii) a fixed point free involution $\tau: H(\Ga)\rar H(\Ga)$ whose
set of orbits is denoted by $E(\Ga)$ and is called the set of edges, and (iii)
 a permutation $\sigma: H(\Ga)\rar H(\Ga)$ whose set of orbits is denoted by $V(\Ga)$
 and is is called the set of vertices of the ribbon graph $\Ga$. The permutation $\sigma$ can be uniquely represented as the product of cycles, $\sigma=\prod_{v\in V(\Ga)} \sigma_v$; the elements $h_{i_1},\ldots h_{i_v}\in H(\Ga)$  constituting the cycle $\sigma_v=(h_{i_1}\cdots h_{i_v})$  are called half-edges attached to the vertex $v$, they come equipped with a cyclic ordering for each $v\in V(\Ga)$. Next,
the orbits of the permutation $\sigma^{-1}\circ \tau$
are called {\em boundaries}\, of the ribbon graph $\Ga$; they are better understood as real circle-like boundaries of the 2-dimensional surface obtained from $\Ga$ by thickening of edges into 2-dimensional strips. The set of boundaries of $\Ga$ is
 denoted by $B(\Ga)$. For example, the graph
$\hspace{-1mm}
\Ba{c}\resizebox{7mm}{!}{ \xy
(0,5)*{\circ}="1";
(0,-4)*{\circ}="3";
"1";"3" **\crv{(4,0) & (4,1)};
"1";"3" **\crv{(-4,0) & (-4,-1)};
\ar @{-} "1";"3" <0pt>
\endxy}\Ea\hspace{-1mm}
$
has two vertices and three boundaries, while the ribbon graph $\hspace{-1mm}\Ba{c}\resizebox{8mm}{!}{ \xy
(0,5)*{\circ}="1";
(0,-4)*{\circ}="3";
"1";"3" **\crv{(-5,2) & (5,2)};
"1";"3" **\crv{(5,2) & (-5,2)};
"1";"3" **\crv{(-7,7) & (-7,-7)};
\endxy}\Ea\hspace{-1mm}
$ has two vertices and one boundary.

 \sip

  The integer
$$
 g= 1+\frac{1}{2}\left(\# E(\Ga) - \# V(\Ga)- \# B(\Ga)\right)
$$
is called the {\em genus}\, of $\Ga$. This is the genus of a closed surface obtained from $\Ga$ by first thickening edges into strips and then gluing disks into each $S^1$-like boundary of $\Ga$. A ribbon graph $\Ga$ is called {\em directed}\, if its every edge is equipped with the fixed total ordering of its half-edges.

\sip

 Let $RG_{m,n;l}$ be the set of (isomorphism classes of)
directed connected
ribbon graphs $\Ga$  whose edges, vertices and boundaries are distinguished, say, enumerated via some bijections $E(\Ga)\rar [l]$, $V(\Ga) \rar [n]=\{1,2,\ldots, n\}$ and $B(\Ga) \rar [\bar{m}]=\{\bar{1},\ldots, \bar{m}\}$.
The  group $\bS_l \ltimes \bS_l^2\times \bS_m\times \bS_n$ acts on elements of $RG_{m,n;l}$ by relabelling
the edges, vertices and boundaries, and by flipping directions on edges. For any integer $d\in \Z$
we define an $\bS_m^{op}\times \bS_n$-module
$$
\cR\cG ra_{d}(m,n):= \prod_{l\geq 0} \K \langle RG_{m,n:l}\rangle \ot_{ \bS_l \ltimes
 (\bS_2)^l}  \sgn_l^{|d-1|}\ot \sgn_2^{\ot l|d|} [l(d-1)]
$$
Thus a ribbon graph $\Ga$ of $\cR\cG ra_d(m,n)$ has distinguished vertices and distinguished boundaries and
is equipped with an
{\em orientation}\, which depends on the parity of $d$:
\Bi
 \item for even $d$ the orientation of $\Ga$ is defined as the choice
 of ordering of the set of edges $E(\Ga)$ up to a sign action of $\bS_{l}$,
\item for odd $d$ the orientation of $\Ga$ is defined as the choice
of a direction on each edge up to a sign action of $\bS_2$.
\Ei
Note that every ribbon graph $\Ga$ has precisely two possible orientations, $or$ and $or^{opp}$,
and $\Ga$ vanishes if it admits an automorphism which changes the orientation. The group $\bS_m^{op}\times \bS_n$ acts on $\cR\cG ra_{d}(m,n)$ by relabelling the vertices and the boundaries.

\sip

The $\bS$-bimodule $\cR\cG ra_d=\{ \cR\cG ra_{d}(m,n)\}_{m,n\geq 1}$ has a natural structure of a properad described in full details in \cite{MW}; horizontal compositions are given by disjoint unions of ribbon graphs while the vertical compositions are given by substituting  boundaries $b\in B(\Ga_1)$ of one graph into (blown into disks) vertices $v\in V(\Ga_2)$  of another graph and then re-distributing edges attached earlier to $v$ among the vertices belonging to the closed oriented path $b$ in all possible ways while respecting the cyclic orderings.

\sip

The main point of introducing the prop $\cR \cG ra_d$ is that it comes equipped with several remarkable morphisms of properads \cite{MW}, and that these morphisms possess in turn quite interesting  deformation complexes. The first important morphism is the following one,
$$
i: \LB_{d,d} \lon \cR\cG ra_{d}
$$
which is given on generators by,
\Beq\label{4: maps i from LoBdd to RGra}
i:
\Ba{c}\begin{xy}
 <0mm,-0.55mm>*{};<0mm,-2.5mm>*{}**@{-},
 <0.5mm,0.5mm>*{};<2.2mm,2.5mm>*{}**@{-},
 <-0.48mm,0.48mm>*{};<-2.2mm,2.5mm>*{}**@{-},
 <0mm,0mm>*{\bu};<0mm,0mm>*{}**@{},
 <0mm,-0.55mm>*{};<0mm,-3.8mm>*{_{_1}}**@{},
 <0.5mm,0.5mm>*{};<2.7mm,3.1mm>*{^{^{\bar{2}}}}**@{},
 <-0.48mm,0.48mm>*{};<-2.7mm,3.1mm>*{^{^{\bar{1}}}}**@{},
 \end{xy}\Ea
 \lon
 \xy
 (0.5,1)*{^{{^{\bar{1}}}}},
(0.5,5)*{^{{^{\bar{2}}}}},
 (0,-2)*+{_{_1}}*\frm{o}="A";
"A"; "A" **\crv{(7,7) & (-7,7)};
\endxy
 \ \ \  , \ \ \
i:
\Ba{c}\begin{xy}
 <0mm,0.66mm>*{};<0mm,3mm>*{}**@{-},
 <0.39mm,-0.39mm>*{};<2.2mm,-2.2mm>*{}**@{-},
 <-0.35mm,-0.35mm>*{};<-2.2mm,-2.2mm>*{}**@{-},
 <0mm,0mm>*{\bu};<0mm,4.1mm>*{^{^{\bar{1}}}}**@{},
   <0.39mm,-0.39mm>*{};<2.9mm,-4mm>*{^{_2}}**@{},
   <-0.35mm,-0.35mm>*{};<-2.8mm,-4mm>*{^{_1}}**@{},
\end{xy}\Ea
 \lon
 \Ba{c}\resizebox{10mm}{!}{  \xy
 (3.5,4)*{^{\bar{1}}};
 (7,0)*+{_2}*\frm{o}="A";
 (0,0)*+{_1}*\frm{o}="B";
 \ar @{-} "A";"B" <0pt>
\endxy} \Ea.
\Eeq
For $d$ even we assume (implicitly) that  a particular arrow on the unique edge of graphs in the r.h.s.\ of the map $i$ is chosen up to a flip and multiplying by $-1$ (we often omit showing too many details in our pictures such as arrows on edges and  labels of vertices or boundaries). A remarkable fact about the prop $\RGra_d$ is that, due to automorphisms of ribbon graphs from $\RGra_d$ which
 reverse orientations, the graphs controlling Jacobi, co-Jacobi, Drinfeld and involutivity relations in the prop of Lie bialgebras vanish {\em automatically}\, under the map $i$ and hence ensure that the above map $i$ is well-defined. For example, substituting the two labelled boundaries (i.e. the two outputs) of the ribbon graph
 $
\xy
 (0.5,1)*{^{{^{\bar{1}}}}},
(0.5,5)*{^{{^{\bar{2}}}}},
 (0,-2)*+{_{_1}}*\frm{o}="A";
"A"; "A" **\crv{(7,7) & (-7,7)};
\endxy$ into into the respective vertices (i.e. the two inputs) of the graph  $\Ba{c}\resizebox{10mm}{!}{  \xy
 (3.5,4)*{^{\bar{1}}};
 (7,0)*+{_2}*\frm{o}="A";
 (0,0)*+{_1}*\frm{o}="B";
 \ar @{-} "A";"B" <0pt>
\endxy} \Ea$ one obtains a ribbon graph with one vertex and one boundary given by
$$
 \xy
(0,-2)*{\bu}="A";
(0,-2)*{\bu}="B";
"A"; "B" **\crv{(-4,6) & (10,6)};
"A"; "B" **\crv{(-10,6) & (4,6)};
\endxy \in \RGra_d.
$$
However the latter graph admits an automorphism reversing its orientation and hence vanishes identically, i.e.  $\xy
(0,-2)*{\bu}="A";
(0,-2)*{\bu}="B";
"A"; "B" **\crv{(-4,6) & (10,6)};
"A"; "B" **\crv{(-10,6) & (4,6)};
\endxy\equiv 0$  in $\RGra_d$.

\sip

Another important morphism is given by Theorem 4.3.3 in \cite{MW} which says  that there is a properad map
\Beq\label{3: s from HoLB0,1 to RGra}
s: \HoLB_{0,1} \lon \cR\cG ra_{1},
\Eeq
defined on generators by
\begin{align*}
s: \left(
\resizebox{14mm}{!}{\begin{xy}
 <0mm,0mm>*{\bu};<0mm,0mm>*{}**@{},
 <-0.6mm,0.44mm>*{};<-8mm,5mm>*{}**@{-},
 <-0.4mm,0.7mm>*{};<-4.5mm,5mm>*{}**@{-},
 <0mm,0mm>*{};<-1mm,5mm>*{\ldots}**@{},
 <0.4mm,0.7mm>*{};<4.5mm,5mm>*{}**@{-},
 <0.6mm,0.44mm>*{};<8mm,5mm>*{}**@{-},
   <0mm,0mm>*{};<-8.5mm,5.5mm>*{^1}**@{},
   <0mm,0mm>*{};<-5mm,5.5mm>*{^2}**@{},
   <0mm,0mm>*{};<4.5mm,5.5mm>*{^{m\hspace{-0.5mm}-\hspace{-0.5mm}1}}**@{},
   <0mm,0mm>*{};<9.0mm,5.5mm>*{^m}**@{},
 <-0.6mm,-0.44mm>*{};<-8mm,-5mm>*{}**@{-},
 <-0.4mm,-0.7mm>*{};<-4.5mm,-5mm>*{}**@{-},
 <0mm,0mm>*{};<-1mm,-5mm>*{\ldots}**@{},
 <0.4mm,-0.7mm>*{};<4.5mm,-5mm>*{}**@{-},
 <0.6mm,-0.44mm>*{};<8mm,-5mm>*{}**@{-},
   <0mm,0mm>*{};<-8.5mm,-6.9mm>*{^1}**@{},
   <0mm,0mm>*{};<-5mm,-6.9mm>*{^2}**@{},
   <0mm,0mm>*{};<4.5mm,-6.9mm>*{^{n\hspace{-0.5mm}-\hspace{-0.5mm}1}}**@{},
   <0mm,0mm>*{};<9.0mm,-6.9mm>*{^n}**@{},
 \end{xy}}
 \right)=
 \begin{cases}
 \frac 1 m
\sum_{\sigma\in\bS_m}
 \underbrace{
\Ba{c}
{\xy
(5,0)*{...},
   \ar@/^1pc/(0,0)*{\bullet};(10,0)*{\bullet}
   \ar@/^{-1pc}/(0,0)*{\bullet};(10,0)*{\bullet}
   \ar@/^0.6pc/(0,0)*{\bullet};(10,0)*{\bullet}
   \ar@/^{-0.6pc}/(0,0)*{\bullet};(10,0)*{\bullet}
 \endxy}
 \Ea}_{m\ \mathrm{edges}}
&\text{if $n=2$ and $m$ odd}
\\
0 & \text{otherwise}
\end{cases},
\end{align*}
where the symbol $\sum_{\sigma\in \bS_m}$ means the summation over all ways of numbering the $m$
 boundaries on r.h.s. More examples of morphisms and their deformation complexes can be found in \cite{MW}.

 \subsubsection{\bf Remark} There are other constructions of props based on the idea of ribbon graph
 which can be found in the literature (see e.g. \cite{TZ,WW}). That constructions are very different from the prop $\RGra_d$ introduced in \cite{MW} because of the different orientations defined on ribbon graphs (due to different extra structures on that graphs). In particular, all these props  differ in the very definition of what is {\em zero}\, in the associated graded vector spaces of graphs:
 graphs admitting automorphisms  which reverse their orientations must be set identically to zero. This vanishing condition assures the key property of $\RGra_d$ of being a prop under $\LB_{d,d}$ --- the ribbon graphs controlling Jacobi, co-Jacobi and Drinfeld identities vanish {\em automatically}\, in $\RGra_d$ while they stay non-zero in other known constructions of props of ribbon graphs.

 \subsubsection{\bf On graph complexes associated with the properad $\cR\cG ra_d$} The deformation complex of the morphism $i$ in (\ref{4: maps i from LoBdd to RGra}) is given by,
$$
(\RGC_d, \delta+\Delta):=\Def(\LB_{d,d}\stackrel{i}{\lon}\cR\cG ra_d)\simeq
\bigoplus_{m,n\geq 1} \left(\cR\cG ra_{d}(m,n)\ot_{\bS_m\times\bS_n}
(sgn_m^{|d|}\ot
sgn_n^{|d|})\right)[d(2-m-n)],
$$
 It is spanned by ribbon graphs with unidentifiable boundaries and vertices and its differential consists of two operations $\delta$ and $\Delta$ which are governed by the values of $i$ on the Lie and, respectively, coLie generators of $\LBcd$. The Jacobi and co-Jacobi relations in the ideal (\ref{3: R for LieB}) assure that these operations are differentials each,
 $$
 \delta^2=0, \ \ \Delta^2=0,
 $$
  while the V.\ Drinfeld's relation in  (\ref{3: R for LieB}) guarantees their compatibility,
  $$
  \delta\Delta + \Delta\delta=0.
  $$
The differential $\delta$ splits vertices of ribbon graphs $\Ga$ into two new vertices connected by an edge,
while the second  differential $\Delta$ --- introduced first by T.Bridgeland ---
attaches a new edge to $\Ga$  in such a way that the number of boundaries is increased by $1$ but no new loop-like boundaries (i.e.\ the ones consisting of single edges only) are created.
The reduced complex
$$
 (\RGC_d, \delta)
$$
(more precisely, its subcomplex $\RGC_d^{\geq 3}$ spanned by at least trivalent graphs)  is identical to the ribbon graph complex studied earlier by  Robert Penner \cite{Pe} and Maxim Kontsevich \cite{Ko1}; it computes the compactly supported  cohomology groups of moduli spaces of algebraic curves of arbitrary genus with skew-symmetrized marked points,
$$
H^\bu(\RGC_d^{\geq 3}, \delta)=
\prod_{g,m} H^{\bu +2dg-m}_c(\cM_{g,m})\ot_{\bS_m}\sgn_m^{|d|+1}.
$$
The full complex $(\RGC_d, \delta+\Delta)$ computes conjecturally the totality  $\prod_{g\geq 1} H^{\bu}_c(\cM_{g})$
 of moduli spaces of algebraic curves of arbitrary genus without punctures \cite{AWZ}. It was proven in \cite{MW} that there is an explicit morphism of cohomology groups
 $$
 f: H^\bu(\GC_2) \lon H^{\bu+1}\left(\RGC_1, \delta+\Delta\right)
 $$
which is conjecturally an injection; here $\GC_2$ is the famous graph complex introduced by Maxim Kontsevich in \cite{Ko2} in the context of the deformation quantization of Poisson structures (we remind its definition in \S {\ref{4: subseq on reminder about GC_d}} below); it has a remarkable property established by Thomas Willwacher in \cite{W} that
$$
H^0(\GC_2)=\grt_1,
$$
the Lie algebra of the mysterious Grothendieck-Teichm\"uller group $GRT_1$. It has been proven recently by Melody Chan, Soren Galatius and Sam Payne in  \cite{CGP} that there is an injection of the cohomology groups
$$
H^\bu(\GC_2) \lon \prod_{g}H^{\bu+2g}_c(\cM_g;\K).
$$
Conjecturally, this injection coincides with the purely combinatorial morphism of cohomology groups $f$
considered just above.

\sip

The deformation complex of another important morphism (\ref{3: s from HoLB0,1 to RGra}),
$$
(\RGC_{0,1}, \delta+\sum_{i\geq 1}\Delta_{2i+1}):=\Def(\HoLB_{0,1}\stackrel{s}{\lon}\cR gra_1)\simeq
{\bigoplus_{m,n\geq 1\atop} \left(\cR\cG ra_{1}(m,n)\ot_{\bS_m\times\bS_n}
(\id_m\ot
sgn_n)\right)[(1-n)]}
$$
has a differential which is also a deformation of the standard differential $\delta$, but this time the deformed part is an infinite series $\sum_{i\geq 1}\Delta_{2i+1}$ of operators, each increasing the number of boundaries by $2i$ as $i\rar \infty$.
The cohomology of this complex comes equipped with a morphism from the cohomology of  M.\ Kontsevich's {\it odd}\, graph complex,
 $$
 H^\bu(\GC_1) \lon H^{\bu+1}\left(\RGC_{0,1}, \delta+\sum_{i\geq 1}\Delta_{2i+1}\right)
 $$
which is non-trivial on infinitely many polytope classes. Probably this complex $\RGC_{0,1}$ from \cite{MW} and  the complex (\ref{1: totality of Mg,n+m as Def}) with 3-terms differential introduced in this paper are related to each other.

\subsection{Twisted properad of ribbon graphs}\label{3: subsec on twRgra_d} Let $\tw$ be 
Thomas Willwacher's twisting endofunctor \cite{W} in the category of dg properads under 
$\Lie_d$, and
let $\tw\cR\cG ra_d=\{\tw\cR\cG ra_d(m,n)\}_{m\geq 1, n\geq 0}$ be the twisted version of the 
above properad $\cR\cG ra_d$ viewed just as a properad under $\Lie_d$. The properad 
$\tw\cR\cG ra_d$ is generated by elements of $\cR \cG ra_d$ and one extra MC generator with one output (i.e.\ one boundary) and no inputs which we represent as a ribbon graph consisting of just one black vertex $\bu$ of cohomological 
degree $d$ (cf.\ \S {\ref{2: subsec on grav operad and  Rtrees}}); it is not allowed (by its very definition) to substitute
boundaries of other graphs into that vertex, but it is allowed to substitute its unique boundary  into the vertices  of ribbon graphs from $\cR\cG ra_d$. Thus $\tw\cR\cG ra_d(m,n)$ is generated by  connected  ribbon graphs with $m$ labelled boundaries, $n$ labelled white vertices and any number of black vertices which are assigned the cohomological degree $d$, for example
$$
 \Ba{c}\resizebox{9mm}{!}{ \xy
  (0,5)*+{_1}*\frm{o}="1";
(0,-4)*{\bu}="3";
"1";"3" **\crv{(4,0) & (4,1)};
"1";"3" **\crv{(-4,0) & (-4,-1)};
\ar @{-} "1";"3" <0pt>
\endxy}\Ea \in \tw\cR\cG ra_d(3,1),\  \ \Ba{c}\resizebox{9mm}{!}{ \xy
 (0,8)*+{_1}*\frm{o}="1";
(0,-4)*{\bu}="3";
"1";"3" **\crv{(-5,2) & (5,2)};
"1";"3" **\crv{(5,2) & (-5,2)};
"1";"3" **\crv{(-7,7) & (-7,-7)};
\endxy}\Ea \in \tw\cR\cG ra_d(1,1),
 \ \ \
 \Ba{c}\resizebox{10mm}{!}{  \xy
 (3.5,4)*{^{\bar{1}}};
 (7,0)*+{_2}*\frm{o}="A";
 (0,0)*+{_1}*\frm{o}="B";
 \ar @{-} "A";"B" <0pt>
\endxy} \Ea \in \tw\cR\cG ra_d(1,2),
\ \
 \xy
 (0.5,1)*{^{{^{\bar{1}}}}},
(0.5,5)*{^{{^{\bar{2}}}}},
 (0,-2)*+{_{_1}}*\frm{o}="A";
"A"; "A" **\crv{(7,7) & (-7,7)};
\endxy\in \tw\cR\cG ra_d(2,1),
%
%
$$
the first two graphs having degree $3(1-d)+d=3-2d$, the last two graph having degree $1-d$.
We denote the set of white (resp.\ black)  vertices of a ribbon graph $\Ga\in \tw\cR\cG ra_d$ by $V_\circ(\Ga)$ (resp.\ $V_\bu(\Ga)$).
Note that all boundaries of these graphs are labelled but we omit often their labellings in our pictures. The orientation of graphs $\Ga$ from $\tw\RGra_d$ is defined as a choice of
one of the two unit vectors in the 1-dimensional real vector space $Or_d(\Ga)$ (equipped with the standard Euclidean metric) which depends on the parity of $d$,
\Bi
 \item for even $d$, $Or_d(\Ga):=\det(E(\Ga))$,
\item for odd $d$, $Or_d(\Ga):=\det(V_\bu(\Ga))\bigotimes_{e\in E(\Ga)} \det(H(e))$ where
$H(e)$ is the set of two half-edges constituting the edge $e$.
\Ei
Here, for a finite set $I$, we denote $\det(I):=\wedge^{\# I}(\R[I])$ with $\R[I]$ being the linear span of $I$ over the field $\R$ of real numbers.

\sip

 The properadic compositions in $\tw\cR\cG ra_d$ are given by substituting boundaries $b\in B(\Ga')$ into white vertices $v\in V_\circ(\Ga'')$ using exactly same rules as in the case of $\cG\cR ra_d$; it is worth noting that, when distributing dangling half-edges attached to $v$ to the vertices belonging to the boundary $b$, we allow their reconnections to vertices of both types, white and black ones.

\sip

 Note that $\cR\cG ra_d$ has, by its very definition, the trivial differential so that  the induced differential $\sd_\centerdot$ in $\tw\cR\cG ra_d$ comes solely from the twisting procedure and is given,  on an arbitrary ribbon graph $\Ga\in \tw\cR\cG ra_d(m,n)$, by the formula\footnote{Given a prop $\cP=\{\cP(m,n)\}$, we denote by $_i\hspace{-0.0mm}\circ_j$ the partial vertical composition   $\cP(m_1,n_1)\ot \cP(m_2,n_2)\rar \cP(m_1+m_2-1, n_1+n_2-1)$ given, for any $i\in [n_1]$ and $j\in [m_2]$, by inserting $j$-th output of $\cP(m_2,n_2)$ into $i$-th input of $\cP(m_1,n_1)$.}
\Beq\label{4: d in Tw(RGra)}
d_\centerdot\Ga:=
\sum_{i=1}^m \Ba{c}\resizebox{3mm}{!}{  \xy
 (0,6)*{\bu}="A";
 (0,0)*+{_1}*\frm{o}="B";
 \ar @{-} "A";"B" <0pt>
\endxy} \Ea
\  _1\hspace{-0.7mm}\circ_j \Ga\ \
- \ \ (-1)^{|\Ga|} 
\sum_{j=1}^n\Ga\  _j\hspace{-0.7mm}\circ_1
\Ba{c}\resizebox{3mm}{!}{  \xy
 (0,6)*{\bu}="A";
 (0,0)*+{_1}*\frm{o}="B";
 \ar @{-} "A";"B" <0pt>
\endxy} \Ea
\  -(-1)^{|\Ga|}\ \frac{1}{2} \sum_{v\in V_\bu(\Ga)} \Ga\circ_v  \left(\xy
 (0,0)*{\bullet}="a",
(5,0)*{\bu}="b",
\ar @{-} "a";"b" <0pt>
\endxy\right)
\Eeq
where the symbol $\Ga\circ_v  \left(\xy
 (0,0)*{\bullet}="a",
(5,0)*{\bu}="b",
\ar @{-} "a";"b" <0pt>
\endxy\right)$ means substituting the graph  $\xy
 (0,0)*{\bullet}="a",
(5,0)*{\bu}="b",
\ar @{-} "a";"b" <0pt>
\endxy$ into the black vertex $v$ of the graph $\Ga$ and then taking the sum over all possible re-attachments of the half-edges attached earlier  to $v$ among the the two black vertices in a way which respects their cyclic ordering. Note that for almost all graphs the new univalent black vertices arising in the first part of $\sd_\centerdot$ cancel out the new univalent black vertices arising in the second and the third part of that differential.

\subsection{On the role of the integer parameter $d$}\label{3: subsec on orientations for different d}
It is clear from the very definition that dg props $\tw\cR\cG ra_d$ with $d$ of the  same
{\it parity}\, are isomorphic to each other (up to degree shifts). However it is not immediately evident that dg props $\tw\cR\cG ra_d$ for $d$ with opposite parity control essentially one and the same mathematics because the definitions of {\it orientation}\ of ribbon graphs for $d$ even and $d$ odd look quite different.

\sip

Let $\Ga$ be any connected ribbon graph
 and $\Sigma_\Ga$ be the associated  connected compact oriented Riemann surface canonically associated to by $\Gamma$ by thickening edges of $\Ga$ into strips and filling every boundary with a disk. The ribbon graph gives us a cell decomposition of $\Sigma_\Gamma$ and hence the associated cell complex,
 $$
 C_2(\Sigma_\Gamma)\stackrel{\p_2}{\lon} C_1(\Sigma_\Gamma) \stackrel{\p_1}{\lon} C_0(\Sigma_\Gamma)
 $$
 computing the singular homology $H_\bu(\Sigma_\Ga)$ with coefficients in the field $\K$. We can identify $C_2(\Sigma_\Gamma)$ with the vector space $\K[B(\Ga)]$ generated by the set of boundaries of $\Gamma$, $C_0(\Sigma_\Gamma)$ with the vector space $\K[V(\Ga)]$ generated by its set of vertices, and
 $$
 C_1(\Sigma_\Gamma)=\bigoplus_{e\in E(\Gamma)} \K[\det( H(e))]
 $$
Then the three  exact sequences,
$$
0\lon H_2(\Sigma_\Gamma) \lon C_2(\Sigma_\Gamma) \lon \Img\p_2 \lon 0
$$
$$
0\lon \Img\p_2 \lon \Ker \p_1 \lon H_1(\Sigma_\Gamma) \lon 0
$$
$$
0\lon \Ker\p_1 \lon  C_1(\Sigma_\Gamma) \lon C_0(\Sigma_\Gamma) \lon H_0(S_\Gamma)\lon 0
$$
imply the following canonical isomorphisms
$$
\det(B(\Ga))=\det(H_2(\Sigma_\Gamma))\ot \det(\Img\p_2),\ \ \ \ \ \ \   \det(\Ker\p_1)=\det(\Img\p_2)\ot \det(H_1(\Sigma_\Ga)),
$$
$$
\det(\Ker\p_1)\ot \det(V(\Ga))=\det (E(\Ga))\ot \bigotimes_{e\in E(\Ga)}\det(H(e))\ot \det(H_0(S_\Ga)).
$$
As spaces $H_2(\Sigma_\Gamma)$ and $H_0(\Sigma_\Gamma)$ are 1-dimensional and come equipped with canonical bases, they can be omitted in the above equalities. The vector space $H_1(\Sigma_\Gamma)$ can be equipped with a so called {\em canonical basis of cycles}\,
$\{a_i,b_i\}_{i\in [g]}$ characterized by the standard intersection properties, $a_i\cdot a_j=0$, $b_i\cdot b_j=0$, $a_i\cdot b_j=\delta_{ij}$. This basis is not defined uniquely, but up to a multiplication by a matrix
from the symplectic group $Sp({2g},\Z)$. Therefore such a non-unique canonical basis of cycles gives us a well-defined distinguished basis vector of $\det(H_1(\Sigma_\Ga))$
so that this tensor factor can also be omitted in the above equalities. The remaining equalities assemble into the canonical isomorphism
$$
\det (E(\Ga))\ot \det(V(\Ga))\ot \det(B(\Ga))\equiv   \bigotimes_{e\in E(\Ga)}\det(H(e)).
$$
or, equivalently,
$$
\det (E(\Ga))\ot \det(V_\circ(\Ga))\ot \det(V_\bu(\Ga))\ot \det(B(\Ga))\equiv   \bigotimes_{e\in E(\Ga)}\det(H(e)).
$$
As boundaries and white vertices of $\Ga$ are labelled, the vector spaces $\det(V_\circ(\Ga))$ and $\det(B(\Ga))$ can be identified canonically with $\K$ so that the above identity simplifies to the following one,
$$
\det(E(\Ga))=\det(V_\bu(\Ga)) \bigotimes_{e\in E(\Ga)}\det(H(e))
$$
which identifies orientation spaces $Or_d$ for any $d$. Thus we have proved the following

\subsubsection{\bf Lemma} {\em Up to degree shifts, the dg properads $\tw\RGra_d$ for different $d$ are isomorphic to each other.}

\mip

Let $\tw\RGra_d(m,n;g)\subset \tw\RGra_d(m,n)$ be the subcomplex spanned by twisted ribbon graphs of genus $g$. Then above mentioned degree shift is given explicitly by
\Beq\label{3: degree shift iso for twGra_d}
\tw\RGra_d(m,n;g)=\tw\RGra_0(m,n;g)[d(2g-2+m+n)],  \ \ \  \forall d\in \Z,
\Eeq
%
%
that is, an oriented genus g ribbon graph $\Ga$ of cohomological degree $k$ in $\tw\RGra_0(m,n)$ will have degree
$k-d(2g-2+m+n)$ when considered as an oriented graph in $\tw\RGra_d(m,n)$.

\sip

We conclude that it is enough to understand the meaning of the dg prop $\tw\cR\cG ra_0$ (though in applications, especially in string topology \cite{Me}, it is often useful to work with its degree shifted version $\tw\cR\cG ra_d$ for a suitable $d$, as $d$ is determined by the dimension of a manifold under study).

\subsection{The (chain) gravity properad}
It is easy to see that the subspace $\tw\cR\cG ra_d^{\geq 3}$  of    $\tw\cR\cG ra_d$ spanned by ribbon graphs {\it with every black vertex at least trivalent}\, is a dg sub-properad which is called the {\it chain gravity properad}, while its cohomology
$$
\GRav_d:= H^\bu\left( \tw\cR\cG ra_d^{\geq 3} \right)
$$
is called the {\it gravity properad}\, (for the reasons explained below).

\sip

The properad $\GRav_0$ is abbreviated by $\GRav$.

\subsubsection{\bf Proposition}\label{3: prop on tw^3Rgra to twRgra} {\it The inclusion
$$
\tw\cR\cG ra_d^{\geq 3}(m,n) \lon \tw\cR\cG ra_d(m,n)
$$
is  a quasi-isomorphism for any $m,n\geq 1$. For $n=0$ on has
$$
H^\bu(\tw\cR\cG ra_d(m,0))= H^\bu\left(\tw\cR\cG ra_d^{\geq 3}(m,0)\right) \oplus \ \bigoplus_{p\geq 1\atop
p\equiv 2d+1 \bmod 4} \K[-p(d-1)]
$$
where the summand $\K[-p]$ is generated by the polytope with $p$ edges and $p$ bivalent vertices which are all black.
}

\begin{proof} The argument is the same as in the case of ordinary (i.e.\ non-ribbon) graph complexes. We shall sketch it  following \cite{W}.

\sip

Let  $\tw\cR\cG ra_d(m,n)^1 \subset \tw\cR\cG ra_d(m,n)$ be the subcomplex spanned by ribbon graphs with at least one univalent black vertex. The differential preserves the total number of univalent black vertices so we can assume without loss of generality that they are distinguished, say, enumerated.
Call {\it antenna}\, the maximal subgraph of a ribbon graph from  $\tw\cR\cG ra_d(m,n)^1$ consisting of the univalent black vertex labelled by the minimal integer and the maximal connected chain of bivalent black vertices attached to it; denote such a chain with $n\geq 0$ bivalent black vertices by $A_n$.
Consider a filtration of $\tw\cR\cG ra_d(m,n)^1$ by the total number of non-antenna vertices. The induced differential acts only on $A_n$ antennas as follows,
$$
d: A_{2n}\rar A_{2n+1}, \ \  d: A_{2n+1}\rar 0.
$$
Hence the associated graded complex is acyclic implying acyclicity of  $\tw\cR\cG ra_d(m,n)^1$.

\sip

Let $\tw\cR\cG ra_d(m,n)^{2}\subset  \tw\cR\cG ra_d(m,n)$ be the subcomplex spanned by ribbon graphs with at least one bivalent black vertex. It contains in turn a subcomplex spanned by graphs with at least one trivalent black vertex or at least one white vertex. It is easy to see that this subcomplex is acyclic (see the proof of Proposition 3.4 in \cite{W}). It remains to consider a subcomplex of $\tw\cR\cG ra_d(m,n)$
spanned by ribbon graphs whose all vertices are bivalent and black. Its cohomology is given by the polytopes described in the Proposition.
\end{proof}


\subsection{What does the properad $\GRav$ measure and control?} An answer about the part, $\GRav(m,0)$, of the  gravity properad with no white vertices  is immediate --- the complex  $\tw\cR\cG ra^{\geq 3}(m,0)$ is precisely the (degree shifted) R.\ Penner's ribbon graph complex \cite{Pe,Ko1} with marked boundaries, so that its cohomology  is given by
$$
\GRav(m,0)=\prod_{2g+m\geq 3} H^{\bu-m}_c(\cM_{g,m})
$$
where $\cM_{g,m}$ is the moduli spaces of genus $g$ algebraic curves with $m$ marked points, and $H^\bu_c$ stands for the compactly supported cohomology functor. Our main purpose in this section is to understand the algebro-geometric meaning of $\GRav(m,n)$ for $n\geq 1$.

\sip

Let $\cM_{g, m+n}$  the moduli space of genus $g$ algebraic curves with the set $S$ of distinct marked  points decomposed into the disjoint union
$$
S=S^{out}\sqcup S_{in},\ \ \ \# S^{out}=m\geq 1,\  \# S_{in}=n\geq 0.
$$
The permutation group $\bS_m^{op}\times \bS_n$ acts naturally on $\cM_{g, m+n}$ by relabelling (separately) the out- and in-points. Note that $\cM_{g,m_1+n_1}$ and $\cM_{g,m_2+n_2}$ with the
same total number of marked points, $n_1+m_1=n_2+m_2$, have to be considered as different topological
$\bS$-bimodules if $m_1\neq m_2$.

\subsubsection{\bf Theorem}\label{3: Main theorem on H(twRgra)} {\em For any $m\geq 1, n\geq 0$ there is an isomorphism of $\bS_m^{op}\times \bS_n$-modules}.
$$
\GRav(m,n)=\prod_{2g+m+n\geq 3\atop  m\geq 1, n\geq 0}H_c^{\bu -m}(\cM_{g, m+n})
$$

\sip

For $n=0$ this result says nothing new. For $n\geq 1$ it is less obvious;
the proof is based on the theory of K.\ Costello's moduli spaces of nodal disks \cite{Co1,Co2} and his wonderful homotopy equivalence; it is presented in the subsection  \S {\ref{3: subsec on Costello moduli spaces}} below.

\sip

The chain gravity properad contains the dg operad of twisted ribbon trees from \S {\ref{2: subsec on grav operad and  Rtrees}} as a dg sub-operad. Thus the gravity operad is a sub-operad of $\GRav$.

\subsection{Quasi-Lie bialgebra sub-properad of the gravity properad} Before proving the above Theorem, let us study three simplest non-trivial cohomology classes of the chain gravity properad, and their inter-relations.

\sip

\sip

It is easy to check that the graph  $\Ba{c}\resizebox{10mm}{!}{  \xy
 (7,0)*+{_2}*\frm{o}="A";
 (0,0)*+{_1}*\frm{o}="B";
 \ar @{-} "A";"B" <0pt>
\endxy}\Ea$ is a cycle in $\tw\cR\cG ra_0^{\geq 3}$ so that there is a morphism of dg properads
$$
\Ba{rccc}
\imath: & \Lie_0 & \lon & \tw\cR\cG ra_0  \vspace{1mm}\\
        &
\Ba{c}\begin{xy}
 <0mm,0.66mm>*{};<0mm,3mm>*{}**@{-},
 <0.39mm,-0.39mm>*{};<2.2mm,-2.2mm>*{}**@{-},
 <-0.35mm,-0.35mm>*{};<-2.2mm,-2.2mm>*{}**@{-},
 <0mm,0mm>*{\bu};
   <0.39mm,-0.39mm>*{};<2.9mm,-4mm>*{^{_2}}**@{},
   <-0.35mm,-0.35mm>*{};<-2.8mm,-4mm>*{^{_1}}**@{},
\end{xy}\Ea
& \lon &
 \Ba{c}\resizebox{10mm}{!}{  \xy
 (7,0)*+{_2}*\frm{o}="A";
 (0,0)*+{_1}*\frm{o}="B";
 \ar @{-} "A";"B" <0pt>
\endxy} \Ea
\Ea
$$
given by the same formula as in the case of $\cR\cG ra_0$. Moreover, this class is not a coboundary so that it induces a non-trivial map of properads
$$
\Lie_0 \lon \GRav.
$$
Another important graph
$ \xy
 (0.5,1)*{^{{^{\bar{1}}}}},
(0.5,5)*{^{{^{\bar{2}}}}},
 (0,-2)*+{_{_1}}*\frm{o}="A";
"A"; "A" **\crv{(7,7) & (-7,7)};
\endxy$ in $\tw\RGra_0$  is {\it not}\, a cocycle,
$$
\sd_\cdot \xy
 (0.5,1)*{^{{^{\bar{1}}}}},
(0.5,5)*{^{{^{\bar{2}}}}},
 (0,-2)*+{_{_1}}*\frm{o}="A";
"A"; "A" **\crv{(7,7) & (-7,7)};
\endxy=
\Ba{c}\resizebox{6mm}{!}{
\mbox{$\xy
 (0.5,0.9)*{^{{^{\bar{1}}}}},
(0.5,5)*{^{{^{\bar{2}}}}},
 (0,-8)*+{_{_1}}*\frm{o}="C";
(0,-2)*{\bu}="A";
(0,-2)*{\bu}="B";
"A"; "B" **\crv{(6,6) & (-6,6)};
 \ar @{-} "A";"C" <0pt>
\endxy$}}
\Ea
+
\Ba{c}\resizebox{6mm}{!}{
\mbox{$\xy
 (0.5,0.9)*{^{{^{\bar{2}}}}},
(0.5,5)*{^{{^{\bar{1}}}}},
 (0,-8)*+{_{_1}}*\frm{o}="C";
(0,-2)*{\bu}="A";
(0,-2)*{\bu}="B";
"A"; "B" **\crv{(6,6) & (-6,6)};
 \ar @{-} "A";"C" <0pt>
\endxy$}}
\Ea
$$
so that the above map $\imath$ does {\it not}\, extend to a morphism
$\LB_{0,0}\rar \tw\cR\cG ra_0$. It is easy to see that the graph
$$
\Ba{c}\resizebox{6mm}{!}{
\mbox{$\xy
 (0.5,0.9)*{^{{^{\bar{1}}}}},
(0.5,5)*{^{{^{\bar{2}}}}},
 (0,-8)*+{_{_1}}*\frm{o}="C";
(0,-2)*{\bu}="A";
(0,-2)*{\bu}="B";
"A"; "B" **\crv{(6,6) & (-6,6)};
 \ar @{-} "A";"C" <0pt>
\endxy$}}
\Ea
$$
is a cocycle in $\tw\cR\cG ra_0^{\geq 3}$ which is not a coboundary. Moreover, at the cohomology
level it defines a cohomology class in $\GRav(2,1)$ of  degree $2$ which is, by the above formula, skew-symmetric with respect to the out-labels $\bar{1}$ and $\bar{2}$.

\subsubsection{\bf Proposition}\label{3: Prop on map j} {\it There is a morphism of properads
$$
j:q\LB_{-1,0} \lon \GRav
$$
given on the generators of $\LB_{-1,0}$ by}
\Beq\label{3: map j from qLB to H(TwGra)}
j:  \
\Ba{c}\begin{xy}
<0mm,-0.55mm>*{};<0mm,-2.5mm>*{}**@{-},
 <0.5mm,0.5mm>*{};<2.2mm,2.5mm>*{}**@{-},
 <-0.48mm,0.48mm>*{};<-2.2mm,2.5mm>*{}**@{-},
 <0mm,0mm>*{\bu};<0mm,0mm>*{}**@{},
 <0mm,-0.55mm>*{};<0mm,-3.8mm>*{_{_1}}**@{},
 <0.5mm,0.5mm>*{};<2.7mm,3.1mm>*{^{^{\bar{2}}}}**@{},
 <-0.48mm,0.48mm>*{};<-2.7mm,3.1mm>*{^{^{\bar{1}}}}**@{},
 \end{xy}\Ea
 \lon  \left[\Ba{c}\resizebox{6mm}{!}{
\mbox{$\xy
 (0.5,0.9)*{^{{^{\bar{1}}}}},
(0.5,5)*{^{{^{\bar{2}}}}},
 (0,-8)*+{_{_1}}*\frm{o}="C";
(0,-2)*{\bu}="A";
(0,-2)*{\bu}="B";
"A"; "B" **\crv{(6,6) & (-6,6)};
 \ar @{-} "A";"C" <0pt>
\endxy$}}
\Ea\right], \ \ \ \ \
\Ba{c}\begin{xy}
 <0mm,0.66mm>*{};<0mm,3mm>*{}**@{-},
 <0.39mm,-0.39mm>*{};<2.2mm,-2.2mm>*{}**@{-},
 <-0.35mm,-0.35mm>*{};<-2.2mm,-2.2mm>*{}**@{-},
 <0mm,0mm>*{\bu};<0mm,4.1mm>*{^{^{\bar{1}}}}**@{},
   <0.39mm,-0.39mm>*{};<2.9mm,-4mm>*{^{_2}}**@{},
   <-0.35mm,-0.35mm>*{};<-2.8mm,-4mm>*{^{_1}}**@{},
\end{xy}\Ea
 \lon
 \Ba{c}\resizebox{10mm}{!}{  \xy
 (3.5,4)*{^{\bar{1}}};
 (7,0)*+{_2}*\frm{o}="A";
 (0,0)*+{_1}*\frm{o}="B";
 \ar @{-} "A";"B" <0pt>
\endxy} \Ea,
\ \ \ \ \
\Ba{c}\begin{xy}
 <0mm,-1mm>*{\bu};<-4mm,3mm>*{^{_1}}**@{-},
 <0mm,-1mm>*{\bu};<0mm,3mm>*{^{_2}}**@{-},
 <0mm,-1mm>*{\bu};<4mm,3mm>*{^{_3}}**@{-},
 \end{xy}\Ea
 \lon
 \frac{1}{2}\left(-\hspace{-1mm}
 \Ba{c}\resizebox{11mm}{!}{
\xy
(2.0,-3.5)*{^{{^{\bar{2}}}}},
 (-2,-3.5)*{^{{^{\bar{1}}}}},
(0.5,4)*{^{{^{\bar{3}}}}},
 (0,-8)*{\bu}="C";
(0,3)*{\bu}="A1";
(0,3)*{\bu}="A2";
"C"; "A1" **\crv{(-5,-9) & (-5,4)};
"C"; "A2" **\crv{(5,-9) & (5,4)};
 \ar @{-} "A1";"C" <0pt>
\endxy}
\Ea
+\hspace{-1mm}
\Ba{c}\resizebox{11mm}{!}{
\xy
(2.0,-3.5)*{^{{^{\bar{1}}}}},
 (-2,-3.5)*{^{{^{\bar{2}}}}},
(0.5,4)*{^{{^{\bar{3}}}}},
 (0,-8)*{\bu}="C";
(0,3)*{\bu}="A1";
(0,3)*{\bu}="A2";
"C"; "A1" **\crv{(-5,-9) & (-5,4)};
"C"; "A2" **\crv{(5,-9) & (5,4)};
 \ar @{-} "A1";"C" <0pt>
\endxy}
\Ea\right)
\Eeq

\begin{proof} Let us lift the above values of the map $j$ to the gravity chain complex
by the same formulae as above except the first one
$$
j: \ \ \
\Ba{c}\begin{xy}
<0mm,-0.55mm>*{};<0mm,-2.5mm>*{}**@{-},
 <0.5mm,0.5mm>*{};<2.2mm,2.5mm>*{}**@{-},
 <-0.48mm,0.48mm>*{};<-2.2mm,2.5mm>*{}**@{-},
 <0mm,0mm>*{\bu};<0mm,0mm>*{}**@{},
 <0mm,-0.55mm>*{};<0mm,-3.8mm>*{_{_1}}**@{},
 <0.5mm,0.5mm>*{};<2.7mm,3.1mm>*{^{^{\bar{2}}}}**@{},
 <-0.48mm,0.48mm>*{};<-2.7mm,3.1mm>*{^{^{\bar{1}}}}**@{},
 \end{xy}\Ea
 \lon  \frac{1}{2}\left(
 \Ba{c}\resizebox{6mm}{!}{
\mbox{$\xy
 (0.5,0.9)*{^{{^{\bar{1}}}}},
(0.5,5)*{^{{^{\bar{2}}}}},
 (0,-8)*+{_{_1}}*\frm{o}="C";
(0,-2)*{\bu}="A";
(0,-2)*{\bu}="B";
"A"; "B" **\crv{(6,6) & (-6,6)};
 \ar @{-} "A";"C" <0pt>
\endxy$}}
\Ea
-
 \Ba{c}\resizebox{6mm}{!}{
\mbox{$\xy
 (0.5,0.9)*{^{{^{\bar{2}}}}},
(0.5,5)*{^{{^{\bar{1}}}}},
 (0,-8)*+{_{_1}}*\frm{o}="C";
(0,-2)*{\bu}="A";
(0,-2)*{\bu}="B";
"A"; "B" **\crv{(6,6) & (-6,6)};
 \ar @{-} "A";"C" <0pt>
\endxy$}}
\Ea
\right)=: \Ba{c}\begin{xy}
 <0mm,-0.55mm>*{};<0mm,-2.5mm>*{}**@{-},
 <0.5mm,0.5mm>*{};<2.2mm,2.5mm>*{}**@{-},
 <-0.48mm,0.48mm>*{};<-2.2mm,2.5mm>*{}**@{-},
 <0mm,0mm>*{\circ};<0mm,0mm>*{}**@{},
 <0mm,-0.55mm>*{};<0mm,-3.8mm>*{_{_1}}**@{},
 <0.5mm,0.5mm>*{};<2.7mm,3.1mm>*{^{^{\bar{2}}}}**@{},
 <-0.48mm,0.48mm>*{};<-2.7mm,3.1mm>*{^{^{\bar{1}}}}**@{},
 \end{xy}\Ea
$$
and check that all the relations $\cR_q$ in $q\LBcd$ are respected by the map $j$ modulo $\sd_\centerdot$-exact terms. The element $\Ba{c}\begin{xy}
 <0mm,-0.55mm>*{};<0mm,-2.5mm>*{}**@{-},
 <0.5mm,0.5mm>*{};<2.2mm,2.5mm>*{}**@{-},
 <-0.48mm,0.48mm>*{};<-2.2mm,2.5mm>*{}**@{-},
 <0mm,0mm>*{\circ};<0mm,0mm>*{}**@{},
 <0mm,-0.55mm>*{};<0mm,-3.8mm>*{}**@{},
 <0.5mm,0.5mm>*{};<2.7mm,3.1mm>*{^{^{\bar{2}}}}**@{},
 <-0.48mm,0.48mm>*{};<-2.7mm,3.1mm>*{^{^{\bar{1}}}}**@{},
 \end{xy}\Ea$ is a cycle in $\tw\RGra_0^{\geq 3}$.

 \sip

 We assume tacitly in all our pictures that the orientation of any twisted ribbon graph is given by {\it  ordering the edges from the bottom to the top, and, if on the same level, from left to right}. This rule explains the signs in the formulae below.

\mip

{\sc{Step 1}}: first we  check the Drinfeld compatibility condition, that is, the equality
$$
j\left(\Ba{c}\resizebox{5mm}{!}{\begin{xy}
 <0mm,2.47mm>*{};<0mm,0.12mm>*{}**@{-},
 <0.5mm,3.5mm>*{};<2.2mm,5.2mm>*{}**@{-},
 <-0.48mm,3.48mm>*{};<-2.2mm,5.2mm>*{}**@{-},
 <0mm,3mm>*{\bu};<0mm,3mm>*{}**@{},
  <0mm,-0.8mm>*{\bu};<0mm,-0.8mm>*{}**@{},
<-0.39mm,-1.2mm>*{};<-2.2mm,-3.5mm>*{}**@{-},
 <0.39mm,-1.2mm>*{};<2.2mm,-3.5mm>*{}**@{-},
     <0.5mm,3.5mm>*{};<2.8mm,5.7mm>*{^2}**@{},
     <-0.48mm,3.48mm>*{};<-2.8mm,5.7mm>*{^1}**@{},
   <0mm,-0.8mm>*{};<-2.7mm,-5.2mm>*{^1}**@{},
   <0mm,-0.8mm>*{};<2.7mm,-5.2mm>*{^2}**@{},
\end{xy}}\Ea
-
\Ba{c}\resizebox{7mm}{!}{\begin{xy}
 <0mm,-1.3mm>*{};<0mm,-3.5mm>*{}**@{-},
 <0.38mm,-0.2mm>*{};<2.0mm,2.0mm>*{}**@{-},
 <-0.38mm,-0.2mm>*{};<-2.2mm,2.2mm>*{}**@{-},
<0mm,-0.8mm>*{\bu};<0mm,0.8mm>*{}**@{},
 <2.4mm,2.4mm>*{\bu};<2.4mm,2.4mm>*{}**@{},
 <2.77mm,2.0mm>*{};<4.4mm,-0.8mm>*{}**@{-},
 <2.4mm,3mm>*{};<2.4mm,5.2mm>*{}**@{-},
     <0mm,-1.3mm>*{};<0mm,-5.3mm>*{^1}**@{},
     <2.5mm,2.3mm>*{};<5.1mm,-2.6mm>*{^2}**@{},
    <2.4mm,2.5mm>*{};<2.4mm,5.7mm>*{^2}**@{},
    <-0.38mm,-0.2mm>*{};<-2.8mm,2.5mm>*{^1}**@{},
    \end{xy}}\Ea
 -
\Ba{c}\resizebox{7mm}{!}{\begin{xy}
 <0mm,-1.3mm>*{};<0mm,-3.5mm>*{}**@{-},
 <0.38mm,-0.2mm>*{};<2.0mm,2.0mm>*{}**@{-},
 <-0.38mm,-0.2mm>*{};<-2.2mm,2.2mm>*{}**@{-},
<0mm,-0.8mm>*{\bu};<0mm,0.8mm>*{}**@{},
 <2.4mm,2.4mm>*{\bu};<2.4mm,2.4mm>*{}**@{},
 <2.77mm,2.0mm>*{};<4.4mm,-0.8mm>*{}**@{-},
 <2.4mm,3mm>*{};<2.4mm,5.2mm>*{}**@{-},
     <0mm,-1.3mm>*{};<0mm,-5.3mm>*{^2}**@{},
     <2.5mm,2.3mm>*{};<5.1mm,-2.6mm>*{^1}**@{},
    <2.4mm,2.5mm>*{};<2.4mm,5.7mm>*{^2}**@{},
    <-0.38mm,-0.2mm>*{};<-2.8mm,2.5mm>*{^1}**@{},
    \end{xy}}\Ea
  +
\Ba{c}\resizebox{7mm}{!}{\begin{xy}
 <0mm,-1.3mm>*{};<0mm,-3.5mm>*{}**@{-},
 <0.38mm,-0.2mm>*{};<2.0mm,2.0mm>*{}**@{-},
 <-0.38mm,-0.2mm>*{};<-2.2mm,2.2mm>*{}**@{-},
<0mm,-0.8mm>*{\bu};<0mm,0.8mm>*{}**@{},
 <2.4mm,2.4mm>*{\bu};<2.4mm,2.4mm>*{}**@{},
 <2.77mm,2.0mm>*{};<4.4mm,-0.8mm>*{}**@{-},
 <2.4mm,3mm>*{};<2.4mm,5.2mm>*{}**@{-},
     <0mm,-1.3mm>*{};<0mm,-5.3mm>*{^2}**@{},
     <2.5mm,2.3mm>*{};<5.1mm,-2.6mm>*{^1}**@{},
    <2.4mm,2.5mm>*{};<2.4mm,5.7mm>*{^1}**@{},
    <-0.38mm,-0.2mm>*{};<-2.8mm,2.5mm>*{^2}**@{},
    \end{xy}}\Ea
 +
\Ba{c}\resizebox{7mm}{!}{\begin{xy}
 <0mm,-1.3mm>*{};<0mm,-3.5mm>*{}**@{-},
 <0.38mm,-0.2mm>*{};<2.0mm,2.0mm>*{}**@{-},
 <-0.38mm,-0.2mm>*{};<-2.2mm,2.2mm>*{}**@{-},
<0mm,-0.8mm>*{\bu};<0mm,0.8mm>*{}**@{},
 <2.4mm,2.4mm>*{\bu};<2.4mm,2.4mm>*{}**@{},
 <2.77mm,2.0mm>*{};<4.4mm,-0.8mm>*{}**@{-},
 <2.4mm,3mm>*{};<2.4mm,5.2mm>*{}**@{-},
     <0mm,-1.3mm>*{};<0mm,-5.3mm>*{^1}**@{},
     <2.5mm,2.3mm>*{};<5.1mm,-2.6mm>*{^2}**@{},
    <2.4mm,2.5mm>*{};<2.4mm,5.7mm>*{^1}**@{},
    <-0.38mm,-0.2mm>*{};<-2.8mm,2.5mm>*{^2}**@{},
    \end{xy}}\Ea
\right)=\sd_\centerdot(\text{something})
$$
Remarkably, it holds true at the level of cocycles, i.e.\ we can write zero on the r.h.s. The rules for the properadic composition in $\tw\RGra_0^{\geq 3}$ give us  the following formulae,
$$
\Ba{c}\resizebox{7mm}{!}{\begin{xy}
 <0mm,2.47mm>*{};<0mm,0.12mm>*{}**@{-},
 <0.5mm,3.5mm>*{};<2.2mm,5.2mm>*{}**@{-},
 <-0.48mm,3.48mm>*{};<-2.2mm,5.2mm>*{}**@{-},
 <0mm,3mm>*{\circ};<0mm,3mm>*{}**@{},
  <0mm,-0.8mm>*{\circ};<0mm,-0.8mm>*{}**@{},
<-0.39mm,-1.2mm>*{};<-2.2mm,-3.5mm>*{}**@{-},
 <0.39mm,-1.2mm>*{};<2.2mm,-3.5mm>*{}**@{-},
     <0.5mm,3.5mm>*{};<2.8mm,5.7mm>*{^{\bar{2}}}**@{},
     <-0.48mm,3.48mm>*{};<-2.8mm,5.7mm>*{^{\bar{1}}}**@{},
   <0mm,-0.8mm>*{};<-2.7mm,-5.2mm>*{^1}**@{},
   <0mm,-0.8mm>*{};<2.7mm,-5.2mm>*{^2}**@{},
\end{xy}}\Ea=
\Ba{c}\resizebox{6mm}{!}{
\mbox{$\xy
 (0.5,0.9)*{^{{^{\bar{1}}}}},
(0.5,5)*{^{{^{\bar{2}}}}},
 (0,-7)*+{_{_2}}*\frm{o}="C";
  (0,-14)*+{_{_1}}*\frm{o}="D";
(0,-2)*{\bu}="A";
(0,-2)*{\bu}="B";
"A"; "B" **\crv{(6,6) & (-6,6)};
 \ar @{-} "A";"C" <0pt>
 \ar @{-} "D";"C" <0pt>
\endxy$}}
\Ea
+
\Ba{c}\resizebox{6mm}{!}{
\mbox{$\xy
 (0.5,0.9)*{^{{^{\bar{1}}}}},
(0.5,5)*{^{{^{\bar{2}}}}},
 (0,-7)*+{_{_1}}*\frm{o}="C";
  (0,-14)*+{_{_2}}*\frm{o}="D";
(0,-2)*{\bu}="A";
(0,-2)*{\bu}="B";
"A"; "B" **\crv{(6,6) & (-6,6)};
 \ar @{-} "A";"C" <0pt>
 \ar @{-} "D";"C" <0pt>
\endxy$}}
\Ea
-
\Ba{c}\resizebox{6mm}{!}{
\mbox{$\xy
 (0.5,0.9)*{^{{^{\bar{2}}}}},
(0.5,5)*{^{{^{\bar{1}}}}},
 (0,-7)*+{_{_2}}*\frm{o}="C";
  (0,-14)*+{_{_1}}*\frm{o}="D";
(0,-2)*{\bu}="A";
(0,-2)*{\bu}="B";
"A"; "B" **\crv{(6,6) & (-6,6)};
 \ar @{-} "A";"C" <0pt>
 \ar @{-} "D";"C" <0pt>
\endxy$}}
\Ea
-
\Ba{c}\resizebox{6mm}{!}{
\mbox{$\xy
 (0.5,0.9)*{^{{^{\bar{2}}}}},
(0.5,5)*{^{{^{\bar{1}}}}},
 (0,-7)*+{_{_2}}*\frm{o}="C";
  (0,-14)*+{_{_1}}*\frm{o}="D";
(0,-2)*{\bu}="A";
(0,-2)*{\bu}="B";
"A"; "B" **\crv{(6,6) & (-6,6)};
 \ar @{-} "A";"C" <0pt>
 \ar @{-} "D";"C" <0pt>
\endxy$}}
\Ea
$$
and
$$
\Ba{c}\resizebox{10.0mm}{!}{\begin{xy}
 <0mm,-1.3mm>*{};<0mm,-3.5mm>*{}**@{-},
 <0.38mm,-0.2mm>*{};<2.0mm,2.0mm>*{}**@{-},
 <-0.38mm,-0.2mm>*{};<-2.2mm,2.2mm>*{}**@{-},
<0mm,-0.8mm>*{\circ};<0mm,0.8mm>*{}**@{},
 <2.4mm,2.4mm>*{\circ};<2.4mm,2.4mm>*{}**@{},
 <2.77mm,2.0mm>*{};<4.4mm,-0.8mm>*{}**@{-},
 <2.4mm,3mm>*{};<2.4mm,5.2mm>*{}**@{-},
     <0mm,-1.3mm>*{};<0mm,-5.3mm>*{^1}**@{},
     <2.5mm,2.3mm>*{};<5.1mm,-2.6mm>*{^2}**@{},
    <2.4mm,2.5mm>*{};<2.4mm,5.7mm>*{^{\bar{2}}}**@{},
    <-0.38mm,-0.2mm>*{};<-2.8mm,2.5mm>*{^{\bar{1}}}**@{},
    \end{xy}}\Ea
    =
    \Ba{c}\resizebox{6mm}{!}{
\mbox{$\xy
 (0.5,0.9)*{^{{^{\bar{1}}}}},
(0.5,5)*{^{{^{\bar{2}}}}},
 (0,-7)*+{_{_1}}*\frm{o}="C";
  (0,-14)*+{_{_2}}*\frm{o}="D";
(0,-2)*{\bu}="A";
(0,-2)*{\bu}="B";
"A"; "B" **\crv{(6,6) & (-6,6)};
 \ar @{-} "A";"C" <0pt>
 \ar @{-} "D";"C" <0pt>
\endxy$}}
\Ea
-
\Ba{c}\resizebox{9.8mm}{!}{
\mbox{$\xy
 (0.5,0.9)*{^{{^{\bar{1}}}}},
(0.5,5)*{^{{^{\bar{2}}}}},
 (-4,-7)*+{_{_1}}*\frm{o}="C";
  (4,-7)*+{_{_2}}*\frm{o}="D";
(0,-2)*{\bu}="A";
(0,-2)*{\bu}="B";
"A"; "B" **\crv{(6,6) & (-6,6)};
 \ar @{-} "A";"C" <0pt>
 \ar @{-} "D";"A" <0pt>
\endxy$}}
\Ea
+
\Ba{c}\resizebox{9.8mm}{!}{
\mbox{$\xy
 (0.5,0.9)*{^{{^{\bar{1}}}}},
(0.5,5)*{^{{^{\bar{2}}}}},
 (-4,-7)*+{_{_2}}*\frm{o}="C";
  (4,-7)*+{_{_1}}*\frm{o}="D";
(0,-2)*{\bu}="A";
(0,-2)*{\bu}="B";
"A"; "B" **\crv{(6,6) & (-6,6)};
 \ar @{-} "A";"C" <0pt>
 \ar @{-} "D";"A" <0pt>
\endxy$}}
\Ea
$$
The symmetrization over $(12)$ kills the last two terms in the last formula
so that the Drinfeld compatibility condition follows.

\mip

{\sc{Step 2}}: we check next the claim
$$
\oint_{123} j\left(\hspace{-2mm} \Ba{c}\resizebox{9.4mm}{!}{
\begin{xy}
 <0mm,0mm>*{\bu};<0mm,0mm>*{}**@{},
 <0mm,-0.49mm>*{};<0mm,-3.0mm>*{}**@{-},
 <0.49mm,0.49mm>*{};<1.9mm,1.9mm>*{}**@{-},
 <-0.5mm,0.5mm>*{};<-1.9mm,1.9mm>*{}**@{-},
 <-2.3mm,2.3mm>*{\bu};<-2.3mm,2.3mm>*{}**@{},
 <-1.8mm,2.8mm>*{};<0mm,4.9mm>*{}**@{-},
 <-2.8mm,2.9mm>*{};<-4.6mm,4.9mm>*{}**@{-},
   <0.49mm,0.49mm>*{};<2.7mm,2.3mm>*{^3}**@{},
   <-1.8mm,2.8mm>*{};<0.4mm,5.3mm>*{^2}**@{},
   <-2.8mm,2.9mm>*{};<-5.1mm,5.3mm>*{^1}**@{},
 \end{xy}}\Ea
+
\Ba{c}\resizebox{12.5mm}{!}{
\begin{xy}
(0,0)*{\bu};(-4,5)*{^{1}}**@{-},
(0,0)*{\bu};(0,5)*{^{2}}**@{-},
(0,0)*{\bu};(4,5)*{\bu}**@{-},
(4,5)*{\bu};(4,10)*{^{3}}**@{-},
(4,5)*{\bu};(8,0)*{_{\, 1}}**@{-},
\end{xy}
}
\Ea
\hspace{-2mm}
\right) \ \text{is}\ \sd_\centerdot\text{-exact}.
$$

The properadic composition in $\tw\cR\cG ra_0^{\geq 3}$ gives the following result
$$
\Ba{c}
\begin{xy}
 <0mm,0mm>*{\circ};<0mm,0mm>*{}**@{},
 <0mm,-0.49mm>*{};<0mm,-3.0mm>*{}**@{-},
 <-0.49mm,0.49mm>*{};<-1.9mm,1.9mm>*{}**@{-},
 <0.5mm,0.5mm>*{};<1.9mm,1.9mm>*{}**@{-},
 <2.3mm,2.3mm>*{\circ};<2.3mm,2.3mm>*{}**@{},
 <1.8mm,2.8mm>*{};<0mm,4.9mm>*{}**@{-},
 <2.8mm,2.9mm>*{};<4.6mm,4.9mm>*{}**@{-},
   <-0.49mm,0.49mm>*{};<-2.7mm,2.3mm>*{^{\bar{1}}}**@{},
   <1.8mm,2.8mm>*{};<-0.4mm,5.3mm>*{^{\bar{2}}}**@{},
   <2.8mm,2.9mm>*{};<5.1mm,5.3mm>*{^{\bar{3}}}**@{},
   <0mm,-4.3mm>*{_{_1}}**@{},
 \end{xy}
 \Ea
 =\frac{1}{4}\left(\Id - (\bar{2}\bar{3})\right)\left(
-
\Ba{c}{
\mbox{$\xy
(3.0,2.5)*{^{{^{\bar{1}}}}},
 (0.8,-1.1)*{^{{^{\bar{3}}}}},
(0.5,5)*{^{{^{\bar{2}}}}},
 (0,-8)*+{_{_1}}*\frm{o}="C";
(0,-3)*{\bu}="A1";
(0,-3)*{\bu}="A2";
(0,4)*{\bu}="B1";
(0,4)*{\bu}="B2";
"A1"; "A2" **\crv{(6,3) & (-6,3)};
"B1"; "B2" **\crv{(6,9) & (-6,9)};
 \ar @{-} "B1";"C" <0pt>
\endxy$}}
\Ea
+
\Ba{c}
\mbox{$\xy
(2.0,0.5)*{^{{^{\bar{3}}}}},
 (-3,-4.0)*{^{{^{\bar{1}}}}},
(0.5,5)*{^{{^{\bar{2}}}}},
 (0,-8)*+{_{_1}}*\frm{o}="C";
(0,-3)*{\bu}="A1";
(0,-3)*{\bu}="A2";
(0,4)*{\bu}="B1";
(0,4)*{\bu}="B2";
"A1"; "A2" **\crv{(-6,-9) & (-6,3)};
"B1"; "B2" **\crv{(6,9) & (-6,9)};
 \ar @{-} "B1";"C" <0pt>
\endxy$}
\Ea
+
\Ba{c}
\mbox{$\xy
(2.1,0.0)*{^{{^{\bar{3}}}}},
 (3,-4.0)*{^{{^{\bar{1}}}}},
(0.5,5)*{^{{^{\bar{2}}}}},
 (0,-8)*+{_{_1}}*\frm{o}="C";
(0,-3)*{\bu}="A1";
(0,-3)*{\bu}="A2";
(0,4)*{\bu}="B1";
(0,4)*{\bu}="B2";
"A1"; "A2" **\crv{(6,-9) & (6,3)};
"B1"; "B2" **\crv{(6,9) & (-6,9)};
 \ar @{-} "B1";"C" <0pt>
\endxy$}
\Ea
-
\Ba{c}
\mbox{$\xy
(2.0,2.5)*{^{{^{\bar{3}}}}},
 (0.0,5.8)*{^{{^{\bar{1}}}}},
(0.0,-8.6)*{^{{^{\bar{2}}}}},
 (0,0)*+{_{_1}}*\frm{o}="C";
(0,5)*{\bu}="A1";
(0,5)*{\bu}="A2";
(0,-5)*{\bu}="B1";
(0,-5)*{\bu}="B2";
"A1"; "A2" **\crv{(6,11) & (-6,11)};
"B1"; "B2" **\crv{(-6,-11) & (6,-11)};
 \ar @{-} "A1";"C" <0pt>
  \ar @{-} "B1";"C" <0pt>
\endxy$}
\Ea
\right)
$$
 One has
$$
\oint_{123} \left( -
\Ba{c}
\mbox{$\xy
(2.0,2.5)*{^{{^{\bar{3}}}}},
 (0.0,5.8)*{^{{^{\bar{1}}}}},
(0.0,-8.6)*{^{{^{\bar{2}}}}},
 (0,0)*+{_{_1}}*\frm{o}="C";
(0,5)*{\bu}="A1";
(0,5)*{\bu}="A2";
(0,-5)*{\bu}="B1";
(0,-5)*{\bu}="B2";
"A1"; "A2" **\crv{(6,11) & (-6,11)};
"B1"; "B2" **\crv{(-6,-11) & (6,-11)};
 \ar @{-} "A1";"C" <0pt>
  \ar @{-} "B1";"C" <0pt>
\endxy$}
\Ea
+
\Ba{c}
\mbox{$\xy
(2.0,2.5)*{^{{^{\bar{2}}}}},
 (0.0,5.8)*{^{{^{\bar{1}}}}},
(0.0,-8.6)*{^{{^{\bar{3}}}}},
 (0,0)*+{_{_1}}*\frm{o}="C";
(0,5)*{\bu}="A1";
(0,5)*{\bu}="A2";
(0,-5)*{\bu}="B1";
(0,-5)*{\bu}="B2";
"A1"; "A2" **\crv{(6,11) & (-6,11)};
"B1"; "B2" **\crv{(-6,-11) & (6,-11)};
 \ar @{-} "A1";"C" <0pt>
  \ar @{-} "B1";"C" <0pt>
\endxy$}
\Ea
    \right)=0
$$
because the linear combination of graphs on the l.h.s.\ admits an automorphism which
reverses  the orientation. Hence one can ignore the last summand in the formula for
$\Ba{c}
\begin{xy}
 <0mm,0mm>*{\circ};<0mm,0mm>*{}**@{},
 <0mm,-0.49mm>*{};<0mm,-3.0mm>*{}**@{-},
 <-0.49mm,0.49mm>*{};<-1.9mm,1.9mm>*{}**@{-},
 <0.5mm,0.5mm>*{};<1.9mm,1.9mm>*{}**@{-},
 <2.3mm,2.3mm>*{\circ};<2.3mm,2.3mm>*{}**@{},
 <1.8mm,2.8mm>*{};<0mm,4.9mm>*{}**@{-},
 <2.8mm,2.9mm>*{};<4.6mm,4.9mm>*{}**@{-},
   <-0.49mm,0.49mm>*{};<-2.7mm,2.3mm>*{^{\bar{1}}}**@{},
   <1.8mm,2.8mm>*{};<-0.4mm,5.3mm>*{^{\bar{2}}}**@{},
   <2.8mm,2.9mm>*{};<5.1mm,5.3mm>*{^{\bar{3}}}**@{},
 \end{xy}\Ea$ above.

 \sip

One has
$$
\sd_\centerdot
\Ba{c}
\mbox{$\xy
 (0.0,-8.8)*{^{{^{\bar{1}}}}},
(0.0,3.3)*{^{{^{\bar{2}}}}},
(2.0,-1.6)*{^{{^{\bar{3}}}}},
 (0,-4)*+{_{_1}}*\frm{o}="C1";
 (0,-4)*+{_{_1}}*\frm{o}="C2";
(0,2)*{\bu}="A1";
(0,2)*{\bu}="A2";
"A1"; "A2" **\crv{(6,8) & (-6,8)};
"C1"; "C2" **\crv{(-7,-12) & (7,-12)};
 \ar @{-} "A1";"C1" <0pt>
\endxy$}
\Ea
=
\Ba{c}{
\mbox{$\xy
(3.0,2.5)*{^{{^{\bar{1}}}}},
 (0.8,-1.1)*{^{{^{\bar{3}}}}},
(0.5,5)*{^{{^{\bar{2}}}}},
 (0,-8)*+{_{_1}}*\frm{o}="C";
(0,-3)*{\bu}="A1";
(0,-3)*{\bu}="A2";
(0,4)*{\bu}="B1";
(0,4)*{\bu}="B2";
"A1"; "A2" **\crv{(6,3) & (-6,3)};
"B1"; "B2" **\crv{(6,9) & (-6,9)};
 \ar @{-} "B1";"C" <0pt>
\endxy$}}
\Ea
+
\Ba{c}
\mbox{$\xy
(2.0,0.5)*{^{{^{\bar{3}}}}},
 (-3,-4.0)*{^{{^{\bar{1}}}}},
(0.5,5)*{^{{^{\bar{2}}}}},
 (0,-8)*+{_{_1}}*\frm{o}="C";
(0,-3)*{\bu}="A1";
(0,-3)*{\bu}="A2";
(0,4)*{\bu}="B1";
(0,4)*{\bu}="B2";
"A1"; "A2" **\crv{(-6,-9) & (-6,3)};
"B1"; "B2" **\crv{(6,9) & (-6,9)};
 \ar @{-} "B1";"C" <0pt>
\endxy$}
\Ea
+
\Ba{c}
\mbox{$\xy
(2.1,0.0)*{^{{^{\bar{3}}}}},
 (3,-4.0)*{^{{^{\bar{1}}}}},
(0.5,5)*{^{{^{\bar{2}}}}},
 (0,-8)*+{_{_1}}*\frm{o}="C";
(0,-3)*{\bu}="A1";
(0,-3)*{\bu}="A2";
(0,4)*{\bu}="B1";
(0,4)*{\bu}="B2";
"A1"; "A2" **\crv{(6,-9) & (6,3)};
"B1"; "B2" **\crv{(6,9) & (-6,9)};
 \ar @{-} "B1";"C" <0pt>
\endxy$}
\Ea
+
\Ba{c}
\mbox{$\xy
(2.0,2.5)*{^{{^{\bar{3}}}}},
 (0.0,5.8)*{^{{^{\bar{1}}}}},
(0.0,-8.6)*{^{{^{\bar{2}}}}},
 (0,0)*+{_{_1}}*\frm{o}="C";
(0,5)*{\bu}="A1";
(0,5)*{\bu}="A2";
(0,-5)*{\bu}="B1";
(0,-5)*{\bu}="B2";
"A1"; "A2" **\crv{(6,11) & (-6,11)};
"B1"; "B2" **\crv{(-6,-11) & (6,-11)};
 \ar @{-} "A1";"C" <0pt>
  \ar @{-} "B1";"C" <0pt>
\endxy$}
\Ea
$$
and
$$
\sd_\centerdot\left(
\Ba{c}
\mbox{$\xy
(0.0,0.5)*{^{{^{\bar{3}}}}},
 (-3,-2.0)*{^{{^{\bar{1}}}}},
(3,-2)*{^{{^{\bar{2}}}}},
 (0,-9)*+{_{_1}}*\frm{o}="C";
(0,-1)*{\bu}="A1";
(0,-1)*{\bu}="A2";
"A1"; "A2" **\crv{(-6,-7) & (-6,5)};
"A1"; "A2" **\crv{(6,-7) & (6,5)};
 \ar @{-} "A1";"C" <0pt>
\endxy$}
\Ea
-
\Ba{c}
\mbox{$\xy
(0.0,0.5)*{^{{^{\bar{3}}}}},
 (-3,-2.0)*{^{{^{\bar{2}}}}},
(3,-2)*{^{{^{\bar{1}}}}},
 (0,-9)*+{_{_1}}*\frm{o}="C";
(0,-1)*{\bu}="A1";
(0,-1)*{\bu}="A2";
"A1"; "A2" **\crv{(-6,-7) & (-6,5)};
"A1"; "A2" **\crv{(6,-7) & (6,5)};
 \ar @{-} "A1";"C" <0pt>
\endxy$}
\Ea
\right)
=\Ba{c}
\mbox{$\xy
(2.0,0.5)*{^{{^{\bar{3}}}}},
 (-3,-4.0)*{^{{^{\bar{1}}}}},
(0.5,5)*{^{{^{\bar{2}}}}},
 (0,-8)*+{_{_1}}*\frm{o}="C";
(0,-3)*{\bu}="A1";
(0,-3)*{\bu}="A2";
(0,4)*{\bu}="B1";
(0,4)*{\bu}="B2";
"A1"; "A2" **\crv{(-6,-9) & (-6,3)};
"B1"; "B2" **\crv{(6,9) & (-6,9)};
 \ar @{-} "B1";"C" <0pt>
\endxy$}
\Ea
+
\Ba{c}
\mbox{$\xy
(2.1,0.0)*{^{{^{\bar{3}}}}},
 (3,-4.0)*{^{{^{\bar{1}}}}},
(0.5,5)*{^{{^{\bar{2}}}}},
 (0,-8)*+{_{_1}}*\frm{o}="C";
(0,-3)*{\bu}="A1";
(0,-3)*{\bu}="A2";
(0,4)*{\bu}="B1";
(0,4)*{\bu}="B2";
"A1"; "A2" **\crv{(6,-9) & (6,3)};
"B1"; "B2" **\crv{(6,9) & (-6,9)};
 \ar @{-} "B1";"C" <0pt>
\endxy$}
\Ea
-
\Ba{c}
\mbox{$\xy
(2.0,0.5)*{^{{^{\bar{3}}}}},
 (-3,-4.0)*{^{{^{\bar{2}}}}},
(0.5,5)*{^{{^{\bar{1}}}}},
 (0,-8)*+{_{_1}}*\frm{o}="C";
(0,-3)*{\bu}="A1";
(0,-3)*{\bu}="A2";
(0,4)*{\bu}="B1";
(0,4)*{\bu}="B2";
"A1"; "A2" **\crv{(-6,-9) & (-6,3)};
"B1"; "B2" **\crv{(6,9) & (-6,9)};
 \ar @{-} "B1";"C" <0pt>
\endxy$}
\Ea
-
\Ba{c}
\mbox{$\xy
(2.1,0.0)*{^{{^{\bar{3}}}}},
 (3,-4.0)*{^{{^{\bar{2}}}}},
(0.5,5)*{^{{^{\bar{1}}}}},
 (0,-8)*+{_{_1}}*\frm{o}="C";
(0,-3)*{\bu}="A1";
(0,-3)*{\bu}="A2";
(0,4)*{\bu}="B1";
(0,4)*{\bu}="B2";
"A1"; "A2" **\crv{(6,-9) & (6,3)};
"B1"; "B2" **\crv{(6,9) & (-6,9)};
 \ar @{-} "B1";"C" <0pt>
\endxy$}
\Ea
+ 2\left(
\Ba{c}
\mbox{$\xy
(2.0,-1.5)*{^{{^{\bar{2}}}}},
 (-2,-1.5)*{^{{^{\bar{1}}}}},
(0.5,4)*{^{{^{\bar{3}}}}},
 (0,-4)*{\bu}="C";
  (0,-10)*+{_{_1}}*\frm{o}="D";
(0,3)*{\bu}="A1";
(0,3)*{\bu}="A2";
"C"; "A1" **\crv{(-6,-6) & (-6,5)};
"C"; "A2" **\crv{(6,-6) & (6,5)};
 \ar @{-} "A1";"C" <0pt>
  \ar @{-} "D";"C" <0pt>
\endxy$}
\Ea
-
\Ba{c}
\mbox{$\xy
(2.0,-1.5)*{^{{^{\bar{1}}}}},
 (-2,-1.5)*{^{{^{\bar{2}}}}},
(0.5,4)*{^{{^{\bar{3}}}}},
 (0,-4)*{\bu}="C";
  (0,-10)*+{_{_1}}*\frm{o}="D";
(0,3)*{\bu}="A1";
(0,3)*{\bu}="A2";
"C"; "A1" **\crv{(-6,-6) & (-6,5)};
"C"; "A2" **\crv{(6,-6) & (6,5)};
 \ar @{-} "A1";"C" <0pt>
  \ar @{-} "D";"C" <0pt>
\endxy$}
\Ea
\right)
$$
These results imply after a straightforward calculation that
$$
\oint_{123}\left(
\Ba{c}
\begin{xy}
 <0mm,0mm>*{\circ};<0mm,0mm>*{}**@{},
 <0mm,-0.49mm>*{};<0mm,-3.0mm>*{}**@{-},
 <-0.49mm,0.49mm>*{};<-1.9mm,1.9mm>*{}**@{-},
 <0.5mm,0.5mm>*{};<1.9mm,1.9mm>*{}**@{-},
 <2.3mm,2.3mm>*{\circ};<2.3mm,2.3mm>*{}**@{},
 <1.8mm,2.8mm>*{};<0mm,4.9mm>*{}**@{-},
 <2.8mm,2.9mm>*{};<4.6mm,4.9mm>*{}**@{-},
   <-0.49mm,0.49mm>*{};<-2.7mm,2.3mm>*{^{\bar{1}}}**@{},
   <1.8mm,2.8mm>*{};<-0.4mm,5.3mm>*{^{\bar{2}}}**@{},
   <2.8mm,2.9mm>*{};<5.1mm,5.3mm>*{^{\bar{3}}}**@{},
   <0mm,-4.3mm>*{_{_1}}**@{},
 \end{xy}
 \Ea
 +
 \frac{1}{2}
\Ba{c}
\mbox{$\xy
(2.0,-1.5)*{^{{^{\bar{2}}}}},
 (-2,-1.5)*{^{{^{\bar{1}}}}},
(0.5,4)*{^{{^{\bar{3}}}}},
 (0,-4)*{\bu}="C";
  (0,-10)*+{_{_1}}*\frm{o}="D";
(0,3)*{\bu}="A1";
(0,3)*{\bu}="A2";
"C"; "A1" **\crv{(-6,-6) & (-6,5)};
"C"; "A2" **\crv{(6,-6) & (6,5)};
 \ar @{-} "A1";"C" <0pt>
  \ar @{-} "D";"C" <0pt>
\endxy$}
\Ea
-\frac{1}{2}
\Ba{c}
\mbox{$\xy
(2.0,-1.5)*{^{{^{\bar{1}}}}},
 (-2,-1.5)*{^{{^{\bar{2}}}}},
(0.5,4)*{^{{^{\bar{3}}}}},
 (0,-4)*{\bu}="C";
  (0,-10)*+{_{_1}}*\frm{o}="D";
(0,3)*{\bu}="A1";
(0,3)*{\bu}="A2";
"C"; "A1" **\crv{(-6,-6) & (-6,5)};
"C"; "A2" **\crv{(6,-6) & (6,5)};
 \ar @{-} "A1";"C" <0pt>
  \ar @{-} "D";"C" <0pt>
\endxy$}
\Ea
\right)
=
\sd_\centerdot   \oint \left(\frac{1}{4}
\Ba{c}
\mbox{$\xy
 (0.0,-8.8)*{^{{^{\bar{1}}}}},
(0.0,3.3)*{^{{^{\bar{3}}}}},
(2.0,-1.6)*{^{{^{\bar{2}}}}},
 (0,-4)*+{_{_1}}*\frm{o}="C1";
 (0,-4)*+{_{_1}}*\frm{o}="C2";
(0,2)*{\bu}="A1";
(0,2)*{\bu}="A2";
"A1"; "A2" **\crv{(6,8) & (-6,8)};
"C1"; "C2" **\crv{(-7,-12) & (7,-12)};
 \ar @{-} "A1";"C1" <0pt>
\endxy$}
\Ea
- \frac{1}{4}
\Ba{c}
\mbox{$\xy
 (0.0,-8.8)*{^{{^{\bar{1}}}}},
(0.0,3.3)*{^{{^{\bar{2}}}}},
(2.0,-1.6)*{^{{^{\bar{3}}}}},
 (0,-4)*+{_{_1}}*\frm{o}="C1";
 (0,-4)*+{_{_1}}*\frm{o}="C2";
(0,2)*{\bu}="A1";
(0,2)*{\bu}="A2";
"A1"; "A2" **\crv{(6,8) & (-6,8)};
"C1"; "C2" **\crv{(-7,-12) & (7,-12)};
 \ar @{-} "A1";"C1" <0pt>
\endxy$}
\Ea
- \frac{1}{2}
\Ba{c}
\mbox{$\xy
(0.0,0.5)*{^{{^{\bar{3}}}}},
 (-3,-2.0)*{^{{^{\bar{1}}}}},
(3,-2)*{^{{^{\bar{2}}}}},
 (0,-9)*+{_{_1}}*\frm{o}="C";
(0,-1)*{\bu}="A1";
(0,-1)*{\bu}="A2";
"A1"; "A2" **\crv{(-6,-7) & (-6,5)};
"A1"; "A2" **\crv{(6,-7) & (6,5)};
 \ar @{-} "A1";"C" <0pt>
\endxy$}
\Ea
+
\frac{1}{2}
\Ba{c}
\mbox{$\xy
(0.0,0.5)*{^{{^{\bar{3}}}}},
 (-3,-2.0)*{^{{^{\bar{2}}}}},
(3,-2)*{^{{^{\bar{1}}}}},
 (0,-9)*+{_{_1}}*\frm{o}="C";
(0,-1)*{\bu}="A1";
(0,-1)*{\bu}="A2";
"A1"; "A2" **\crv{(-6,-7) & (-6,5)};
"A1"; "A2" **\crv{(6,-7) & (6,5)};
 \ar @{-} "A1";"C" <0pt>
\endxy$}
\Ea
\right)
$$

Finally one checks the equality
$$
\oint_{123} \frac{1}{2}\left(
\Ba{c}
\mbox{$\xy
(2.0,-1.5)*{^{{^{\bar{2}}}}},
 (-2,-1.5)*{^{{^{\bar{1}}}}},
(0.5,4)*{^{{^{\bar{3}}}}},
 (0,-4)*{\bu}="C";
  (0,-10)*+{_{_1}}*\frm{o}="D";
(0,3)*{\bu}="A1";
(0,3)*{\bu}="A2";
"C"; "A1" **\crv{(-6,-6) & (-6,5)};
"C"; "A2" **\crv{(6,-6) & (6,5)};
 \ar @{-} "A1";"C" <0pt>
  \ar @{-} "D";"C" <0pt>
\endxy$}
\Ea
-
\Ba{c}
\mbox{$\xy
(2.0,-1.5)*{^{{^{\bar{1}}}}},
 (-2,-1.5)*{^{{^{\bar{2}}}}},
(0.5,4)*{^{{^{\bar{3}}}}},
 (0,-4)*{\bu}="C";
  (0,-10)*+{_{_1}}*\frm{o}="D";
(0,3)*{\bu}="A1";
(0,3)*{\bu}="A2";
"C"; "A1" **\crv{(-6,-6) & (-6,5)};
"C"; "A2" **\crv{(6,-6) & (6,5)};
 \ar @{-} "A1";"C" <0pt>
  \ar @{-} "D";"C" <0pt>
\endxy$}
\Ea
\right)
=\oint_{123}
\Ba{c}\resizebox{12.5mm}{!}{
\begin{xy}
(0,0)*{\circ};(-4,5)*{^{1}}**@{-},
(0,0)*{\circ};(0,5)*{^{2}}**@{-},
(0,0)*{\circ};(4,5)*{\circ}**@{-},
(4,5)*{\circ};(4,10)*{^{3}}**@{-},
(4,5)*{\circ};(8,0)*{_{\, 1}}**@{-},
\end{xy}}
\Ea
$$
where on the r.h.s.\ we use the notation
$$
\Ba{c}\begin{xy}
 <0mm,-1mm>*{\circ};<-4mm,3mm>*{^{_1}}**@{-},
 <0mm,-1mm>*{\circ};<0mm,3mm>*{^{_2}}**@{-},
 <0mm,-1mm>*{\circ};<4mm,3mm>*{^{_3}}**@{-},
 \end{xy}\Ea
 :=
 \frac{1}{2}\left(-
 \Ba{c}\resizebox{11mm}{!}{
\xy
(2.0,-3.5)*{^{{^{\bar{2}}}}},
 (-2,-3.5)*{^{{^{\bar{1}}}}},
(0.5,4)*{^{{^{\bar{3}}}}},
 (0,-8)*{\bu}="C";
(0,3)*{\bu}="A1";
(0,3)*{\bu}="A2";
"C"; "A1" **\crv{(-5,-9) & (-5,4)};
"C"; "A2" **\crv{(5,-9) & (5,4)};
 \ar @{-} "A1";"C" <0pt>
\endxy}
\Ea
+
\Ba{c}\resizebox{11mm}{!}{
\xy
(2.0,-3.5)*{^{{^{\bar{1}}}}},
 (-2,-3.5)*{^{{^{\bar{2}}}}},
(0.5,4)*{^{{^{\bar{3}}}}},
 (0,-8)*{\bu}="C";
(0,3)*{\bu}="A1";
(0,3)*{\bu}="A2";
"C"; "A1" **\crv{(-5,-9) & (-5,4)};
"C"; "A2" **\crv{(5,-9) & (5,4)};
 \ar @{-} "A1";"C" <0pt>
\endxy}
\Ea\right)
$$
and the properadic composition in  $\tw\RGra_0^{\geq 3}$.   Hence the claim in {\sc Step 2} follows.

\mip

{\sc Step 3}: we check the last non-obvious claim that the twisted ribbon graph
$$
\oint_{123}j\left(\hspace{-1.5mm}
\Ba{c}\resizebox{11.4mm}{!}{
\begin{xy}
(0,0)*{\bu};(-4,5)*{^{1}}**@{-},
(0,0)*{\bu};(0,5)*{^{2}}**@{-},
(0,0)*{\bu};(4,5)*{\bu}**@{-},
(4,5)*{\bu};(1.5,10)*{^{3}}**@{-},
(4,5)*{\bu};(6.5,10)*{^{4}}**@{-},
\end{xy}}
\Ea
-
\Ba{c}\resizebox{11.4mm}{!}{
\begin{xy}
(0,0)*{\bu};(-4,5)*{^{4}}**@{-},
(0,0)*{\bu};(0,5)*{^{1}}**@{-},
(0,0)*{\bu};(4,5)*{\bu}**@{-},
(4,5)*{\bu};(1.5,10)*{^{2}}**@{-},
(4,5)*{\bu};(6.5,10)*{^{3}}**@{-},
\end{xy}}
\Ea\hspace{-1mm}\right)\equiv\left(\Id - \zeta +\zeta^2 - \zeta^3 - (23) - (24)\right)
 j
\left(\Ba{c}\resizebox{11.4mm}{!}{
\begin{xy}
(0,0)*{\bu};(-4,5)*{^{1}}**@{-},
(0,0)*{\bu};(0,5)*{^{2}}**@{-},
(0,0)*{\bu};(4,5)*{\bu}**@{-},
(4,5)*{\bu};(1.5,10)*{^{3}}**@{-},
(4,5)*{\bu};(6.5,10)*{^{4}}**@{-},
\end{xy}}
\Ea\right) \ \  \ \text{is $\sd_\centerdot$-exact},
$$
where $\zeta:=(1234)\in \bS_4$.
Rules for the properadic composition in $\tw \RGra_0^{\geq 3}$ give
$$
\Ba{c}\resizebox{11.4mm}{!}{
\begin{xy}
(0,0)*{\circ};(-4,5)*{^{1}}**@{-},
(0,0)*{\circ};(0,5)*{^{2}}**@{-},
(0,0)*{\circ};(4,5)*{\circ}**@{-},
(4,5)*{\circ};(1.5,10)*{^{3}}**@{-},
(4,5)*{\circ};(6.5,10)*{^{4}}**@{-},
\end{xy}}
\Ea
=\frac{1}{2}
\left(-\Ba{c}
\mbox{$\xy
(2.0,-1.1)*{^{{^{\bar{2}}}}},
 (-2,-1.1)*{^{{^{\bar{1}}}}},
(2.5,5)*{^{{^{\bar{3}}}}},
(0,10)*{^{{^{\bar{4}}}}},
 (0,4)*{\bu}="C";
  (0,7)*{\bu}="D";
(0,-3)*{\bu}="A1";
(0,-3)*{\bu}="A2";
"C"; "A1" **\crv{(-6,6) & (-6,-5)};
"C"; "A2" **\crv{(6,6) & (6,-5)};
"D"; "D" **\crv{(-6,15) & (6,15)};
 \ar @{-} "A1";"C" <0pt>
  \ar @{-} "D";"C" <0pt>
\endxy$}
\Ea
+
\Ba{c}
\mbox{$\xy
(2.0,-1.1)*{^{{^{\bar{1}}}}},
 (-2,-1.1)*{^{{^{\bar{2}}}}},
(2.5,5)*{^{{^{\bar{3}}}}},
(0,10)*{^{{^{\bar{4}}}}},
 (0,4)*{\bu}="C";
  (0,7)*{\bu}="D";
(0,-3)*{\bu}="A1";
(0,-3)*{\bu}="A2";
"C"; "A1" **\crv{(-6,6) & (-6,-5)};
"C"; "A2" **\crv{(6,6) & (6,-5)};
"D"; "D" **\crv{(-6,15) & (6,15)};
 \ar @{-} "A1";"C" <0pt>
  \ar @{-} "D";"C" <0pt>
\endxy$}
\Ea
+
\Ba{c}
\mbox{$\xy
(2.0,-1.1)*{^{{^{\bar{2}}}}},
 (-2,-1.1)*{^{{^{\bar{1}}}}},
(2.5,5)*{^{{^{\bar{4}}}}},
(0,10)*{^{{^{\bar{3}}}}},
 (0,4)*{\bu}="C";
  (0,7)*{\bu}="D";
(0,-3)*{\bu}="A1";
(0,-3)*{\bu}="A2";
"C"; "A1" **\crv{(-6,6) & (-6,-5)};
"C"; "A2" **\crv{(6,6) & (6,-5)};
"D"; "D" **\crv{(-6,15) & (6,15)};
 \ar @{-} "A1";"C" <0pt>
  \ar @{-} "D";"C" <0pt>
\endxy$}
\Ea
-
\Ba{c}
\mbox{$\xy
(2.0,-1.1)*{^{{^{\bar{2}}}}},
 (-2,-1.1)*{^{{^{\bar{1}}}}},
(2.5,5)*{^{{^{\bar{4}}}}},
(0,10)*{^{{^{\bar{3}}}}},
 (0,4)*{\bu}="C";
  (0,7)*{\bu}="D";
(0,-3)*{\bu}="A1";
(0,-3)*{\bu}="A2";
"C"; "A1" **\crv{(-6,6) & (-6,-5)};
"C"; "A2" **\crv{(6,6) & (6,-5)};
"D"; "D" **\crv{(-6,15) & (6,15)};
 \ar @{-} "A1";"C" <0pt>
  \ar @{-} "D";"C" <0pt>
\endxy$}
\Ea
\right)
$$

A straightforward but tedious calculation shows that the above claim is true indeed with the $\sd_\centerdot$-exact term given by
$$
\oint_{123}j\left(\hspace{-1.5mm}
\Ba{c}\resizebox{11.4mm}{!}{
\begin{xy}
(0,0)*{\bu};(-4,5)*{^{1}}**@{-},
(0,0)*{\bu};(0,5)*{^{2}}**@{-},
(0,0)*{\bu};(4,5)*{\bu}**@{-},
(4,5)*{\bu};(1.5,10)*{^{3}}**@{-},
(4,5)*{\bu};(6.5,10)*{^{4}}**@{-},
\end{xy}}
\Ea
-
\Ba{c}\resizebox{11.4mm}{!}{
\begin{xy}
(0,0)*{\bu};(-4,5)*{^{4}}**@{-},
(0,0)*{\bu};(0,5)*{^{1}}**@{-},
(0,0)*{\bu};(4,5)*{\bu}**@{-},
(4,5)*{\bu};(1.5,10)*{^{2}}**@{-},
(4,5)*{\bu};(6.5,10)*{^{3}}**@{-},
\end{xy}}
\Ea\hspace{-1mm}\right)=\hspace{110mm}
$$
$$
\hspace{23mm}
=
\frac{1}{2}\left(\Id - \zeta +\zeta^2 - \zeta^3 - (23) - (24)\right)\sd_\centerdot
\left(
-
\Ba{c}
\mbox{$\xy
(2.0,4.0)*{^{{^{\bar{3}}}}},
 (-3,-0.2)*{^{{^{\bar{1}}}}},
(3,-0.2)*{^{{^{\bar{2}}}}},
(0,10)*{^{{^{\bar{4}}}}},
 (0,7)*{\bu}="C";
(0,1)*{\bu}="A1";
(0,1)*{\bu}="A2";
"A1"; "A2" **\crv{(-6,7) & (-6,-5)};
"A1"; "A2" **\crv{(6,7) & (6,-5)};
"C"; "C" **\crv{(-6,15) & (6,15)};
 \ar @{-} "A1";"C" <0pt>
\endxy$}
\Ea
+
\Ba{c}
\mbox{$\xy
(2.0,4.0)*{^{{^{\bar{3}}}}},
 (-3,-0.2)*{^{{^{\bar{2}}}}},
(3,-0.2)*{^{{^{\bar{1}}}}},
(0,10)*{^{{^{\bar{4}}}}},
 (0,7)*{\bu}="C";
(0,1)*{\bu}="A1";
(0,1)*{\bu}="A2";
"A1"; "A2" **\crv{(-6,7) & (-6,-5)};
"A1"; "A2" **\crv{(6,7) & (6,-5)};
"C"; "C" **\crv{(-6,15) & (6,15)};
 \ar @{-} "A1";"C" <0pt>
\endxy$}
\Ea
+
\Ba{c}
\mbox{$\xy
(2.0,4.0)*{^{{^{\bar{4}}}}},
 (-3,-0.2)*{^{{^{\bar{1}}}}},
(3,-0.2)*{^{{^{\bar{2}}}}},
(0,10)*{^{{^{\bar{3}}}}},
 (0,7)*{\bu}="C";
(0,1)*{\bu}="A1";
(0,1)*{\bu}="A2";
"A1"; "A2" **\crv{(-6,7) & (-6,-5)};
"A1"; "A2" **\crv{(6,7) & (6,-5)};
"C"; "C" **\crv{(-6,15) & (6,15)};
 \ar @{-} "A1";"C" <0pt>
\endxy$}
\Ea
-
\Ba{c}
\mbox{$\xy
(2.0,4.0)*{^{{^{\bar{4}}}}},
 (-3,-0.2)*{^{{^{\bar{1}}}}},
(3,-0.2)*{^{{^{\bar{2}}}}},
(0,10)*{^{{^{\bar{3}}}}},
 (0,7)*{\bu}="C";
(0,1)*{\bu}="A1";
(0,1)*{\bu}="A2";
"A1"; "A2" **\crv{(-6,7) & (-6,-5)};
"A1"; "A2" **\crv{(6,7) & (6,-5)};
"C"; "C" **\crv{(-6,15) & (6,15)};
 \ar @{-} "A1";"C" <0pt>
\endxy$}
\Ea
\right)
$$
The proof is completed.
\end{proof}

\subsection{K.\ Costello's moduli spaces as a dg properad}\label{3: subsec on Costello moduli spaces}
For integers $g,n,r \geq 0$ and $m > 1$, let  ${\cN}_{g,m,r,n}$ be the moduli space of connected stable Riemann
surfaces $\Sigma$ of genus $g$ with $m$ {\em marked}\
boundary components, $n$ {\em marked}\, points
in the interior $\Sigma\setminus \p\Sigma$ and with $r$ {\em marked}\, points on the boundary $\p\Sigma$. It is assumed that ${\cN}_{g,m,r,n}=\emptyset$ if $g=0$, $m=1$ and $r+2n<3$ or $g=0$, $m=2$, $r=n=0$.

\sip

 Kevin Costello introduces in \cite{Co1,Co2} a certain partial compactification $\overline{\cN}_{g,m,r,n}$ of ${\cN}_{g,m,r,n}$ by allowing stable Riemann surfaces $\Sigma$ with nodes {\it on the boundary}; the marked points are required to be distinct from the
nodes and each other, and the stability condition requires that the normalization --- cutting each node into two separate ``half-nodes" ---  of each singular surface $\Sigma$ gives a disconnected surface whose every connected component is stable in the above sense.  The moduli space $\overline{\cN}_{g,m,r,n}$ is a  $(6g-6 + 3m + r+ 2n)$-dimensional orbifold with corners so that the inclusion ${\cN}_{g,m,r,n}\hook \overline{\cN}_{g,m,r,n}$ is a homotopy equivalence. The boundary $\p\overline{\cN}_{g,m,r,n}$ consists of singular Riemann surfaces, the ones with nodes. It is proven in \cite{Co2} that the inclusion $\p\overline{\cN}_{g,m,r,n}\hook \overline{\cN}_{g,m,r,n}$ is also a homotopy equivalence.

\sip

Let  ${\caD}_{g,m,r,n}\subset \overline{\cN}_{g,m,r,n}$ be the locus where all irreducible components of surfaces are discs (with {\em at most one}\, internal marked point). A remarkable theorem by  K.\ Costello \cite{Co2} says that this inclusion is a homotopy equivalence so that
\Beq\label{3: H(D)=H(N)}
H_\bu({\caD}_{g,m,r,n})=H_\bu({\cN}_{g,m,r,n})
\Eeq
 Moreover, the compact orbi-space ${\caD}_{g,m,r,n}$   of disks admits a decomposition  into orbi-cells which can be parameterized by ribbon graphs.
We are interested in this paper in the case $r=0$ only; the general case $r\geq 1$ employs twisted ribbon graphs {\it with hairs}\, which we do not need in this paper.

\sip

The space ${\caD}_{g,m,0,n}$ is decomposed into orbi-cells via an equivalence relation,
$\Sigma \sim \Sigma'$,  if there exists a homeomorphism $\Sigma\rar \Sigma'$ which preserves the marked points and the orientation.
Each orbi-cell $O$ is uniquely determined by \cite{Co2}
\Bi
\item[-] the set
$V_{\bu}$ of its irreducible components (disks) with no marked points inside,
\item[-] the set  $V_{\circ}$ of its irreducible components (disks) with precisely one marked point inside; this set comes equipped with a fixed isomorphism $V_\circ\rar [n]$
    which associates to an  irreducible component the marking of its unique internal point;
\item[-] the set of ``half-nodes" $H(v)$ on each irreducible component $v\in V_{\circ}\sqcup V_{\bu}$; the stability condition requires $\# H(v)\geq 3$ for any $v\in V_{\bu}$ and
    $\# H(v)\geq 1$ for any $v\in V_{\circ}$; each such ``half-node" is  a point on the boundary ($\simeq S^1$) of the irreducible component $v$; hence each set $H(v)$ comes equipped with a fixed cyclic ordering;
    \item[-] the fixed point free involution $\tau$ on the disjoint union
$$
H:=\coprod_{v\in V_{\circ}\sqcup V_{\bu}} H(v)
$$
 controlling which ``half-nodes"  in the boundaries of disks are glued into genuine nodes; the set of orbits of this involution, i.e.\ the set of nodes, is denoted by $E$; for each $e\in E$ the associated 2-element orbit of $\tau$ is denoted by $H(e)$;

\item[-] the orientation of cell $O$ can be encoded, according to \cite{Co2}, into a total ordering of the set
$$
V_\bu \sqcup H
$$
defined  up to an even permutation.
\Ei

As the involution $\tau$ decomposes the set $H$ into the disjoint union of its orbits,
$$
H=\coprod_{e\in E}H(e).
$$
an orientation of the cell $O$ can be equivalently encoded into the choice
of one of the two possible unit vectors of the 1-dimensional Euclidean space (cf.\ \S {\ref{3: subsec on orientations for different d}})
$$
\det(V_\bu(\Ga)) \bigotimes_{e\in E(\Ga)}\det(H(e))
$$

\sip

Thus each orbi-cell $O\subset {\caD}_{g,m,0,n}$ gives us a uniquely defined oriented ribbon graph $\Ga$ such that
\Bi
\item[(i)] its set of vertices comes decomposed into the disjoint union of black and white vertices,
    $$
    V_{\circ}(\Ga)\sqcup V_{\bu}(\Ga)
    $$
and a marking isomorphism $V_{\circ}(\Ga)\rar [n]$ is fixed; each vertex $v$ has $H(v)$ half-edges attached with $\# H(v)\geq 3$ for $v\in V_\bu(\Ga)$ and $\#H(v)\geq 1$ for $v\in V_\circ(\Ga)$;
\item[(ii)] it has $m$ labelled boundaries;
\item[(iii)] its set of edges $E(\Ga)$ is in one-to-one correspondence with the set of nodes $E$ of surfaces in the cell $O$, or, equivalently, with the set of orbits of the involution $\tau$ on $H$;
\item[(iv)] it is assigned the cohomological degree
$$
|\Ga|_{geom}:=-\dim O= 3\#V_\bu(\Ga) + \# V_\circ(\Ga) - 2 \# E(\Ga)
$$
\item[(v)] it is equipped with an orientation in exactly the same sense as graphs from $\tw\RGra_d(m,n)$ for $d$ odd (see \S {\ref{3: subsec on twRgra_d}});
\item[(vi)] the equality
$$
 g= 1+\frac{1}{2}\left(\# E(\Ga) - \# V_\bu(\Ga)- n-m\right)
$$
holds true.
\Ei

\sip

The correspondence between the set of cells $O$ of $\caD_{g,m,0,n}$  and the set  oriented ribbon graphs $\Ga$ satisfying the conditions (i)-(iv) above is one-to-one. Given such a graph $\Ga$,
the associated cell $O_\Ga$ is given by
$$
O_\Ga=\left(\prod_{v\in V_\bu(\Ga)} \frac{\{H(v)\hook S^1\}}{SL(2,\R)} \times \prod_{v\in V_\circ(\Ga)} \frac{\{H(v)\hook S^1\}}{S^1}   \right)_{Aut\Ga}
$$
where $Aut\Ga$ is the automorphism groups of $\Ga$. This formula explains the above formula for the geometric degree
$|\Ga|=-\dim O_\Ga=3\#V_\bu(\Ga)  - 2 \# E(\Ga) +n$. On the other hand, a ribbon graph in  $\tw\RGra_3$ with exactly the same characteristics would have the cohomological degree
$3\#V_\bu(\Ga)  - 2 \# E(\Ga)$. Therefore we observe a one-to-one correspondence
between the set of generators (=cells) of the K.\ Costello's  cell complex, $\cC hains(\caD_{g,m,0,n})$, and the set of ribbon graph generators of the complex $\tw\RGra_3(m,n;g)[-n]$ introduced in \S {\ref{3: subsec on twRgra_d}} implying the canonical isomorphism,
$$
\cC hains(\caD_{g,m,0,n})= \tw\RGra_3^{\geq 3}(m,n;g)[-n]
$$
of graded vector spaces. It is immediate to see that the differentials of both sides are also identical to each other so that, using (\ref{3: degree shift iso for twGra_d}),  we get a canonical isomorphism of complexes (in fact, of dg $\bS_m^{op}\times \bS_n$-bimodules),
$$
\cC hains(\caD_{g,m,0,n})= \tw\RGra_3^{\geq 3}(m,n;g)[-n]=\tw\RGra_0^{\geq 3}(m,n;g)[6g-6+3m+2n]
$$
or, equivalently,
$$
\tw\RGra_0^{\geq 3}(m,n;g)= \cC hains(\caD_{g,m,0,n})[6-6g-3m-2n]
$$
Thus we proved the following

\subsubsection{\bf Proposition-definition} {\it The $\bS$-bimodule
$$
\left\{\cC hains(\caD_{g,m,0,n})[6-6g-3m-2n]\right\}_{m\geq 1, n\geq 0, 2g+m+n \geq 3}
$$
has a canonical structure of dg properad with properadic compositions given by the explicit combinatorial rules for $\tw\RGra_0^{\geq 3}$.}

\mip

K.\ Costello's isomorphism (\ref{3: H(D)=H(N)}) and the Poincare duality imply
$$
H_\bu(\caD_{g,m,0,n})=H_{\bu}(\cN_{g,m,0,n})=H_c^{\bu+ 6g-6 +3m+2n}(\cN_{g,m,0,n})
$$
On the other hand, the uniformization theorem gives us an equality of compactly supported cohomology groups
$$
H_c^{\bu}(\cN_{g,m,0,n})=H_c^{\bu}(\cM_{g,m+n}\times \R_{>0}^m)=H_c^{\bu -m}(\cM_{g,m+n})
$$
Combining all the above isomorphisms we finally obtain
$$
H^\bu(\tw\RGra_0^{\geq 3}(m,n;g))=H_\bu(\caD_{g,m,0,n})[6-6g-3m-2n]=H_c^{\bu}(\cN_{g,m,0,n})=H_c^{\bu -m}(\cM_{g,m+n})
$$
which proves finally Theorem {\ref{3: Main theorem on H(twRgra)}} and hence Theorem A in the Introduction.

\subsection{Examples of non-trivial properadic compositions of cohomology classes}\label{3: subsec on non-triv of prop compositions}
Let us call the graph with one vertex $v$ and $m$ (resp. $n$) outgoing (resp.\ incoming) edges
$$
v=\Ba{c}\resizebox{18mm}{!}
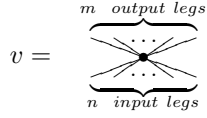
{
\begin{xy}
 <0mm,0mm>*{\bullet};
 <-0.5mm,0.2mm>*{};<-8mm,3mm>*{}**@{-},
 <-0.4mm,0.3mm>*{};<-4.5mm,3mm>*{}**@{-},
 <0mm,0mm>*{};<0mm,2.6mm>*{\ldots}**@{},
 <0.4mm,0.3mm>*{};<4.5mm,3mm>*{}**@{-},
 <0.5mm,0.2mm>*{};<8mm,3mm>*{}**@{-},
<-0.4mm,-0.2mm>*{};<-8mm,-3mm>*{}**@{-},
 <-0.5mm,-0.3mm>*{};<-4.5mm,-3mm>*{}**@{-},
 <0mm,0mm>*{};<0mm,-2.6mm>*{\ldots}**@{},
 <0.5mm,-0.3mm>*{};<4.5mm,-3mm>*{}**@{-},
 <0.4mm,-0.2mm>*{};<8mm,-3mm>*{}**@{-};
<0mm,5mm>*{\overbrace{\ \ \ \ \ \ \ \ \ \ \ \ \ \  }};
<0mm,-5mm>*{\underbrace{\ \ \ \ \ \ \ \ \ \ \ \ \ \ }};
<0mm,7mm>*{^{m\ \ output\ legs}};
<0mm,-7mm>*{_{n\ \ input\ legs}};
 \end{xy}}\Ea
$$
the {\em  $(m,n)$-corolla}; we denote $|v|_{in}:=n$ and $|v|_{out}=m$. A generic {\em directed graph $\Ga$ with legs}\, is built from such corollas by gluing some (or all) out-legs of one  corolla to the some (or all) in-legs of another corolla and creating thereby {\it internal edges}\, of $\Ga$; those in-legs and out-legs of corollas which remain free become  {\em in-}\, and respectively {\em out-legs of the directed graph $\Ga$}. The corollas used in the construction of $\Ga$ are called {\em vertices}\ of $\Ga$ and their set is denoted by $V(\Ga)$; the set of internal edges of $\Ga$ is denoted by $E(\Ga)$, while the set of in-legs (resp., out-legs) of $\Ga$ is denoted by $L_{in}(\Ga)$ (resp., $L_{out}(\Ga)$. We always assume that in- and out legs are distinguished, e.g.\ some isomorphisms
 $$
L_{in}(\Ga)\lon [\# L_{in}(\Ga)]
\ \ , \ \ L_{out}(\Ga)\lon [\# L_{out}(\Ga)]
$$
are fixed. The number
$$
g(\Ga):= \# E(\Ga) - \# V(\Ga) +1
$$
is called the {\em loop number}\, of $\Ga$.
 For example,
 $$
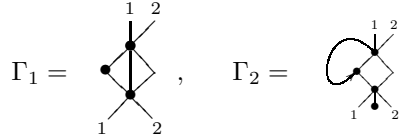
\Ga_1= \Ba{c}\resizebox{10mm}{!}{
\begin{xy}
 <0mm,0mm>*{\bullet};
<0.39mm,0.39mm>*{};<3.4mm,3.4mm>*{}**@{-},
<0mm,0.39mm>*{};<0mm,3.4mm>*{}**@{-},
 <0.39mm,-0.39mm>*{};<3.4mm,-3.7mm>*{}**@{-},
 <-0.35mm,-0.35mm>*{};<-2.9mm,-2.9mm>*{}**@{-},
 <-3.4mm,-3.4mm>*{\bullet};
 <-3.4mm,-3.4mm>*{};<0mm,-6.8mm>*{}**@{-},
 <3.4mm,-3.7mm>*{};<0mm,-6.8mm>*{}**@{-},
  <0mm,-6.8mm>*{\bullet};
  <0mm,-6.8mm>*{};<0mm,0mm>*{\bullet}**@{-},
  <0mm,-6.8mm>*{};<3mm,-10mm>*{}**@{-},
<0mm,-6.8mm>*{};<-3mm,-10mm>*{}**@{-},
  <0.4mm,-8.5mm>*{};<3.8mm,-12mm>*{^2}**@{},
 <0.4mm,-6.5mm>*{};<-3.8mm,-12mm>*{^1}**@{},
 <0mm,0mm>*{};<0mm,4.5mm>*{^1}**@{},
 <0mm,0mm>*{};<3.6mm,4.5mm>*{^2}**@{},
 \end{xy}}
 \Ea
 ,\ \ \ \ \
 \Ga_2= \Ba{c}\resizebox{12mm}{!}{
\begin{xy}
 <0mm,0mm>*{\bullet};
<0.39mm,0.39mm>*{};<3.4mm,3.4mm>*{}**@{-},
<0mm,0.39mm>*{};<0mm,3.4mm>*{}**@{-},
 <0.39mm,-0.39mm>*{};<3.4mm,-3.7mm>*{}**@{-},
 <-0.35mm,-0.35mm>*{};<-2.9mm,-2.9mm>*{}**@{-},
 <-3.4mm,-3.4mm>*{\bullet};
 <-3.4mm,-3.4mm>*{};<0mm,-6.8mm>*{}**@{-},
 <3.4mm,-3.7mm>*{};<0mm,-6.8mm>*{}**@{-},
  <0mm,-6.8mm>*{\bullet};
  <0mm,-6.8mm>*{};<0mm,-10mm>*{\bullet}**@{-},
  <0mm,-6.8mm>*{};<3mm,-10mm>*{}**@{-},
<0mm,-6.8mm>*{};<-3mm,-10mm>*{}**@{-},
  <0.4mm,-8.5mm>*{};<3.8mm,-12mm>*{^2}**@{},
 <0.4mm,-6.5mm>*{};<-3.8mm,-12mm>*{^1}**@{},
 <0mm,0mm>*{};<0mm,4.5mm>*{^1}**@{},
 <0mm,0mm>*{};<3.6mm,4.5mm>*{^2}**@{},
(-0.39,0.39)*{}
   \ar@{->}@(ul,dl) (-3.6,-3.6)*{}
 \end{xy}}
\Ea
$$
are directed graphs with loop number $2$, two in-legs, two out-legs and with $\#V(\Ga_1)=3$ and $\#V(\Ga_2)=4$; every leg and edge above is assumed to be directed from bottom to the top unless the direction shown explicitly. A directed graph is called {\em oriented}\, if it has no {\em closed}\, paths of directed edges. The first graph $\Ga_1$ is oriented, while the second one $\Ga_2$ is not. A directed graph is called {\em connected}\, if every two vertices are can be connected by a (not necessarily directed) path of edges.

\sip

A properad structure in $\Grav$ means that for any oriented graph $\Ga$ with $m$ out-legs, $n$ in-legs
and with loop number $g$ there is an associated linear map
$$
\Phi_\Ga: \bigotimes_{v\in V(\Ga)} H_c^{k_v}\left(\cM_{g_v; |v|_{out} + |v|_{in}\times \R_+^{|v|_{out}}}\right) \lon H^{\sum_v k_v}\left(\cM_{g+\sum_v g_v; m+n}\times \R_+^m \right)
$$
 satisfying suitable compatibility conditions. The automorphism group $Aut(\Ga)$ acts naturally on the (unordered) tensor product in the l.h.s., and the map $\Phi_\Ga$ must be invariant under that action.
  A class from $H_c^{k_v}(\cM_{g_v; |v|_{out} + |v|_{in}\times \R_+^{|v|_{out}}})$ can be understood in this context as a decoration of the vertex $v\in V(\Ga)$.

\sip

Let us test non-triviality of gravity properadic compositions for oriented graphs $\Ga$ with $g(\Ga)>0$ in the simplest possible case when all vertices are decorated by genus zero cohomology classes $g_v=0$ with $|v|_{out} + |v|_{in}=3$; all such oriented graphs are trivalent and can be identified with elements of $q\LB_{-1,0}$. As a first example consider a loop number 1 graph with 4 out-legs and no in-legs,
$$
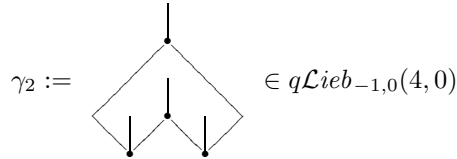
\ga_{2}:= \Ba{c}\resizebox{20mm}{!}{\xy
 (-8,-8)*{\bullet}="d1",
(8,-8)*{\bu}="d2",
(-16,0)*{}="l",
(16,0)*{}="r",
(0,16)*{\bu}="top",
(0,0)*{\bu}="u1",
 (0,24)*{}="1",
(0,8)*{}="2",
(-8,0)*{}="3",
(8,0)*{}="4",
\ar @{-} "d1";"l" <0pt>
\ar @{-} "d1";"u1" <0pt>
\ar @{-} "d2";"u1" <0pt>
\ar @{-} "d2";"r" <0pt>
\ar @{-} "l";"top" <0pt>
\ar @{-} "r";"top" <0pt>
\ar @{-} "top";"1" <0pt>
\ar @{-} "u1";"2" <0pt>
\ar @{-} "d1";"3" <0pt>
\ar @{-} "d2";"4" <0pt>
\endxy}\Ea \in q\LB_{-1,0}(4,0)
$$
The l.h.s\ of the associated properadic composition
$$
\Phi_{\ga_2}: H^3_c(\cM_{0,3+0}\times \R_+^3)\ot H^1_c(\cM_{0,1+2}\times \R_+)\ot H^3_c(\cM_{0,3+0}\times \R_+^3)\ot H^1_c(\cM_{0,1+2}\times \R_+) \lon  H^8_c(\cM_{1,4+0}\times \R_+^4)
$$
is a one-dimensional vector space $\K[-8]$. The image of a generator of $\K[-8]$ under $\Phi_{\ga_2}$ is given by composing images under $j$ of the second and third corollas in (\ref{3: map j from qLB to H(TwGra)}) in accordance with the combinatorial composition rules in $\tw\RGra_0$; in this case the answer is straightforward to compute --- it is given by the  ribbon graph
$$
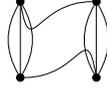
\Ba{c}\resizebox{16mm}{!}{
 \begin{tikzpicture}[baseline=-.65ex]
  \node[int] (v0) at (-45:1) {};
  \node[int] (v1) at (45:1) {};
\node[int] (v2) at (135:1) {};
\node[int] (v3) at (225:1) {};
\draw (v0) edge (v1) edge[bend left] (v1) edge[bend right] (v1);
\draw (v2) edge (v3) edge[bend left] (v3) edge[bend right] (v3);
\draw (v1) to[ out=180, in=-80, looseness=2] (v2);
\draw (v0) to[ out=100, in=0, looseness=2] (v3);
 \end{tikzpicture}
 }\Ea
$$
which represents (in Penner's ribbon graph complex) a non-zero cohomology class in $H^4(\cM_{1,4})$
(cf.\ \S 7.4 in \cite{MW}). In a similar way one shows that images of the following infinite sequence of elements in $q\LB_{-1,0}$,
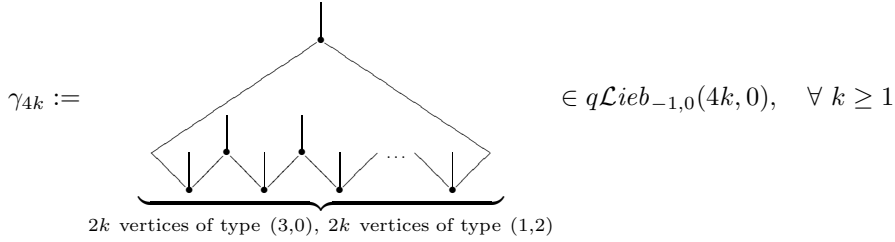
$$
\ga_{4k}:= \underbrace{\Ba{c}\resizebox{45mm}{!}{\xy
 (-8,-8)*{\bullet}="d1",
(8,-8)*{\bu}="d2",
(24,-8)*{\bu}="d3",
(-16,0)*{}="l",
(20,24)*{\bu}="top",
(0,0)*{\bu}="u1",
(16,0)*{\bu}="u2",
 (20,32)*{}="1",
(0,8)*{}="2",
(-8,0)*{}="3",
(8,0)*{}="4",
(16,8)*{}="5",
(24,0)*{}="6",
(32,0)*{}="7",
(36,-1)*{\ldots},
(40,0)*{}="8",
(48,0)*{}="9",
(56,0)*{}="10",
(48,-8)*{\bu}="dk",
\ar @{-} "d1";"l" <0pt>
\ar @{-} "d1";"u1" <0pt>
\ar @{-} "d2";"u1" <0pt>
\ar @{-} "d2";"u2" <0pt>
\ar @{-} "d3";"u2" <0pt>
\ar @{-} "l";"top" <0pt>
\ar @{-} "10";"top" <0pt>
\ar @{-} "top";"1" <0pt>
\ar @{-} "u1";"2" <0pt>
\ar @{-} "u2";"5" <0pt>
\ar @{-} "d1";"3" <0pt>
\ar @{-} "d2";"4" <0pt>
\ar @{-} "d3";"6" <0pt>
\ar @{-} "d3";"7" <0pt>
\ar @{-} "dk";"8" <0pt>
\ar @{-} "dk";"9" <0pt>
\ar @{-} "dk";"10" <0pt>
\endxy}\Ea}_{2k\ \text{vertices of type}\ (3,0),\  2k\ \text{vertices of type}\ (1,2)   }
 \in q\LB_{-1,0}(4k,0), \ \ \ \forall\ k\geq 1
$$
are given by non-zero cohomology classes in  $H^{4k}(\cM_{1,4k})$ with $2k$ theta-like subgraphs as in
the ribbon graph above (which corresponds to the case $k=1$). These examples prove Theorem B
in the Introduction.


\bip

\bip

{\large
\section{\bf Totality of cohomology groups of moduli spaces as a complex}
}

\bip

\subsection{Deformation complex of properads under $q\LBcd$} Let $(\cP,\sd)$ be a dg properad.
We represent its elements as decorated $(m,n)$-corollas, see (\ref{2: generic elements of cP as (m,n)-corollas}).
Assume $\cP$ comes equipped with a morphism of properads
$$
\Ba{rccc}
i: & (\HoqLBcd, \delta) & \lon & (\cP,\sd)\vspace{1mm}\\
&
\Ba{c}\resizebox{14mm}{!}{  \xy
(-7,7)*+{_1}="U1";
(-3,7)*+{_2}="U2";
(2,5)*{...};
(7,7)*+{_m}="U3";
(0,0)*{\bu}="C";
(-7,-7)*+{_1}="L1";
(-3,-7)*+{_2}="L2";
(2,-5)*{...};
(7,-7)*+{_n}="L3";
\ar @{-} "C";"L1" <0pt>
\ar @{-} "C";"L2" <0pt>
\ar @{-} "C";"L3" <0pt>
\ar @{-} "C";"U1" <0pt>
\ar @{-} "C";"U2" <0pt>
\ar @{-} "C";"U3" <0pt>
 \endxy}
 \Ea
 &\lon &
 \Ba{c}\resizebox{14mm}{!}{  \xy
(-7,7)*+{_1}="U1";
(-3,7)*+{_2}="U2";
(2,5)*{...};
(7,7)*+{_m}="U3";
(0,0)*{\circledcirc}="C";
(-7,-7)*+{_1}="L1";
(-3,-7)*+{_2}="L2";
(2,-5)*{...};
(7,-7)*+{_n}="L3";
\ar @{-} "C";"L1" <0pt>
\ar @{-} "C";"L2" <0pt>
\ar @{-} "C";"L3" <0pt>
\ar @{-} "C";"U1" <0pt>
\ar @{-} "C";"U2" <0pt>
\ar @{-} "C";"U3" <0pt>
 \endxy}
 \Ea
\Ea
$$
The images of the generators of $\HoqLBcd$ are special elements in $\cP$ and hence are distinguished by a special notation $\circledcirc$ for the vertices.

\sip

According to the general theory developed in \cite{MV}, the deformation complex of the morphism $i$ is given by the following $\Z$-graded vector space,
$$
\Def(\HoqLBcd \stackrel{i}{\rar} \cP)=\prod_{m\geq 1, n\geq 0} \cP(m,n)\otimes_{\bS_m\times \bS_m} \left(\sgn_m^{\ot |c|}\otimes \sgn_n^{\ot |d|}\right)[c(1-m)+d(1-n)]
$$
Its elements can also be represented as  corollas whose vertices are decorated by elements of $\cP$ and whose out- and in-going legs (if any) are (skew)symmetrized so that their labellings can be omitted. To distinguish these elements from elements of $\cP$ and $\HoqLBcd$, we denote vertices of such corollas by $\circledast$, e.g. the corolla
$$
 \Ba{c}\resizebox{14mm}{!}
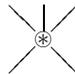
{  \xy
(-6,6)*+{}="U1";
(-0,6)*+{}="U2";
(6,6)*+{}="U3";
(0,0)*{\circledast}="C";
(-6,-6)*+{}="L1";
(6,-6)*+{}="L3";
\ar @{-} "C";"L1" <0pt>
\ar @{-} "C";"L3" <0pt>
\ar @{-} "C";"U1" <0pt>
\ar @{-} "C";"U2" <0pt>
\ar @{-} "C";"U3" <0pt>
 \endxy}
 \Ea
$$
can stand for an element in the deformation complex of the form
$$
\Ba{c}
\resizebox{15mm}{!}{ \xy
(0,0)*{\circ}="d1",
(10,0)*{\circ}="d2",
(-5,-5)*{}="dl",
(5,-5)*{}="dc",
(15,-5)*{}="dr",
(0,10)*{\circ}="u1",
(10,10)*{\circ}="u2",
(5,15)*{}="uc",
(5,15)*{}="uc",
(15,15)*{}="ur",
(0,15)*{}="ul",
\ar @{->} "d1";"d2" <0pt>
\ar @{<-} "d2";"dc" <0pt>
\ar @{<-} "d2";"dr" <0pt>
\ar @{<-} "u1";"d1" <0pt>
\ar @{->} "u1";"u2" <0pt>
\ar @{<-} "u1";"d2" <0pt>
\ar @{->} "u2";"d2" <0pt>
\ar @{<-} "u2";"d1" <0pt>
\ar @{<-} "uc";"u2" <0pt>
\ar @{<-} "ur";"u2" <0pt>
\ar @{<-} "ul";"u1" <0pt>
\endxy}
\Ea \in \Def(\HoqLBcd\rar \cP).
$$
which should be understood as the decoration of the $\circledast$-vertex of the
above corolla.
 Then the differential in the deformation complex $\delta$ is given  by
\Beq\label{3: differential in Def(Hoqlieb to P)}
\delta: \hspace{-3mm}
\Ba{c}\resizebox{14mm}{!}{  \xy
(-7,7)*+{}="U1";
(-3,7)*+{}="U2";
(1.5,5)*{...};
(7,7)*+{}="U3";
%
%
(0,0)*{\circledast}="C";
(-7,-7)*+{}="L1";
(-3,-7)*+{}="L2";
(1.5,-5)*{...};
(7,-7)*+{}="L3";
\ar @{-} "C";"L1" <0pt>
\ar @{-} "C";"L2" <0pt>
\ar @{-} "C";"L3" <0pt>
\ar @{-} "C";"U1" <0pt>
\ar @{-} "C";"U2" <0pt>
\ar @{-} "C";"U3" <0pt>
 \endxy}
 \Ea
 \lon
\sd \hspace{-3mm}
\Ba{c}\resizebox{14mm}{!}{  \xy
(-7,7)*+{}="U1";
(-3,7)*+{}="U2";
(1.5,5)*{...};
(7,7)*+{}="U3";
%
%
(0,0)*{\circledast}="C";
(-7,-7)*+{}="L1";
(-3,-7)*+{}="L2";
(1.5,-5)*{...};
(7,-7)*+{}="L3";
\ar @{-} "C";"L1" <0pt>
\ar @{-} "C";"L2" <0pt>
\ar @{-} "C";"L3" <0pt>
\ar @{-} "C";"U1" <0pt>
\ar @{-} "C";"U2" <0pt>
\ar @{-} "C";"U3" <0pt>
 \endxy}
 \Ea
+
 \sum
%
\left(
\pm
\Ba{c}\resizebox{15mm}{!}{  \xy
(-10,8)*{}="11";
(-7,8)*{}="12";
(0,8)*{}="13";
(-3,-5)*{...};
(-13,-13)*{...};
    (-5,+2)*{\circledcirc}="L";
 (-14,-6)*{\circledast}="B";
 (-20,-13)*{}="b1";
 (-17,-13)*{}="b2";
 (-10,-13)*{}="b3";
 (-20,1)*{}="a1";
 (-17,1)*{}="a2";
 (-11,1)*{}="a3";
  (-8,-5)*{}="C";
   (1,-5)*{}="D";
\ar @{-} "D";"L" <0pt>
\ar @{-} "C";"L" <0pt>
\ar @{-} "B";"L" <0pt>
\ar @{-} "B";"b1" <0pt>
\ar @{-} "B";"b2" <0pt>
\ar @{-} "B";"b3" <0pt>
\ar @{-} "B";"a1" <0pt>
\ar @{-} "B";"a2" <0pt>
\ar @{-} "B";"a3" <0pt>
\ar @{-} "11";"L" <0pt>
\ar @{-} "12";"L" <0pt>
\ar @{-} "13";"L" <0pt>
%
 %
 \endxy}
 \Ea
\mp
\Ba{c}\resizebox{15mm}{!}{  \xy
(-10,8)*{}="11";
(-7,8)*{}="12";
(0,8)*{}="13";
(-3,-5)*{...};
(-13,-13)*{...};
    (-5,+2)*{\circledast}="L";
 (-14,-6)*{\circledcirc}="B";
 (-20,-13)*{}="b1";
 (-17,-13)*{}="b2";
 (-10,-13)*{}="b3";
 (-20,1)*{}="a1";
 (-17,1)*{}="a2";
 (-11,1)*{}="a3";
  (-8,-5)*{}="C";
   (1,-5)*{}="D";
\ar @{-} "D";"L" <0pt>
\ar @{-} "C";"L" <0pt>
\ar @{-} "B";"L" <0pt>
\ar @{-} "B";"b1" <0pt>
\ar @{-} "B";"b2" <0pt>
\ar @{-} "B";"b3" <0pt>
\ar @{-} "B";"a1" <0pt>
\ar @{-} "B";"a2" <0pt>
\ar @{-} "B";"a3" <0pt>
\ar @{-} "11";"L" <0pt>
\ar @{-} "12";"L" <0pt>
\ar @{-} "13";"L" <0pt>
%
 %
 \endxy}
 \Ea
 \right)
\Eeq
where the rule of signs depends on $d$ and is read from (\ref{3: d in qLBcd_infty}).

\sip

In  practice we work in this paper with properads under $q\LBcd$ rather than under $\HoqLBcd$.
The above machinery applies to such cases via the composition
$$
i: \HoqLBcd \twoheadrightarrow \LBcd \stackrel{\imath}{\lon}\cP
$$
i.e.\ again we work, by definition, with the same graded vector space
 $$
\Def(q\LBcd \stackrel{i}{\rar} \cP)=\prod_{m\geq 1, n\geq 0} \cP(m,n)\otimes_{\bS_m\times \bS_m} \left(\sgn_m^{\ot |c|}\otimes \sgn_n^{\ot |d|}\right)[c(1-m)+d(1-n)]
$$
but the differential in (\ref{3: differential in Def(Hoqlieb to P)}) simplifies --- all the terms in the r.h.s.\ vanish except the ones containing images of the generators of
$q\LBcd$,
\Beq\label{3 delta in Def(qLBcd to P)}
\delta: \hspace{-3mm}
\Ba{c}\resizebox{13mm}{!}{  \xy
(-7,7)*+{}="U1";
(-3,7)*+{}="U2";
(1.5,5)*{...};
(7,7)*+{}="U3";
(0,0)*{\circledast}="C";
(-7,-7)*+{}="L1";
(-3,-7)*+{}="L2";
(1.5,-5)*{...};
(7,-7)*+{}="L3";
\ar @{-} "C";"L1" <0pt>
\ar @{-} "C";"L2" <0pt>
\ar @{-} "C";"L3" <0pt>
\ar @{-} "C";"U1" <0pt>
\ar @{-} "C";"U2" <0pt>
\ar @{-} "C";"U3" <0pt>
 \endxy}
 \Ea
\lon
d \hspace{-3mm}
\Ba{c}\resizebox{13mm}{!}{  \xy
(-7,7)*+{}="U1";
(-3,7)*+{}="U2";
(1.5,5)*{...};
(7,7)*+{}="U3";
%
%
(0,0)*{\circledast}="C";
(-7,-7)*+{}="L1";
(-3,-7)*+{}="L2";
(1.5,-5)*{...};
(7,-7)*+{}="L3";
\ar @{-} "C";"L1" <0pt>
\ar @{-} "C";"L2" <0pt>
\ar @{-} "C";"L3" <0pt>
\ar @{-} "C";"U1" <0pt>
\ar @{-} "C";"U2" <0pt>
\ar @{-} "C";"U3" <0pt>
 \endxy}
 \Ea
+
 \sum
\left(
\pm\hspace{-1mm}
\Ba{c}\resizebox{15mm}{!}{  \xy
(-5,8)*{}="12";
(-13,-13)*{...};
    (-5,+2)*{\circledcirc}="L";
 (-14,-6)*{\circledast}="B";
 (-20,-13)*{}="b1";
 (-17,-13)*{}="b2";
 (-10,-13)*{}="b3";
 (-20,1)*{}="a1";
 (-17,1)*{}="a2";
 (-11,1)*{}="a3";
 %
   (1,-5)*{}="D";
\ar @{-} "D";"L" <0pt>
\ar @{-} "B";"L" <0pt>
\ar @{-} "B";"b1" <0pt>
\ar @{-} "B";"b2" <0pt>
\ar @{-} "B";"b3" <0pt>
\ar @{-} "B";"a1" <0pt>
\ar @{-} "B";"a2" <0pt>
\ar @{-} "B";"a3" <0pt>
\ar @{-} "12";"L" <0pt>
 \endxy}
 \Ea
 \pm\hspace{-1mm}
\Ba{c}\resizebox{15mm}{!}{  \xy
(-9,8)*{}="11";
(-1,8)*{}="12";
(-13,-13)*{...};
%
    (-5,+2)*{\circledcirc}="L";
 (-14,-6)*{\circledast}="B";
 (-20,-13)*{}="b1";
 (-17,-13)*{}="b2";
 (-10,-13)*{}="b3";
 (-20,1)*{}="a1";
 (-17,1)*{}="a2";
 (-11,1)*{}="a3";
\ar @{-} "B";"L" <0pt>
\ar @{-} "B";"b1" <0pt>
\ar @{-} "B";"b2" <0pt>
\ar @{-} "B";"b3" <0pt>
\ar @{-} "B";"a1" <0pt>
\ar @{-} "B";"a2" <0pt>
\ar @{-} "B";"a3" <0pt>
\ar @{-} "12";"L" <0pt>
\ar @{-} "11";"L" <0pt>
 \endxy}
 \Ea
\pm\hspace{-1mm}
\Ba{c}\resizebox{15mm}{!}{  \xy
(-10,8)*{}="11";
(-7,8)*{}="12";
(0,8)*{}="13";
(-3,-5)*{...};
%
%
    (-5,+2)*{\circledast}="L";
 (-14,-6)*{\circledcirc}="B";
 (-18,-11)*{}="a1";
 (-10,-11)*{}="a3";
  (-8,-5)*{}="C";
   (1,-5)*{}="D";
\ar @{-} "D";"L" <0pt>
\ar @{-} "C";"L" <0pt>
\ar @{-} "B";"L" <0pt>

\ar @{-} "B";"a1" <0pt>
\ar @{-} "B";"a3" <0pt>
\ar @{-} "11";"L" <0pt>
\ar @{-} "12";"L" <0pt>
\ar @{-} "13";"L" <0pt>
%
 %
 \endxy}
 \Ea
\pm\hspace{-1mm}
\Ba{c}\resizebox{15mm}{!}{  \xy
(-10,8)*{}="11";
(-7,8)*{}="12";
(0,8)*{}="13";
(-3,-5)*{...};
%
%
    (-5,+2)*{\circledast}="L";
 (-14,-6)*{\circledcirc}="B";
 (-20,0)*{}="a1";
 (-14,-10)*{}="a3";
  (-8,-5)*{}="C";
   (1,-5)*{}="D";
\ar @{-} "D";"L" <0pt>
\ar @{-} "C";"L" <0pt>
\ar @{-} "B";"L" <0pt>

\ar @{-} "B";"a1" <0pt>
\ar @{-} "B";"a3" <0pt>
\ar @{-} "11";"L" <0pt>
\ar @{-} "12";"L" <0pt>
\ar @{-} "13";"L" <0pt>
%
 %
 \endxy}
 \Ea
\pm\hspace{-1mm}
\Ba{c}\resizebox{15mm}{!}{  \xy
(-10,8)*{}="11";
(-7,8)*{}="12";
(0,8)*{}="13";
(-3,-5)*{...};
    (-5,+2)*{\circledast}="L";
 (-14,-6)*{\circledcirc}="B";
 (-20,0)*{}="a1";
 (-14,0)*{}="a3";
  (-8,-5)*{}="C";
   (1,-5)*{}="D";
\ar @{-} "D";"L" <0pt>
\ar @{-} "C";"L" <0pt>
\ar @{-} "B";"L" <0pt>

\ar @{-} "B";"a1" <0pt>
\ar @{-} "B";"a3" <0pt>
\ar @{-} "11";"L" <0pt>
\ar @{-} "12";"L" <0pt>
\ar @{-} "13";"L" <0pt>
%
 %
 \endxy}
 \Ea
 \right)
\Eeq
that is, the following corollas
$$
\Ba{c}\resizebox{8mm}{!}{  \xy
(-5,6)*{}="1";
    (-5,+1)*{\circledcirc}="L";
  (-8,-5)*+{_1}="C";
   (-2,-5)*+{_2}="D";
\ar @{-} "D";"L" <0pt>
\ar @{-} "C";"L" <0pt>
\ar @{-} "1";"L" <0pt>
 \endxy}
 \Ea:= \imath\left( \Ba{c}\begin{xy}
 <0mm,0.66mm>*{};<0mm,4mm>*{}**@{-},
 <0.39mm,-0.39mm>*{};<2.2mm,-2.2mm>*{}**@{-},
 <-0.35mm,-0.35mm>*{};<-2.2mm,-2.2mm>*{}**@{-},
 <0mm,0mm>*{\bu};<0mm,0mm>*{}**@{},
   <0.39mm,-0.39mm>*{};<2.9mm,-4mm>*{^{_2}}**@{},
   <-0.35mm,-0.35mm>*{};<-2.8mm,-4mm>*{^{_1}}**@{},
\end{xy}\Ea \right)
, \ \ \
 \Ba{c}\resizebox{8mm}{!}{  \xy
(-5,-5)*{}="1";
    (-5,-1)*{\circledcirc}="L";
  (-8,5)*+{_1}="C";
   (-2,5)*+{_2}="D";
\ar @{-} "D";"L" <0pt>
\ar @{-} "C";"L" <0pt>
\ar @{-} "1";"L" <0pt>
 \endxy}
 \Ea
:=
\imath\left(
\Ba{c}\begin{xy}
 <0mm,-0.55mm>*{};<0mm,-3mm>*{}**@{-},
 <0.5mm,0.5mm>*{};<2.2mm,2.2mm>*{}**@{-},
 <-0.48mm,0.48mm>*{};<-2.2mm,2.2mm>*{}**@{-},
 <0mm,0mm>*{\bu};<0mm,0mm>*{}**@{},
 <0.5mm,0.5mm>*{};<2.7mm,2.8mm>*{^{_2}}**@{},
 <-0.48mm,0.48mm>*{};<-2.7mm,2.8mm>*{^{_1}}**@{},
 \end{xy}\Ea
 \right)
, \ \ \ \
 \Ba{c}\resizebox{8mm}{!}{  \xy
(-5,6)*{^1}="1";
(0,-1)*{\circledcirc}="C";
(0,6)*{^2}="2";
(5,6)*{^3}="3";
\ar @{-} "C";"1" <0pt>
\ar @{-} "C";"2" <0pt>
\ar @{-} "C";"3" <0pt>
 \endxy}
 \Ea
  :=
 \imath\left( \Ba{c}\begin{xy}
 <0mm,-1mm>*{\bu};<-4mm,3mm>*{^{_1}}**@{-},
 <0mm,-1mm>*{\bu};<0mm,3mm>*{^{_2}}**@{-},
 <0mm,-1mm>*{\bu};<4mm,3mm>*{^{_3}}**@{-},
 \end{xy}\Ea
 \right).
$$
Hence the differential in the deformation complex spits into four parts,
\Beq\label{3: d=d+d_1+d_2+d_3 in Def}
\delta=\sd + \delta_1 + \delta_2 + \delta_3
\Eeq
where $\sd$ is the differential in $\cP$, $\delta_1$ (resp., $\delta_2$ and $\delta_3$) is given by all terms involving $\Ba{c}\resizebox{8mm}{!}{  \xy
(-5,6)*{}="1";
    (-5,+1)*{\circledcirc}="L";
  (-8,-5)*+{_1}="C";
   (-2,-5)*+{_2}="D";
\ar @{-} "D";"L" <0pt>
\ar @{-} "C";"L" <0pt>
\ar @{-} "1";"L" <0pt>
 \endxy}
 \Ea$
 (resp.,  $\Ba{c}\resizebox{8mm}{!}{  \xy
(-5,-5)*{}="1";
    (-5,-1)*{\circledcirc}="L";
  (-8,5)*+{_1}="C";
   (-2,5)*+{_2}="D";
\ar @{-} "D";"L" <0pt>
\ar @{-} "C";"L" <0pt>
\ar @{-} "1";"L" <0pt>
 \endxy}
 \Ea$ and  $\Ba{c}\resizebox{8mm}{!}{  \xy
(-5,6)*{^1}="1";
(0,-1)*{\circledcirc}="C";
(0,6)*{^2}="2";
(5,6)*{^3}="3";
\ar @{-} "C";"1" <0pt>
\ar @{-} "C";"2" <0pt>
\ar @{-} "C";"3" <0pt>
 \endxy}
 \Ea$). It is easy to see that they satisfy the relations
$$
\delta_1^2=0, \ \ \delta_3^2=0,\ \   \delta_2^2+\delta_1\delta_3+ \delta_3\delta_1=0, \ \  \delta_2\delta_3+ \delta_3\delta_2=0,\ \ \delta_1\delta_2 + \delta_2\delta_1=0.
$$
All these relations are just incarnations of the generating relations in the ideal $\cR_q$ used to define $q\LBcd$.

\subsection{An application to the gravity properad}\label{4: subseq on totaligy of Mg,n+m as a complex}

The above Proposition and Theorem {\ref{3: Main theorem on H(twRgra)}} imply that the totality of compactly supported cohomology groups of moduli spaces $\cM_{g,m+n}$ with (skew)symmetrized punctures can be made into a complex.
 Indeed, the deformation complex of the above morphism $j$ has the following structure
\Beq\label{4: totality of Mg,n+m as Def}
\Def(q\LB_{-1,0} \stackrel{j}{\rar} \GRav)\simeq \prod_{g,n\geq 0,m\geq 1\atop 2g+n+m\geq 3} H_c^{\bu-1}(\cM_{g,m+n})\ot_{\bS_m^{op}\times \bS_n} (\sgn_m\ot \id_n)
\Eeq
The induced differential has three parts (cf.\ (\ref{3: d=d+d_1+d_2+d_3 in Def})),
$$
\delta=\delta' + \delta'' + \delta''',
$$
$$
\delta':H_c^{\bu-1}(\cM_{g,m+n})\ot_{\bS_m^{op}\times \bS_n} (\sgn_m\ot \id_n) \lon
H_c^{\bu}(\cM_{g,m+(1+n)})\ot_{\bS_{m}^{op}\times \bS_{n+1}} (\sgn_{m}\ot \id_{n+1})
$$
$$
\delta'':H_c^{\bu-1}(\cM_{g,m+n})\ot_{\bS_m^{op}\times \bS_n} (\sgn_m\ot \id_n) \lon
H_c^{\bu}(\cM_{g,(m+1)+n)})\ot_{\bS_{m+1}^{op}\times \bS_{n}} (\sgn_{m+1}\ot \id_{n})
$$
$$
\delta''':H_c^{\bu-1}(\cM_{g,m+n})\ot_{\bS_m^{op}\times \bS_n} (\sgn_m\ot \id_n) \lon
H_c^{\bu}(\cM_{g,(m+2)+(n-1)})\ot_{\bS_{m+2}^{op}\times \bS_{n-1}} (\sgn_{m+2}\ot \id_{n-1})
$$
with each part  having a geometric interpretation as the attachment of the pair of pants
but in different hyperbolic incarnations --- the first part  employs the pair of pants as the hyperbolic sphere with one geodesic boundary and two cusps,  the second one as the hyperbolic sphere with two geodesic boundaries and one cusp, and the third one as $\bS^2$ with three geodesic boundaries and no cusps. The differential preserves the genus parameter so this complex is a product
of complexes, one for each fixed $g$.

\subsubsection{\bf Theorem C from the Introduction} {\it The complex (\ref{4: totality of Mg,n+m as Def}) has infinitely many cohomology classes.}

\mip

These cohomology classes can described explicitly in terms of twisted ribbon graphs (see below). There are strong indications that they all come from $H_c^\bu(\cM_1)$ so we conjecture that {\it the cohomology of the complex (\ref{4: totality of Mg,n+m as Def}) is equal to $\prod_{g\geq 1} H_c^\bu(\cM_g)$},
where $\cM_g$ is the moduli space of (unpunctured) algebraic curves of genus $g$.

\sip

To proof this Theorem we have to understand better the deformation complex
of the dg properad $\HoqLBcd$ and its relation to the M.\ Kontsevich graph complexes and their oriented versions.

\subsection{Reminder on Maxim Kontsevich's  graph complexes (after \cite{Ko2,W})}\label{4: subseq on reminder about GC_d}
An (ordinary, i.e.\ non-ribbon) {\em graph}\, $\Ga$ can be understood as a 1-dimensional $CW$ complex whose 0-cells are called {\em vertices}\, and 1-cells are called {\em edges}; the set of vertices of $\Ga$ is denoted by $V(\Ga)$ and the set of edges by $E(\Ga)$. A graph $\Ga$ is called {\em directed}\, if each edge $e\in E(\Ga)$ comes equipped with a fixed orientation; we show it in pictures as an arrow. If a vertex $v$ of a directed graph $\Ga$ has
$m\geq 0$ outgoing edges and $n\geq 0$ incoming edges, then we say that $v$ is an $(m,n)$-{\em vertex}.
A $(1,1)$-vertex is called {\em passing}. Here are a few examples of directed graphs
$$
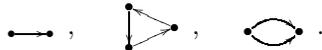
\Ba{c}\resizebox{7mm}{!}{
\xy
 (0,0)*{\bullet}="a",
(8,0)*{\bu}="b",
\ar @{->} "a";"b" <0pt>
\endxy}\Ea, \ \ \
\Ba{c}\resizebox{8mm}{!}{
\xy
 {\ar@{->}(-5,4)*{\bu};(-5,-4)*{\bu}};
   {\ar@{<-}(-5,4)*{\bu};(4,0)*{\bu}};
 {\ar@{<-}(4,0)*{\bu};(-5,-4)*{\bu}};
\endxy}\Ea, \ \ \
\Ba{c}\resizebox{8mm}{!}{\xy
   {\ar@/^0.6pc/(-5,0)*{\bu};(5,0)*{\bu}};
 {\ar@/_0.6pc/(-5,0)*{\bu};(5,0)*{\bu}};
\endxy}\Ea .
$$


\sip

Let $DG_{n,l}$ be the set of directed  graphs $\Ga$ with $n$ vertices and $l$  edges such that
some bijections $V(\Ga)\rar [n]$ and $E(\Ga)\rar [l]$ are fixed, i.e.\ every edge and every vertex of $\Ga$ is marked. The permutation group $\bS_l$ acts on $DG_{n,l}$ by relabeling the edges. Given any integer $d\in \Z$, we consider a collection of $\bS_n$-modules,
$$
\caD\cG ra_{d}=\left\{\caD\cG ra_d(n):= \prod_{l\geq 0} \K \langle DG_{n,l}\rangle \ot_{ \bS_l}  \sgn_l^{\ot |d-1|} [l(d-1)]   \right\}_{n\geq 1}
$$
where the group $\bS_n$ acts on $\caD\cG ra_d(n)$ by relabeling the vertices. Each generator  $\Ga$ of this $\Z$-graded vector space is assigned the degree $(1-d)\# E(\Ga)$ and, if $d$ is even, it is also assumed that $\Ga$  is equipped with a choice of ordering of edges up to an even permutation (an odd permutation acts as the multiplication by $-1$). The latter condition kills some generators from $DG_{n,l}$; for example, the graph $\Ba{c}\resizebox{7mm}{!}{\xy
(-5,2)*{^1},
(5,2)*{^2},
   {\ar@/^0.6pc/(-5,0)*{\bu};(5,0)*{\bu}};
 {\ar@/^{-0.6pc}/(-5,0)*{\bu};(5,0)*{\bu}};
\endxy}
\Ea\in \caD\cG ra_{d\in 2\Z}(2)$ vanishes identically as it admits an automorphism which changed the ordering of edges by an odd permutation, i.e.\ it is equal to minus itself.

\sip

The $\bS$-module $\caD\cG ra_d$ has an operad structure with the operadic composition,
$$
\Ba{rccc}
\circ_i: &  \caD\cG ra_d(n) \times \caD\cG ra_d(m) &\lon & \caD\cG ra_d(m+n-1),  \ \ \forall\ i\in [n]\\
         &       (\Ga_1, \Ga_2) &\lon &      \Ga_1\circ_i \Ga_2,
\Ea
$$
given the substitution of the graph $\Ga_2$ into the $i$-labeled vertex of $\Ga_1$ and taking a sum over re-attachments of dangling edges (attached earlier to $v_i$) to the vertices of $\Ga_2$ in all possible ways. The ordering of edges $\Ga_1\circ_i \Ga_2$ is uniquely determined by the given of orderings of edges in $\Ga_1$ followed by the given ordering of edges in $\Ga_2$. For example,
$$
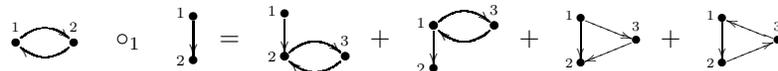
\Ba{c}\resizebox{10mm}{!}{
\xy
(-5,2)*{^1},
(5,2)*{^2},
   {\ar@/^0.6pc/(-5,0)*{\bu};(5,0)*{\bu}};
 {\ar@/^0.6pc/(5,0)*{\bu};(-5,0)*{\bu}};
\endxy}
\Ea
 \ \ \circ_1\
 \Ba{c}\resizebox{4mm}{!}{
\xy
(-2,7)*{^1},
(-2,0)*{_2},
 {\ar@{->}(0,7)*{\bu};(0,0)*{\bu}};
\endxy}\Ea
=
 \Ba{c}\resizebox{12mm}{!}{
\xy
(-7,7)*{^1},
(-7,0)*{_2},
(5,2)*{^3},
 {\ar@{->}(-5,7)*{\bu};(-5,0)*{\bu}};
   {\ar@/^0.6pc/(-5,0)*{\bu};(5,0)*{\bu}};
 {\ar@/^0.6pc/(5,0)*{\bu};(-5,0)*{\bu}};
\endxy}\Ea
+
 \Ba{c}\resizebox{12mm}{!}{
\xy
(-7,0)*{^1},
(-7,-7)*{_2},
(5,2)*{^3},
 {\ar@{->}(-5,0)*{\bu};(-5,-7)*{\bu}};
   {\ar@/^0.6pc/(-5,0)*{\bu};(5,0)*{\bu}};
 {\ar@/^0.6pc/(5,0)*{\bu};(-5,0)*{\bu}};
\endxy}\Ea
+
 \Ba{c}\resizebox{11mm}{!}{
\xy
(-7,4)*{^1},
(-7,-4)*{_2},
(5,2)*{^3},
 {\ar@{->}(-5,4)*{\bu};(-5,-4)*{\bu}};
   {\ar@{->}(-5,4)*{\bu};(5,0)*{\bu}};
 {\ar@{->}(5,0)*{\bu};(-5,-4)*{\bu}};
\endxy}\Ea
+
 \Ba{c}\resizebox{11mm}{!}{
\xy
(-7,4)*{^1},
(-7,-4)*{_2},
(5,2)*{^3},
 {\ar@{->}(-5,4)*{\bu};(-5,-4)*{\bu}};
   {\ar@{<-}(-5,4)*{\bu};(5,0)*{\bu}};
 {\ar@{<-}(5,0)*{\bu};(-5,-4)*{\bu}};
\endxy}
\Ea
$$

\sip

One often works with  an ``undirected" version, $\cG ra_d$, of the above operad  by taking coinvariants with respect of the action of
the group $\bS_l \ltimes
 (\bS_2)^l$ on the linear space $\K\langle  DG_{n,l} \rangle$ given by relabelling the
edges and reversing their  directions,
$$
\cG ra_d:=\left\{
\cG ra_d(n):= \prod_{l\geq 0} \K \langle G_{n,l}\rangle \ot_{ \bS_l \ltimes
 (\bS_2)^l}  \sgn_l^{\ot|d-1|}\ot \sgn_2^{\ot l|d|} [l(d-1)]\right\}_{n\geq 1}
$$
For $d$ even elements
$\cG ra_d(n)$ can be understood as genuine {undirected} graphs, e.g.
 $$
\Ba{c}\resizebox{8mm}{!}{
\xy
(0,2)*{^1},
(8,2)*{^2},
 (0,0)*{\bullet}="a",
(8,0)*{\bu}="b",
\ar @{-} "a";"b" <0pt>
\endxy}\Ea \in \cG ra_{d\in 2\Z}(2),
 \ \ \
\Ba{c}\resizebox{9mm}{!}{
\xy
(-7,4)*{^1},
(-7,-4)*{_2},
(4,2)*{^3},
 {\ar@{-}(-5,4)*{\bu};(-5,-4)*{\bu}};
   {\ar@{-}(-5,4)*{\bu};(4,0)*{\bu}};
 {\ar@{-}(4,0)*{\bu};(-5,-4)*{\bu}};
\endxy}\Ea\in \cG ra_{d\in 2\Z}(3)
$$
 while for $d$ odd the directions on edges are defined up to a flip (and multiplication by $-1$),
$$
\Ba{c}\resizebox{8mm}{!}{
\xy
(0,2)*{^1},
(8,2)*{^2},
 (0,0)*{\bullet}="a",
(8,0)*{\bu}="b",
\ar @{->} "a";"b" <0pt>
\endxy}\Ea=
-
\Ba{c}\resizebox{8mm}{!}{
\xy
(0,2)*{^2},
(8,2)*{^1},
 (0,0)*{\bullet}="a",
(8,0)*{\bu}="b",
\ar @{->} "a";"b" <0pt>
\endxy}\Ea
,
\ \ \
\Ba{c}\resizebox{9mm}{!}{\xy
(-5,2)*{^1},
(5,2)*{^2},
   {\ar@/^0.6pc/(-5,0)*{\bu};(5,0)*{\bu}};
 {\ar@/^{0.6pc}/(5,0)*{\bu};(-5,0)*{\bu}};
\endxy}
\Ea
=-
\Ba{c}\resizebox{9mm}{!}{\xy
(-5,2)*{^1},
(5,2)*{^2},
   {\ar@/^0.6pc/(-5,0)*{\bu};(5,0)*{\bu}};
 {\ar@/^{-0.6pc}/(-5,0)*{\bu};(5,0)*{\bu}};
\endxy}
\Ea\ \ \text{in}\ \ \cG ra_{d\in 2\Z+1}(2).
$$

There is an obvious epimorphism of operads,
$$
\caD\cG ra_d \lon \cG ra_d.
$$
and a less obvious (but easy to check) morphism of operads
$$
\Ba{rccc}
i: & \caL ie_d & \lon & \caD\cG ra_d\vspace{1mm} \\
&
\Ba{c}
\xy
 <0mm,0.55mm>*{};<0mm,3.5mm>*{}**@{-},
 <0.5mm,-0.5mm>*{};<2.2mm,-2.2mm>*{}**@{-},
 <-0.48mm,-0.48mm>*{};<-2.2mm,-2.2mm>*{}**@{-},
 <0mm,0mm>*{\bu};<0mm,0mm>*{}**@{},
 <0.5mm,-0.5mm>*{};<2.7mm,-3.2mm>*{_2}**@{},
 <-0.48mm,-0.48mm>*{};<-2.7mm,-3.2mm>*{_1}**@{},
 \endxy\Ea
& \lon &
\frac{1}{2}\left(
\Ba{c}\resizebox{7mm}{!}{\xy
(0,2)*{_{1}},
(7,2)*{_{2}},
 (0,0)*{\bu}="a",
(7,0)*{\bu}="b",
\ar @{->} "a";"b" <0pt>
\endxy}\Ea
+ (-1)^d \
\Ba{c}\resizebox{7mm}{!}{\xy
(0,2)*{_{2}},
(7,2)*{_{1}},
 (0,0)*{\bu}="a",
(7,0)*{\bu}="b",
\ar @{->} "a";"b" <0pt>
\endxy}\Ea
\right)
\Ea
$$
Hence both operads of graphs $\caD\cG ra_d$  and  $\cG ra_d$ can be twisted with the help of T.\ Willwacher's endofunctor producing new important dg operads \cite{Ko3,W}. Another very useful corollary of this fact is that one can introduce dg Lie algebras, or {\it graph complexes},
$$
\dfGC:=\Def(\Lie_d \rar \caD\cG ra_d)\ \        \text{and} \ \    \fGC_d:= \Def(\Lie_d \rar \cG ra_d)
$$
controlling deformations of the above map $i$ which are given explicitly by
$$
\mathsf{dfGC}_d= \prod_{l\geq 0}\prod_{n\geq 1} \K \langle DG_{n,l}\rangle \ot_{\bS_n\times \bS_l} \left(\sgn_n^{\ot |d|} \ot \sgn_l^{|d-1|}\right) [d(1-n) + l(d-1)]
$$
 $$
\mathsf{fGC}_d= \prod_{l\geq 0}\prod_{n\geq 1} \K \langle DG_{n,l}\rangle \ot_{\bS_n\times  (\bS_l \ltimes
 (\bS_2)^l)} \left(\sgn_n^{\ot |d|} \ot (\sgn_l^{|d-1|}\ot (\sgn_2^{|d|})^{\ot l})\right) [d(1-n) + l(d-1)]
$$
the letter {\it f}\, standing for {\it full}. Thus numerical labels of vertices of graphs get (skew)symmerized and hence can be forgotten.  Graphs $\Ga$ from $\mathsf{dfGC}_{d}$ or $\fGC_d$ are assigned the  cohomological degree,
$$
|\Ga|=d(\# V(\Ga) -1) + (1-d) \# E(\Ga).
$$
The differential is determined by the image of the Lie generator in the corresponding operad of graphs, and is given in, say, $\fGC_d$ explicitly by
$$
\delta \Ga:= \sum_{v\in V(\Ga)} \left(\delta_v'\Ga -(-1)^{|\Ga|} \delta_v''\Ga\right)
$$
where $\delta_v'$ substitutes into the vertex $v$ the graph $\Ba{c}\resizebox{7mm}{!}{
\xy
 (0,1)*{\bullet}="a",
(8,1)*{\bu}="b",
\ar @{-} "a";"b" <0pt>
\endxy}\Ea$ and redistributes the attached to $v$ half-edges (whose set is denoted by $H(v)$) among the two new vertices in all possible ways,
$$
\delta'_v  \Ba{c}\resizebox{10mm}{!}{\xy
 (0,0)*{\bullet}="a",
(-8,0)*{}="1",
(8,0)*{}="2",
(-4,6)*{}="3",
(-4,-6)*{}="4",
(4,6)*{}="5",
(4,-6)*{}="6",
\ar @{-} "a";"1" <0pt>
\ar @{-} "a";"2" <0pt>
\ar @{-} "a";"3" <0pt>
\ar @{-} "a";"4" <0pt>
\ar @{-} "a";"5" <0pt>
\ar @{-} "a";"6" <0pt>
\endxy}
\Ea=\sum_{H(v)=I'\sqcup I''\atop
\# I',\#I''\geq 0}
\Ba{c}\resizebox{11mm}{!}{\xy
%
(0,-10)*{\underbrace{\ \ \ \ \  \ \ \ \ \ }_{I'\ \text{half-edges}}},
(0,18)*{\overbrace{ \ \ \  \ \ \ \ \ }^{I''\ \text{half-edges}}},
 (0,0)*{\bullet}="a",
 (0,8)*{\bullet}="b",
(-8,0)*{}="1",
(-4,14)*{}="2",
(8,0)*{}="3",
(-4,-6)*{}="4",
(4,14)*{}="5",
(4,-6)*{}="6",
\ar @{-} "a";"b" <0pt>
\ar @{-} "a";"1" <0pt>
\ar @{-} "b";"2" <0pt>
\ar @{-} "a";"3" <0pt>
\ar @{-} "a";"4" <0pt>
\ar @{-} "b";"5" <0pt>
\ar @{-} "a";"6" <0pt>
\endxy}
\Ea
,
$$
while $\delta_v''$ attaches to $v$ a new univalent vertex
$$
\delta'_v  \Ba{c}\resizebox{10mm}{!}{\xy
%
 (0,0)*{\bullet}="a",
(-8,0)*{}="1",
(8,0)*{}="2",
(-4,6)*{}="3",
(-4,-6)*{}="4",
(4,6)*{}="5",
(4,-6)*{}="6",
\ar @{-} "a";"1" <0pt>
\ar @{-} "a";"2" <0pt>
\ar @{-} "a";"3" <0pt>
\ar @{-} "a";"4" <0pt>
\ar @{-} "a";"5" <0pt>
\ar @{-} "a";"6" <0pt>
\endxy}
\Ea=
\Ba{c}\resizebox{10mm}{!}{\xy
%
(7,5)*{\bullet}="b",
 (0,0)*{\bullet}="a",
(-8,0)*{}="1",
(7.8,-1)*{}="2",
(-5,5.6)*{}="3",
(-4,-6)*{}="4",
(1.5,7.1)*{}="5",
(4,-6)*{}="6",
\ar @{-} "a";"1" <0pt>
\ar @{-} "a";"2" <0pt>
\ar @{-} "a";"3" <0pt>
\ar @{-} "a";"4" <0pt>
\ar @{-} "a";"5" <0pt>
\ar @{-} "a";"6" <0pt>
\ar @{-} "a";"b" <0pt>
\endxy}\Ea
$$
If  a vertex $v\in V(\Ga)$ has at least one half-edge, then the term with the new univalent vertex in the sum $\delta_v'\Ga$ cancels out with the similar term
in $\delta_v''\Ga$.
The graph complexes $\dfGC_d$ and, resp., $\fGC_d$ contain subcomplexes $\dGC_d$ and $\GC_d^{\geq 2}$ spanned by {\it connected}\, graphs with all vertices at least bivalent and with no passing vertices; the inclusions
$$
\dGC_d \lon \dfGC_d, \ \ \ \GC_d^{\geq 2}\lon \fGC_d
$$
are quasi-isomorphisms \cite{W}. There is also a quasi-isomorphism of dg Lie algebras
\Beq\label{3: from GC to dGC}
\mathsf{GC}^{\geq 2}_d\lon \mathsf{dGC}_d,
\Eeq
which sends a graph with no direction on edges into a sum of graphs with all possible directions on edges \cite{W}. The complex $\dGC_d$ contains a subcomplex $\OGC_d$ spanned by directed graphs with no closed paths of directed edges.
It was proven
 in \cite{W2} that
 $$
 H^\bu(\GC_{d}^{\geq 2}) =H^\bu(\GC_{d+1}^{or})
 $$
 and that for $d=2$
$$
H^0(\dGC_2)=H^0(\GC_2^{\geq 2})=H^0(\GC_3^{or})=\grt_1,
$$
where
$\grt_1$ is the Lie algebra of the Grothendieck-Teichm\"uller group $GRT_1$ introduced in \cite{Dr3}.

\sip

The complex $\GC^{\geq 2}_d$ decomposes into a direct sum of complexes,
$$
\GC_d^{\geq 2}=\GC_d \oplus \GC_d^2,
$$
where $\GC_d$ is generated by graphs with each vertex at least trivalent, and $\GC_d^2$ is generated by graphs containing at least one bivalent vertex. The cohomology of the latter has been computed in \cite{W},
$$
H^\bu(\GC_d^2) = \bigoplus_{j\geq 1\atop j\equiv 2d+1 \mod 4} \K[d-j],
$$
where the summand $\K[d-j]$ is generated by the polytope with $j$ bivalent vertices and $j$ edges. For example, for $d$ odd, the graph,
$$
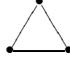
\Ba{c}\resizebox{9mm}{!}{\xy
 (0,0)*{\bullet}="a",
(12,0)*{\bu}="b",
(6,10)*{\bu}="c",
\ar @{-} "a";"b" <0pt>
\ar @{-} "a";"c" <0pt>
\ar @{-} "c";"b" <0pt>
\endxy}\Ea
$$
contributes into that cohomology, while for $d$ even the graph

$$
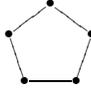
\Ba{c}\resizebox{13mm}{!}{
\xy
(0,8)*{\bullet}="1",
(-8,2)*{\bullet}="5",
(8,2)*{\bullet}="2",
(-5,-7)*{\bullet}="4",
(5,-7)*{\bullet}="3",
\ar @{-} "1";"2" <0pt>
\ar @{-} "2";"3" <0pt>
\ar @{-} "3";"4" <0pt>
\ar @{-} "4";"5" <0pt>
\ar @{-} "5";"1" <0pt>
\endxy}\Ea
$$
is a  cohomology class.

\sip

The complexes $\GC_d$ for $d$ of the same parity are isomorphic to each other (up to degree shift) so that there are  essentially two different complexes, one for $d\in 2\Z$ and one for $d\in 2\Z+1$.
The even complex $\GC_2$ is very useful because its degree zero cohomology is equal to $\grt_1$.
The odd complex $\GC_3$ has also  very nice properties:  its cohomology is concentrated in degrees $\leq -3$; its degree $-3$ cohomology is generated by trivalent graphs and has the structure of a graded commutative algebra which is conjecturally generated by formal variables  $t$, $\om_0$, $\om_1$, $\ldots$ subject to the relations $\om_p\om_q=\om_0\om_{p+q}$ \cite{Kn,V,KWZ}.

\sip

Unfortunately, there is no direct relation between the complexes $\GC_d^{\geq 2}$ and $\OGC_{d+1}$ (however there is a beautiful explicit quasi-isomorphism between their dual subcomplexes \cite{Z} spanned by  graphs containing at least one trivalent vertex). It is not hard too see, however, what happens to the polytope classes in $\GC_d$ in their oriented incarnation in $\OGC_{d+1}$: the polytope cohomology class with $j$ vertices in $\GC_d^{\geq 2}$ gets associated a polytope cohomology class with $j+1$ non-passing vertices in $\OGC_{d+1}$  \cite{MW},
$$
H^\bu(\GC_{odd}^{\geq 2})\, \ni\Ba{c}\resizebox{11mm}{!}{ \begin{tikzpicture}[baseline=-.65ex]
  \node[int] (v0) at (0:1) {};
\node[int] (v1) at (120:1) {};
\node[int] (v2) at (-120:1) {};
\draw (v0) edge (v1) edge (v2) (v2) edge (v1);
 \end{tikzpicture}
 }\Ea
\quad \leftrightsquigarrow   \quad
\Ba{c}\resizebox{15mm}{!}{ \begin{tikzpicture}[baseline=-.65ex]
  \node[int] (v0) at (0:1) {};
\node[int] (v1) at (90:1) {};
\node[int] (v2) at (180:1) {};
\node[int] (v3) at (270:1) {};
\draw[-latex] (v0) edge (v1) edge (v3) (v2) edge (v1) edge (v3);
 \end{tikzpicture}
 }\Ea\in H^\bu(\OGC_{even}),
$$
$$
H^\bu(\GC_{even}^{\geq 2})\, \ni\Ba{c}\resizebox{14mm}{!}{
\begin{tikzpicture}[baseline=-.65ex]
  \node[int] (v0) at (0:1) {};
\node[int] (v1) at (72:1) {};
\node[int] (v2) at (144:1) {};
\node[int] (v3) at (216:1) {};
\node[int] (v4) at (288:1) {};
\draw (v0) edge (v1) edge (v4) (v2) edge (v3) edge (v1) (v3) edge (v4);
 \end{tikzpicture}
 }\Ea
\quad
\leftrightsquigarrow   \quad
\Ba{c}\resizebox{16mm}{!}{
 \begin{tikzpicture}[baseline=-.65ex]
  \node[int] (v0) at (0:1) {};
\node[int] (v1) at (60:1) {};
\node[int] (v2) at (120:1) {};
\node[int] (v3) at (180:1) {};
\node[int] (v4) at (240:1) {};
\node[int] (v5) at (300:1) {};
\draw[-latex] (v0) edge (v1) edge (v5) (v2) edge (v1) edge (v3) (v4) edge (v3) edge (v5);
 \end{tikzpicture}
 }\Ea\in H^\bu(\OGC_{odd}).
 $$
We shall use this correspondence below.

\subsection{Graph complexes and deformations of quasi-Lie bialgebras}\label{4: subseq on Def(HoqLB) and GC}
Let $I_s$ be the differential ideal in $\HoqLBcd$ generated by $(m,0)$-corollas with $m\geq 3$ (we call such corollas {\it sources}). The associated short exact sequence of dg properads,
$$
0 \lon I_s \lon \HoqLBcd \lon \HoLBcd \lon 0
$$
implies an associated short exact sequence of dg Lie algebras,
$$
0 \lon \Ker \pi \lon \Def(\HoqLBcd \stackrel{\Id}{\rar} \HoqLBcd) \stackrel{\pi}{\lon} \Def(\HoqLBcd \stackrel{\Id}{\rar} \HoLBcd)\lon 0
$$
where $\Ker \pi$ is the subcomplex of $\Def(\HoqLBcd \rar \HoqLBcd)$ spanned by graphs with at least one source. As
$$
\Def(\HoqLBcd \stackrel{\Id}{\rar} \HoLBcd)\simeq \Def(\HoLBcd \stackrel{\Id}{\rar} \HoLBcd)
$$
we have a canonical morphism of cohomology groups
$$
\pi: H^\bu\left( \Def(\HoqLBcd \stackrel{\Id}{\rar} \HoqLBcd)\right) \lon
 H^\bu\left( \Def(\HoLBcd \stackrel{\Id}{\rar} \HoLBcd)\right).
$$

\sip

Let $\wHoqLBcd$ (resp., $\wHoLBcd$) be the genus completed version of the properad $\HoqLBcd$ (resp., $\wHoLBcd$) and let  $\Der(\wHoqLBcd)$ (resp.\ $\Der(\wHoqLBcd)$) be the dg Lie algebra of its derivations. There is a canonical isomorphism of  complexes (but not of Lie algebras),
$$
\Der(\wHoqLBcd) \simeq \Def(\HoqLBcd \stackrel{\Id}{\rar} \HoLBcd)[1]
$$
and similarly for $\Der(\wHoLBcd)$.
One of the main results in \cite{MW2} is a construction of an explicit quasi-isomorphism (up to one rescaling class) of dg Lie algebras
\Beq\label{3: Morhism F from OGC to Der(HoLBcd)}
\Ba{rccc}
 F\colon & \mathsf{OGC}_{c+d+1} &\to & \Der(\wHoLBcd)\\
         &   \Ga & \to & F(\Ga)
         \Ea
\Eeq
where the derivation $F(\Ga)$ has, by definition, the following values
on the generators of the completed properad  $\wHoLBcd$
\Beq \label{5:derivation Fstar(Ga)}
\left(\Ba{c}\resizebox{12mm}{!}{\begin{xy}
 <0mm,0mm>*{\circ};<0mm,0mm>*{}**@{},
 <-0.6mm,0.44mm>*{};<-8mm,5mm>*{}**@{-},
 <-0.4mm,0.7mm>*{};<-4.5mm,5mm>*{}**@{-},
 <0mm,0mm>*{};<-1mm,5mm>*{\ldots}**@{},
 <0.4mm,0.7mm>*{};<4.5mm,5mm>*{}**@{-},
 <0.6mm,0.44mm>*{};<8mm,5mm>*{}**@{-},
   <0mm,0mm>*{};<-8.5mm,5.5mm>*{^1}**@{},
   <0mm,0mm>*{};<-5mm,5.5mm>*{^2}**@{},
   <0mm,0mm>*{};<9.0mm,5.5mm>*{^m}**@{},
 <-0.6mm,-0.44mm>*{};<-8mm,-5mm>*{}**@{-},
 <-0.4mm,-0.7mm>*{};<-4.5mm,-5mm>*{}**@{-},
 <0mm,0mm>*{};<-1mm,-5mm>*{\ldots}**@{},
 <0.4mm,-0.7mm>*{};<4.5mm,-5mm>*{}**@{-},
 <0.6mm,-0.44mm>*{};<8mm,-5mm>*{}**@{-},
   <0mm,0mm>*{};<-8.5mm,-6.9mm>*{^1}**@{},
   <0mm,0mm>*{};<-5mm,-6.9mm>*{^2}**@{},
   <0mm,0mm>*{};<9.0mm,-6.9mm>*{^n}**@{},
 \end{xy}}\Ea\right)\cdot F(\Ga)
=
 \sum_{s:[n]\rar V(\Ga)\atop \hat{s}:[m]\rar V(\Ga)}  \Ba{c}\resizebox{11mm}{!}  {\xy
 (-6,7)*{^1},
(-3,7)*{^2},
(2.5,7)*{},
(7,7)*{^m},
(-3,-8)*{_2},
(3,-6)*{},
(7,-8)*{_n},
(-6,-8)*{_1},
(0,4.5)*+{...},
(0,-4.5)*+{...},
(0,0)*+{\Ga}="o",
(-6,6)*{}="1",
(-3,6)*{}="2",
(3,6)*{}="3",
(6,6)*{}="4",
(-3,-6)*{}="5",
(3,-6)*{}="6",
(6,-6)*{}="7",
(-6,-6)*{}="8",
\ar @{-} "o";"1" <0pt>
\ar @{-} "o";"2" <0pt>
\ar @{-} "o";"3" <0pt>
\ar @{-} "o";"4" <0pt>
\ar @{-} "o";"5" <0pt>
\ar @{-} "o";"6" <0pt>
\ar @{-} "o";"7" <0pt>
\ar @{-} "o";"8" <0pt>
\endxy}\Ea\ \ \ \ \ \ \forall\ m,n\geq 1, m+n\geq 3,
\Eeq
 where the sum is taken over all ways of attaching the incoming and outgoing legs to the graph $\Ga$, and setting to zero every graph which has at east one vertex of valency $\leq 2$ or at least one vertex of type $(0,n)$ or $(n,0)$ with $n\geq 3$.

\sip

It is easy to see that the above map $F$ extends to a morphism of dg Lie algebras,

$$
F_q: \OGC_{c+d+1} \lon  \Der(\wHoqLBcd)
$$
which is given by the same formula as in (\ref{3: Morhism F from OGC to Der(HoLBcd)}) as above except that now $m\geq 1, n\geq 0, m+n\geq 3$ and we set to zero on the r.h.s.\ every graph which has at east one vertex of valency $\leq 2$ or at least one vertex of type $(n,0)$ with $n\geq 3$. Therefore we get a commutative diagram,

$$
\Ba{c}\resizebox{90mm}{!}{
\xymatrix{
 H^\bu(\OGC_{c+d+1}) \ar[dr]_-{F} \ar[r]^-{F_q} & H^{\bu+1}\left( \Def(\HoqLBcd \stackrel{\Id}{\rar} \HoqLBcd)\right)\ar[d]^-{\pi}
 \\
& H^{\bu+1}\left( \Def(\HoLBcd \stackrel{\Id}{\rar} \HoLBcd)\right)
}}\Ea
$$
where the diagonal map is an isomorphism (up to one rescaling class). Therefore there is an {\it injection}\, of cohomology groups,
$$
H^\bu(\GC_{c+d}^{\geq 2}) \simeq H^\bu(\OGC_{c+d+1}) \stackrel{F_q}{\lon} H^{\bu+1}\left( \Def(\HoqLBcd \stackrel{\Id}{\rar} \HoqLBcd)\right).
$$
Given a dg properad $\cP$ under $\HoqLBcd$,
$$
i: \HoqLBcd \lon \cP,
$$
the map $i$ induces a morphism of complexes
$$
\Def(\HoqLBcd \stackrel{\Id}{\rar} \HoqLBcd) \lon \Def(\HoqLBcd \stackrel{i}{\rar} \cP)
$$
and hence a morphism of cohomology groups,
$$
H^\bu(\OGC_{c+d+1}) \stackrel{F_q}{\lon} H^{\bu+1}\left( \Def(\HoqLBcd \stackrel{i}{\rar} \cP)\right).
$$
which sometimes can be checked by explicit computations.

\subsection{Proof of Theorem C from the Introduction} Let us apply the above machinery to the morphism
$$
j: q\LB_{-1,0} \lon H^\bu(\GRav)
$$
and check what happens to the polytope cohomology classes in $H^\bu(\OGC_{0})$ under the canonical morphism
$$
H^\bu(\OGC_0) \lon H^{\bu+1}\left( \Def(q\LB_{-1,0} \stackrel{j}{\rar} \GRav)\right).
$$
The simplest polytope class is given by the graph
$$
\Ba{c}\resizebox{14mm}{!}{\xy
 (0,0)*{\bullet}="a1",
(16,0)*{\bu}="a2",
(8,8)*{\bu}="u",
(8,-8)*{\bu}="d",
\ar @{->} "a1";"u" <0pt>
\ar @{->} "a1";"d" <0pt>
\ar @{->} "a2";"u" <0pt>
\ar @{->} "a2";"d" <0pt>
\endxy}\Ea
\in H^4(\OGC_{0}).
$$
The map $F_q$ sends it into the following degree $5$ cohomology class in $\Def(q\LB_{-1,0}\stackrel{\Id}{\rar} q\LB_{-1,0})$,
$$
\Ba{c}\resizebox{22mm}{!}{\xy
 (-8,0)*{\bullet}="a1",
(8,0)*{\bu}="a2",
(0,8)*{\bu}="u",
(0,-8)*{\bu}="d",
 (-14,0)*{}="1",
(14,0)*{}="2",
(0,14)*{}="3",
(0,-14)*{}="4",
\ar @{->} "a1";"u" <0pt>
\ar @{->} "a1";"d" <0pt>
\ar @{->} "a2";"u" <0pt>
\ar @{->} "a2";"d" <0pt>
\ar @{->} "1";"a1" <0pt>
\ar @{->} "2";"a2" <0pt>
\ar @{->} "u";"3" <0pt>
\ar @{->} "d";"4" <0pt>
\endxy}\Ea
+ 2
\Ba{c}\resizebox{19mm}{!}{\xy
 (-8,0)*{\bullet}="a1",
(8,0)*{\bu}="a2",
(0,8)*{\bu}="u",
(0,-8)*{\bu}="d",
 (-2,0)*{}="1",
(14,0)*{}="2",
(0,14)*{}="3",
(0,-14)*{}="4",
\ar @{->} "a1";"u" <0pt>
\ar @{->} "a1";"d" <0pt>
\ar @{->} "a2";"u" <0pt>
\ar @{->} "a2";"d" <0pt>
\ar @{<-} "1";"a1" <0pt>
\ar @{->} "2";"a2" <0pt>
\ar @{->} "u";"3" <0pt>
\ar @{->} "d";"4" <0pt>
\endxy}\Ea
+
\Ba{c}\resizebox{15mm}{!}{\xy
 (-8,0)*{\bullet}="a1",
(8,0)*{\bu}="a2",
(0,8)*{\bu}="u",
(0,-8)*{\bu}="d",
 (-2,0)*{}="1",
(2,0)*{}="2",
(0,14)*{}="3",
(0,-14)*{}="4",
\ar @{->} "a1";"u" <0pt>
\ar @{->} "a1";"d" <0pt>
\ar @{->} "a2";"u" <0pt>
\ar @{->} "a2";"d" <0pt>
\ar @{<-} "1";"a1" <0pt>
\ar @{<-} "2";"a2" <0pt>
\ar @{->} "u";"3" <0pt>
\ar @{->} "d";"4" <0pt>
\endxy}\Ea
$$
which under the explicit  map $j$ given by (\ref{3: map j from qLB to H(TwGra)}) goes into a genus 1 cycle in the complex
 (\ref{4: totality of Mg,n+m as Def}) represented as a linear combination of twisted graphs from $\tw\RGra_0^{\geq 3}$. That linear combination contains   a summand with no white vertices at all which, in this particular case, comes from the last term in the above sum, and is given by
 $$
j\left(\Ba{c}\resizebox{15mm}{!}{\xy
 (-8,0)*{\bullet}="a1",
(8,0)*{\bu}="a2",
(0,8)*{\bu}="u",
(0,-8)*{\bu}="d",
 (-2,0)*{}="1",
(2,0)*{}="2",
(0,14)*{}="3",
(0,-14)*{}="4",
\ar @{->} "a1";"u" <0pt>
\ar @{->} "a1";"d" <0pt>
\ar @{->} "a2";"u" <0pt>
\ar @{->} "a2";"d" <0pt>
\ar @{<-} "1";"a1" <0pt>
\ar @{<-} "2";"a2" <0pt>
\ar @{->} "u";"3" <0pt>
\ar @{->} "d";"4" <0pt>
\endxy}\Ea\right) =
 \begin{tikzpicture}[baseline=-.65ex]
  \node[int] (v0) at (-45:1) {};
  \node[int] (v1) at (45:1) {};
\node[int] (v2) at (135:1) {};
\node[int] (v3) at (225:1) {};
\draw (v0) edge (v1) edge[bend left] (v1) edge[bend right] (v1);
\draw (v2) edge (v3) edge[bend left] (v3) edge[bend right] (v3);
\draw (v1) to[ out=180, in=-80, looseness=2] (v2);
\draw (v0) to[ out=100, in=0, looseness=2] (v3);
 \end{tikzpicture}
 $$
From this explicit representation of that summand we conclude that the cocycle can not be a coboundary; a similar observation holds true for any other polytope class from $\GC_{-1}^2$ --- the summands with no white vertices are all of the similar type and hence can not coboundaries with respect to $\delta_3$ (cf.\ with \S {\ref{3: subsec on non-triv of prop compositions}} above \S 7.4 in \cite{MW}). This completes the proof of Theorem C.

\mip

All these polytope cohomology classes in the complex (\ref{4: totality of Mg,n+m as Def}) have genus 1 and come conjecturally from the cusp series of cohomology classes\footnote{I am grateful to Alexey Kalugin for the reference to the papers \cite{B,T} which describe completely $H_c^\bu(\cM_1)$.} in $H^\bu_c(\cM_1)$ (cf. \cite{B,T}). Thus one can make two conjectures about the complex (\ref{4: totality of Mg,n+m as Def}):  the cohomology of that complex
 is equal to (a)
$\prod_{g\geq 1}H_c^\bu(\cM_g)$, and (b) it contains the cohomology of M.\ Kontsevich's {\it odd}\, graph complex $H^\bu(\GC_{-1}^{\geq 2})$.



\def\cprime{$'$}


\begin{thebibliography}{10}


\bibitem[AWZ]{AWZ} A.\ Andersson, T.\ Willwacher and M.\ \v Zivkovi\' c,
{\it Oriented hairy graphs and moduli spaces of curves},  arXiv:2005.00439 (2020).


\bibitem[Be]{B} K.\ A.\ Behrend, {\it Derived $l$-adic categories for algebraic stacks}\,
Memoirs of the AMS,  {\bf 163} (2003), No. 774.


\bibitem[CFL]{CFL} K. Cieliebak, K. Fukaya and J. Latschev, {\em Homological algebra related to surfaces with boundary}, Quantum Topology, {\bf 11}, No.4 (2020) 691-837.


\bibitem[CGP]{CGP}
M.\ Chan, S.\ Galatius and S.\ Payne,
\newblock {\it Tropical curves, graph homology, and top weight cohomology of $M_g$},
\newblock  J.\ Amer.\ Math.\
Soc. {\bf 34}, No.2  (2021) 565-594.


\bibitem[Co1]{Co1} K.\ Costello, {\em The $A_\infty$ operad and the moduli space of curves},   arXiv:
math.AG/0402015 (2004)

\bibitem[Co2]{Co2} K.\ Costello, {\em
A dual version of the ribbon graph decomposition
of moduli space}, Geometry \& Topology {\bf 11} (2007) 1637-1652.

\bibitem[Dr1]{Dr1}
V. Drinfeld, {\em Quantum groups}, Proc.\ ICM-86 (Berkley) {\bf 1} (1987),  798-820.


\bibitem[Dr2]{Dr2}
V. Drinfeld, {\em Quasi-Hopf algebras}, Leningrad Math. J. {\bf 1} (1990),  1419-1457.

\bibitem[Dr3]{Dr3}
V. Drinfeld, {\em On quasitriangular quasi-Hopf algebras and a group closely connected
with $Gal(\bar{Q}/Q)$}, Leningrad Math. J. {\bf 2}, No.\ 4 (1991),  829-860.



\bibitem[Ge]{Ge} E.\ Getzler, {\em Two-dimensional topological gravity and equivariant cohomology}, Comm.\ Math.\ Phys. {\bf 163} (1994), no. 3,
473-489.



\bibitem[KMS]{KMS}
A.\ Khoroshkin, N.\ Markarian and S.\ Shadrin, {\it Hypercommutative operad as a homotopy quotient of $BV$}, Commun.\ Math.\ Phys.\ {\bf 322} (2013) 697-729.

\bibitem[KWZ]{KWZ}  A.\ Khoroshkin, T.\ Willwacher and
M.\ \v Zivkovi\'c, {\em Differentials on graph complexes},  Adv. Math. {\bf 307} (2017), 1184-1214.



\bibitem[Kn]{Kn} J.\  Kneissler, {\it The number of primitive Vassiliev invariants up to degree 12}, arXiv:q-alg/9706022, (1997).

 \bibitem[Ko1]{Ko1} M.\ Kontsevich, {\em  Formal (non)commutative symplectic geometry}. In Proceedings of the I. M. Gelfand seminar
1990-1992, pages 173-188. Birkhauser, 1993.

\bibitem[Ko2]{Ko2} M.\ Kontsevich, {\em Formality Conjecture}. In: D. Sternheimer et al. (eds.),
Deformation Theory and Symplectic
Geometry, Kluwer 1997, 139-156.


\bibitem[Ko3]{Ko3}  M.\ Kontsevich,  {\em Deformation quantization
 of Poisson manifolds}, Lett.\ Math.\ Phys. {\bf 66} (2003) 157-216.


\bibitem[Me]{Me} S.A.\ Merkulov,  {\em From gravity to string topology}, \ { Lett.\ Math.\ Phys.}, {\bf 113} (2023) 1-23



 \bibitem[MeVa]{MV}  S.A.\ Merkulov and  B.\ Vallette,
{\em Deformation theory of representations of prop(erad)s I \& II},
 J.\ f\"ur die reine und angewandte Mathematik (Qrelle)  {\bf 634}, 51-106,
 \& {\bf 636}, 123-174 (2009)

 \bibitem[MW1]{MW} S.A.\ Merkulov and T.\ Willwacher, {\em Props of ribbon graphs, involutive Lie bialgebras and moduli spaces of curves}, preprint  arXiv:1511.07808 (2015) 51pp.

 \bibitem[MW2]{MW2} S.\ Merkulov and T.\ Willwacher, {\em Deformation theory of  Lie bialgebra properads},   In: Geometry and Physics: A Festschrift in honour of Nigel Hitchin, Oxford University Press 2018, pp. 219-248.


\bibitem[Pe]{Pe}  R.C.\ Penner, {\em The decorated Teichm\"uller space of punctured surfaces}, Comm.\ Math.\
Phys. {\bf 113} (1987), 299-339.


\bibitem[Ta]{T}
L.\ Taelman, {\it Characteristic classes for curves of genus one}, Michigan Math. J. {\bf 64(3)} (2015)
 633-654.


\bibitem[TZ]{TZ} T.\ Tradler and M.\ Zeinalian, {\em  Algebraic string operations},
K-Theory, {\bf 38}(1) (2007) 59-82.

\bibitem[Vo]{V}
P.\ Vogel, {\it Algebraic structures on modules of diagrams}, J.\ Pure Appl.\ Algebra, {\bf 215} (2011) 1292-1339.


\bibitem[WW]{WW} N.\ Wahl and C.\ Westerland, {\em Hochschild homology of structured
algebras},  Adv. Math., 288 (2016) 240-307.

\bibitem[War]{Wa} B.\ Ward, {\em Maurer-Cartan elements and cyclic operads}, J.\ Noncommut.\ Geom.\ {\bf 10} (2016), no. 4, 1403-1464.



 \bibitem[Wil1]{W} T.\ Willwacher, {\em M.\ Kontsevich's graph complex and the Grothendieck-Teichmueller Lie algebra},
Invent. Math. {\bf 200} (2015), 671-760.

\bibitem[Wil2]{W2} T.\ Willwacher, {\em Oriented graph complexes},   Comm. Math. Phys. {\bf 334} (2015), no. 3, 1649-1666.

\bibitem[Zi]{Z} M.\ \v Zivkovi\' c, {\em Multi-oriented graph complexes and quasi-isomorphisms between them I: oriented graphs}, High.\ Struct.\ 4 (2020), no. 1, 266-283.

 \end{thebibliography}
\end{document}